\documentclass[a4paper]{amsart}

\usepackage{amscd,amsthm,amsfonts,amssymb,amsmath}
\usepackage[colorlinks=true]{hyperref}
\usepackage{fullpage}
\usepackage[all]{xy}
\usepackage{mathrsfs}
\usepackage{bm}
\usepackage{tikz}
\usepackage{bbm}
\usepackage{graphics}
\usepackage{enumerate}
\usepackage[utf8]{inputenc}

\newcommand{\BC}{{\mathbb {C}}}
\newcommand{\BG}{{\mathbb {G}}}

\newcommand{\BN}{{\mathbb {N}}}\newcommand{\BP}{{\mathbb {P}}}
\newcommand{\BR}{{\mathbb {R}}}

\newcommand{\BZ}{{\mathbb {Z}}}

\newcommand{\CC}{{\mathcal {C}}}

\newcommand{\CK}{{\mathcal {K}}}
\newcommand{\CN}{{\mathcal {N}}}\newcommand{\CO}{{\mathcal {O}}}

\newcommand{\msd}{\mathscr{D}}
\newcommand{\mse}{\mathscr{E}}\newcommand{\msf}{\mathscr{F}}
\newcommand{\msl}{\mathscr{L}}
\newcommand{\msp}{\mathscr{P}}
\newcommand{\mst}{\mathscr{T}}
\newcommand{\msv}{\mathscr{V}}
\newcommand{\msy}{\mathscr{Y}}

\newcommand{\fa}{{\mathfrak{a}}} \newcommand{\fc}{{\mathfrak{c}}} 
 \newcommand{\fg}{{\mathfrak{g}}} \newcommand{\fh}{{\mathfrak{h}}}
  \newcommand{\fl}{{\mathfrak{l}}}
\newcommand{\fm}{{\mathfrak{m}}} \newcommand{\fn}{{\mathfrak{n}}} \newcommand{\fp}{{\mathfrak{p}}}
 \newcommand{\fs}{{\mathfrak{s}}} \newcommand{\ft}{{\mathfrak{t}}}
 \newcommand{\fv}{{\mathfrak{v}}} 
 \newcommand{\fz}{{\mathfrak{z}}}

  \newcommand{\fD}{{\mathfrak{D}}}

                      \newcommand{\Ad}{{\mathrm{Ad}}}              \newcommand{\ad}{{\mathrm{ad}}}

                      \newcommand{\bs}{\backslash}                 
				\newcommand{\bdd}{\mathrm{bdd}}

\newcommand{\Cl}{{\mathrm{Cl}}}                      		\newcommand{\Cent}{{\mathrm{Cent}}}

                  			\newcommand{\der}{{\mathrm{der}}}

\newcommand{\etale}{\'{e}tale~}			\newcommand{\el}{\mathrm{ell}}

\newcommand{\Hom}{{\mathrm{Hom}}}

\renewcommand{\Im}{{\mathrm{Im}}}                    \newcommand{\Ind}{{\mathrm{Ind}}}

\newcommand{\Lie}{{\mathrm{Lie}}}

                      \renewcommand{\mod}{\, \mathrm{mod}\, }

			\newcommand{\Nrd}{{\mathrm{Nrd}}}		\newcommand{\Norm}{{\mathrm{Norm}}}

\newcommand{\ov}{\overline}

\newcommand{\rank}{{\mathrm{rank}}}                  \newcommand{\ra}{\rightarrow}

                   \newcommand{\reg}{{\mathrm{reg}}}               \newcommand{\Res}{{\mathrm{Res}}}
				 \newcommand{\rs}{{\mathrm{rs}}}

                                  \newcommand{\Supp}{{\mathrm{Supp}}}
 \newcommand{\sk}{\medskip}                      \newcommand{\s}{\sk\noindent}

                   \newcommand{\tr}{{\mathrm{Tr}}}                 
\newcommand{\Trd}{{\mathrm{Trd}}}

\newcommand{\vol}{{\mathrm{vol}}}

\newcommand{\wt}{\widetilde}

\newtheorem{thm}{Theorem}[section]
\newtheorem{coro}[thm]{Corollary}
\newtheorem{lem}[thm]{Lemma}
\newtheorem{prop}[thm]{Proposition}

\theoremstyle{definition}

\theoremstyle{remark}
\newtheorem{remark}[thm]{Remark}

\numberwithin{equation}{subsection}
\def\mat(#1,#2,#3,#4){
  \begin{pmatrix}
  #1 & #2 \\ #3 & #4
  \end{pmatrix}
}

\begin{document}
\title{A local trace formula for $p$-adic infinitesimal symmetric spaces: the case of Guo-Jacquet}
\author{Huajie Li}
\date{\today}
\maketitle

\begin{abstract}

We establish an invariant local trace formula for the tangent space of some symmetric spaces over a non-archimedean local field of characteristic zero. These symmetric spaces are studied in Guo-Jacquet trace formulae and our methods are inspired by works of Waldspurger and Arthur. Some other results are given during the proof including a noninvariant local trace formula, Howe's finiteness for weighted orbital integrals and the representability of the Fourier transform of weighted orbital integrals. These local results are prepared for the comparison of regular semi-simple terms, which are weighted orbital integrals, of an infinitesimal variant of Guo-Jacquet trace formulae. 

\end{abstract}

\tableofcontents


\section{\textbf{Introduction}}

The Guo-Jacquet trace formula \cite{MR1382478} is a promising tool to generalise Waldspurger's result \cite{MR783511} on the relation between toric periods and central values of automorphic $L$-functions for $GL_2$ to higher ranks. It is inspired by Jacquet's new proof \cite{MR868299} of Waldspurger's theorem. Although such a formula has not been established in full generality, its simple form was used by Feigon-Martin-Whitehouse \cite{MR3805647} to obtain some evidence for the conjecture of Guo-Jacquet. For applications, one needs to compare geometric sides of Guo-Jacquet trace formulae for different symmetric pairs. Some local results on the comparison of relative orbital integrals include Guo's fundamental lemma \cite{MR1382478} and Zhang's smooth transfer \cite{MR3414387}. 

In order to study the Guo-Jacquet trace formula and its comparison, one may begin with an infinitesimal variant. That is to say, we replace a symmetric space by its tangent space (called an infinitesimal symmetric space). Such a variant should share some similarities with the geometric side of Guo-Jacquet trace formula. It is simpler than the original formula because spectral objects are replaced by the Fourier transform of geometric objects (cf. \cite{MR1344131} and \cite{MR1893921}). Moreover, by the method of descent dating back to Harish-Chandra's works, the comparison at the infinitesimal level should imply the comparison of geometric sides of original formulae (see \cite{MR3414387} on the transfer of orbital integrals). 

An infinitesimal variant of Guo-Jacquet trace formulae has been established in \cite{zbMATH07499568} and \cite{MR4350885} via an analogue of Arthur's truncation process in \cite{MR518111} (see also \cite{MR1893921} for its Lie algebra variant). We actually consider more general cases suggested by \cite{MR3299843} and \cite{MR2806111}. Most (namely regular semi-simple) terms appearing in these formulae can be written as explicit weighted orbital integrals on infinitesimal symmetric spaces over a number field (see \cite[Theorem 9.2]{zbMATH07499568} and \cite[Theorem 9.2]{MR4350885}). They are noninvariant analogues of ordinary orbital integrals (which can be compared locally thanks to \cite{MR1382478} and \cite{MR3414387}) and should be the next objects to be compared. As the first evidence, the weighted fundamental lemma has been proved in \cite[Theorem 10.9]{MR4350885} thanks to Labesse's work \cite{MR1339717} on the base change for $GL_n$. 

The same philosophy of Waldspurger's work \cite{MR1440722} on the endoscopic transfer has been followed by Zhang \cite{MR3414387} to prove the transfer of local orbital integrals on infinitesimal symmetric spaces of Guo-Jacquet. A simple form of the local trace formula \cite[Lemma 6.5]{MR3414387}, Howe's finiteness for orbital integrals \cite[Theorem 6.1]{MR1375304} and representability of the Fourier transform of orbital integrals \cite[Theorem 6.1]{MR3414387} apart from a simple form of the global trace formula \cite[Theorem 8.4 and p. 1875]{MR3414387} and the fundamental lemma \cite[Lemma 5.18]{MR3414387} at the infinitesimal level have been used in Zhang's proof. It is expected that such a strategy should be extended to the weighted context. In fact, some successful attempts have been made by Chaudouard in \cite{MR2164623} and \cite{MR2332352} on the stable base change. We would like to follow these ideas in the comparison of local weighted orbital integrals on infinitesimal symmetric spaces of Guo-Jacquet. However, further study in noninvariant local harmonic analysis on infinitesimal symmetric spaces is needed to achieve our goal. This paper aims to prepare some essential ingredients such as a noninvariant local trace formula, Howe's finiteness for weighted orbital integrals, representability of the Fourier transform of weighted orbital integrals and an invariant local trace formula. Our methods are mainly inspired by the works of Waldspurger's \cite{MR1344131} and Arthur's \cite{MR1114210}. 

Let $E/F$ be a quadratic field extension of non-archimedean local fields of characteristic zero. Denote by $D$ a central division algebra over $F$ and by $GL_{n,D}$ the reductive group over $F$ whose $F$-points are $GL_n(D)$. We study two generalised cases of Guo-Jacquet trace formulae. The first case is $(G, H)$, where $G:=GL_{2n, D}$ and $H:=GL_{n,D}\times GL_{n,D}$ denotes its diagonal subgroup by diagonal embedding. Denote by $\fs$ the tangent space of the symmetric space $G/H$ at the neutral element, on which $H$ acts by conjugation. The second case is $(G', H')$, where $G'$ is the group of invertible elements in a central simple algebra $\fg'$ over $F$ containing $E$, and $H'$ is the centraliser of $E^\times$ in $G'$. Denote by $\fs'$ the corresponding infinitesimal symmetric space. Notice that $(G,H)$ and $(G',H')$ are the same symmetric pair after a base change to an algebraic closure of $F$ containing $E$. In the rest of the introduction and this paper, we shall focus on results in the first case and provide complete proofs. The second case is similar in statements and proofs, so we shall only state main results, point out additional ingredients and sketch necessary steps for later use. 

Before explaining the main results of this paper, we introduce some necessary notations. Let $\eta$ be the quadratic character of $F^\times/NE^\times$ attached to $E/F$, where $N$ is the norm map $E^\times\ra F^\times$. For $x\in G(F)$, we denote by $\Nrd(x)$ its reduced norm. Let $\CC_c^\infty(\fs(F))$ be the space of locally constant, compactly supported, complex-valued functions on $\fs(F)$. Denote by $\fs_\rs\subseteq\fs$ the Zariski open subset consisting of regular semi-simple elements (see Section \ref{generalcases}). Let $M_0$ be the group of diagonal elements in $G$, which is a common minimal Levi subgroup of $G$ and $H$. A Levi subgroup of $G$ containing $M_0$ is said to be $\omega$-stable if $\omega:=\mat(0,1_n,1_n,0)\in M$. Denote by $\msl^{G,\omega}(M_0)$ the set of $\omega$-stable Levi subgroups of $G$. Suppose that $M\in\msl^{G,\omega}(M_0)$. Write $\fm:=\Lie(M)$ and $M_H:=M\cap H$. Let $A_M$ be the maximal $F$-split torus in the centre of $M$. Denote by $\msf^G(M)$ the set of parabolic subgroups of $G$ containing $M$. This paper is organised in the following way. 

In Section \ref{secnot&pre}, we fix some notations of local harmonic analysis and recall some facts of Arthur's $(G,M)$-families, most of which can be found in \cite[\S I-II]{MR1344131}. In particular, for $M\in\msl^{G,\omega}(M_0)$ and $Q\in\msf^G(M)$, we define the local weight function $v_M^Q$ on $G(F)$ of main interest by \eqref{localwtfn}. 

In Section \ref{secsympai}, we prepare some properties of infinitesimal symmetric spaces. Some of them are stated for a general symmetric pair and most of them are relative avatars of classical works of Harish-Chandra \cite{MR0414797}. Preliminaries on symmetric pairs can be found in \cite{MR1375304} and \cite{MR2553879}. 

In Section \ref{sectwoi}, we define local weighted orbital integrals for the action of $H(F)$ on $\fs(F)$ and study their properties. Let $M\in\msl^{G,\omega}(M_0)$ and $Q\in\msf^G(M)$. For all $f\in\CC_c^\infty(\fs(F))$ and $X\in(\fm\cap\fs_\rs)(F)$, we define 
$$ J_M^Q(\eta, X, f):=|D^\fs(X)|_F^{1/2} \int_{H_X(F)\bs H(F)} f(\Ad(x^{-1})(X)) \eta(\Nrd(x)) v_M^Q(x) dx,   $$
where $|D^\fs(X)|_F$ is the Weyl discriminant factor (see \eqref{weyldisc}) and $H_X$ denotes the centraliser of $X$ in $H$. The distributions $J_M^Q(\eta, X, \cdot)$ on $\fs(F)$ are local analogues of the global weighted orbital integrals obtained in \cite[Theorem 9.2]{zbMATH07499568}.  

In Section \ref{secnontf}, we establish the noninvariant local trace formula. For $f\in\CC_c^\infty(\fs(F))$, we define its Fourier transform $\hat{f}$ by \eqref{fouriertransform1}. For $M\in\msl^{G,\omega}(M_0)$, let $\mst_\el(\fm\cap\fs)$ be a set of representatives for $M_H(F)$-conjugacy classes of $M$-elliptic Cartan subspaces in $\fm\cap\fs$ (see Section \ref{sectsympar1}). 
For $f, f'\in\CC_c^\infty(\fs(F))$, we define 
\[\begin{split}
 J^G(\eta, f, f'):=&\sum_{M\in\msl^{G,\omega}(M_0)} |W_0^{M_n}| |W_0^{GL_n}|^{-1} (-1)^{\dim(A_M/A_G)} \sum_{\fc\in\mst_\el(\fm\cap\fs)} |W(M_H, \fc)|^{-1} \int_{(\fc\cap\fs_\rs)(F)} \\ 
&J_M^G(\eta, X, f, f') dX, 
\end{split}\]
where $W_0^{M_n}$, $W_0^{GL_n}$ and $W(M_H, \fc)$ are certain Weyl groups (see Section \ref{sectsympar1}), and $J_M^G(\eta, X, f, f')$ is defined by \eqref{defV.1.1}. 

\begin{thm}[see Theorems \ref{noninvltf1} and \ref{noninvltf2}]
For all $f, f'\in\CC_c^\infty(\fs(F))$, we have the equality
$$ J^G(\eta, f, \hat{f'})=J^G(\eta, \hat{f}, f'). $$
\end{thm}

This formula results from the Plancherel formula and an analogue of Arthur's truncation process in \cite{MR1114210}. We cannot deduce it via the exponential map as in \cite[\S V]{MR1344131} for lack of a local trace formula for symmetric spaces. One needs to return to the proof of \cite{MR1114210} instead. 

In Section \ref{sechowfin}, we show Howe's finiteness for weighted orbital integrals on $\fs(F)$. For $r\subseteq\fs(F)$ an open compact subgroup, denote by $\CC_c^\infty(\fs(F)/r)$ the subspace of $\CC_c^\infty(\fs(F))$ consisting of the functions invariant by translation of $r$. 

\begin{prop}[see Corollaries \ref{corhowe1} and \ref{corhowe2}]
Let $r$ be an open compact subgroup of $\fs(F)$, $M\in\msl^{G,\omega}(M_0)$ and $\sigma\subseteq (\fm\cap\fs_\rs)(F)$. Suppose that there exists a compact subset $\sigma_0\subseteq(\fm\cap\fs)(F)$ such that $\sigma\subseteq \Ad((M_H)(F))(\sigma_0)$. Then there exists a finite subset $\{X_i: i\in I\}\subseteq\sigma$ and a finite subset $\{f_i: i\in I\}\subseteq\CC_c^\infty(\fs(F)/r)$ such that for all $X\in\sigma$ and all $f\in\CC_c^\infty(\fs(F)/r)$, we have the equality
$$ J_M^G(\eta, X, f)=\sum_{i\in I} J_M^G(\eta, X_i, f) J_M^G(\eta, X, f_i). $$
\end{prop}

The proof originates from Howe's seminal work \cite{MR342645} which is extended to weighted orbital integrals on Lie algebras by Waldspurger. We modify the argument in \cite[\S IV]{MR1344131} to make it apply to our case. 

In Section \ref{secrep}, we show that the Fourier transform of weighted orbital integrals on $\fs(F)$ is represented by locally integrable functions on $\fs(F)$. 

\begin{prop}[see Propositions \ref{repr1} and \ref{repr2}]
Let $M\in\msl^{G,\omega}(M_0)$ and $X\in(\fm\cap\fs_\rs)(F)$. Then there exists a locally constant function $\hat{j}_M^G(\eta, X,\cdot)$ on $\fs_\rs(F)$ such that for all $f\in\CC_c^\infty(\fs(F))$, we have 
$$ J_M^G(\eta, X, \hat{f})=\int_{\fs(F)} f(Y)\hat{j}_M^G(\eta, X,Y)|D^\fs(Y)|_F^{-1/2}dY. $$
\end{prop}

Its proof is similar to that in \cite[\S V]{MR1344131} and makes use of the noninvariant trace formula and Howe's finiteness for weighted orbital integrals. 

In Section \ref{secinvwoi}, we modify weighted orbital integrals to obtain $\eta(\Nrd(\cdot))$-invariant distributions $I_M^G(\eta, X, \cdot)$ on $\fs(F)$ by \eqref{invftofwoi1} and \eqref{invwoi1}, where $M\in\msl^{G,\omega}(M_0)$ and $X\in(\fm\cap\fs_\rs)(F)$. The method is close to Arthur's standard one, but it is simpler here since there is no spectral object involved, which is also a feature of \cite{MR1344131}. By the above proposition and the construction via induction, there also exists a locally constant function $\hat{i}_M^G(\eta, X,\cdot)$ on $\fs_\rs(F)$ such that for all $f\in\CC_c^\infty(\fs(F))$, we have 
$$ I_M^G(\eta, X, \hat{f})=\int_{\fs(F)} f(Y)\hat{i}_M^G(\eta, X,Y)|D^\fs(Y)|_F^{-1/2}dY. $$

In Section \ref{secinvtf}, we establish the invariant local trace formula. For $f, f'\in\CC_c^\infty(\fs(F))$, we define 
\[\begin{split}
 I^G(\eta, f, f'):=&\sum_{M\in\msl^{G, \omega}(M_0)} |W_0^{M_n}| |W_0^{GL_n}|^{-1} (-1)^{\dim(A_M/A_G)} \sum_{\fc\in\mst_\el(\fm\cap\fs)} |W(M_H, \fc)|^{-1} \int_{(\fc\cap\fs_\rs)(F)} \\
 &I_M^G(\eta, X, \hat{f}) I_G^G(\eta, X, f') dX. 
\end{split}\]

\begin{thm}[see Theorems \ref{invltf1} and \ref{invltf2}]
For all $f, f'\in\CC_c^\infty(\fs(F))$, we have the equality
$$ I^G(\eta, f, f')=I^G(\eta, f', f). $$
\end{thm}

This formula is deduced from the noninvariant local trace formula. We mainly consult \cite[\S VII]{MR1344131} for the proof. 

In Section \ref{secvanpro}, we prove a vanishing property at ``infinity'' of the function $\hat{i}_M^G(\eta, \cdot, Y)$ on $(\fm\cap\fs_\rs)(F)$, where $M\in\msl^{G,\omega}(M_0), M\neq G$ and $Y\in\fs_\rs(F)$. Denote by $v_F(\cdot)$ the valuation on $F$. 

\begin{prop}[see Propositions \ref{limitformula1} and \ref{limitformula2}]
Let $M\in\msl^{G,\omega}(M_0), M\neq G$. Let $X\in(\fm\cap\fs_\rs)(F)$ and $Y\in\fs_\rs(F)$. Then there exists $N\in\BN$ such that if $\lambda\in F^\times$ satisfies $v_F(\lambda)<-N$, we have
$$ \hat{i}_M^G(\eta, \lambda X, Y)=0. $$
\end{prop}

It is an analogue of \cite[Proposition 2.2]{MR2164623} and serves as a complement of the limit formula for the function $\hat{i}_G^G(\eta, \cdot, Y)$ on $\fs_\rs(F)$ in \cite[Proposition 7.1]{MR3414387}, where $Y\in\fs_\rs(F)$. In a subsequent paper, we shall prove some relations between the functions $\hat{i}_M^G$ and their analogues for the case of $(G',H')$, which are related to the comparison of Guo-Jacquet trace formulae at the infinitesimal level. The above vanishing property will be used to construct some nice test functions whose Fourier transforms have vanishing nontrivial weighted orbital integrals.  

In the end, we remark that although we concentrate on the case of Guo-Jacquet here, many results in this paper might be extended to other symmetric pairs, which can be seen from their proofs. 

\s{\textbf{Acknowledgement. }}This is part of my thesis under the supervision of Pierre-Henri Chaudouard at the Universit\'{e} de Paris. I am grateful to my thesis advisor for his constant support on this project and particularly for answering my numerous questions about Waldspurger's and Arthur's papers on the local trace formula. I would like to thank Jean-Loup Waldspurger for his valuable comments on my manuscript which pointed out an error in an earlier version. Part of this article was revised when I was a postdoc at the Aix-Marseille Université. I also thank Raphaël Beuzart-Plessis for helpful discussions on general symmetric pairs during the revision. This work was supported by grants from R\'{e}gion Ile-de-France. The project leading to this publication has received funding from Excellence Initiative of Aix–Marseille University–A*MIDEX, a French ``Investissements d’Aveni'' programme. 


\section{\textbf{Notation and preliminaries}}\label{secnot&pre}

\subsection{Fields}
Let $F$ be a non-archimedean local field of characteristic zero. Denote by $|\cdot|_F$ (resp. $v_F(\cdot)$) the normalised absolute value (resp. the valuation) on $F$ and by $\CO_F$ the ring of integers of $F$. Fix a uniformiser $\varpi$ of $\CO_F$. Let $q$ be the cardinality of the residue field of $\CO_F$. 

\subsection{Groups and the map $H_P$}\label{gpandmap}
Let $G$ be a (connected) reductive group defined over $F$. All algebraic groups and algebraic varieties in this article are assumed to be defined over $F$ until further notice. Fix a Levi subgroup $M_0$ of a minimal parabolic subgroup of $G$. 

Denote by $A_G$ the maximal $F$-split torus in the centre of $G$. Define
$$ \fa_G:=\Hom_\BZ(X(G)_F, \BR), $$
where $X(G)_F$ is the group of $F$-rational characters of $G$. Define the homomorphism $H_G: G(F)\ra\fa_G$ by
$$ \langle H_G(x), \chi \rangle=\log(|\chi(x)|_F) $$
for all $x\in G(F)$ and $\chi\in X(G)_F$. Set $\fa_{G,F}:=H_G(G(F))$, which is a lattice in $\fa_G$. 

Fix a maximal compact subgroup $K=K_G$ of $G(F)$ which is admissible relative to $M_0$ in the sense of \cite[p. 9]{MR625344}; in other words, $K$ is the stabiliser in $G(F)$ of a special point in the apartment associated to $A_{M_0}$ of the Bruhat-Tits building of $G$. For a central division algebra $D$ over $F$, we denote by $GL_{n,D}$ the reductive group over $F$ whose $F$-points are $GL_n(D)$. In this paper, when $G=\Res_{E/F} GL_{n,D'}$ with $D'$ being a central division algebra over a finite extension $E$ of $F$, we choose $M_0\simeq(\Res_{E/F} \BG_{m,D'})^n$ to be the subgroup of diagonal elements in $G$, and $K=GL_n(\CO_{D'})$ with $\CO_{D'}$ being the ring of integers of $D'$ (see \cite[p. 191]{MR1344916} for example). Set $W_0^G:=\Norm_{G(F)}(M_0)/M_0(F)$ to be the Weyl group of $(G, M_0)$, where $\Norm_{G(F)}(M_0)$ denotes the normaliser of $M_0$ in $G(F)$. It is known that any element in $W_0^G$ admits a representative in $K$. 

By a Levi subgroup of $G$, we mean a group $M$ containing $M_0$ which is the Levi component of some parabolic subgroup of $G$. For such a group $M$, set $K_M:=M(F)\cap K$. Then the triplet $(M,K_M,M_0)$ satisfies the same hypotheses as $(G,K,M_0)$. Denote by $\msf^G(M)$, $\msp^G(M)$ and $\msl^G(M)$ the set of parabolic subgroups of $G$ containing $M$, parabolic subgroups of $G$ with Levi factor $M$ and Levi subgroups of $G$ containing $M$ respectively. 

For $P\in\msf^G(M_0)$, let $M_P$ be the Levi component containing $M_0$ and $N_P$ the unipotent radical. Denote $A_P:=A_{M_P}$ and $\fa_P:=\fa_{M_P}$ whose dual $\BR$-linear space is denoted by $\fa_P^*$. Define a map $H_P: G(F)\ra\fa_{M_P}$ by
$$ H_P(mnk)=H_{M_P}(m) $$
for all $m\in M_P(F)$, $n\in N_P(F)$ and $k\in K$. Let $\ov{P}\in\msp^G(M_P)$ be the parabolic subgroup opposite to $P$. 

For $P\subseteq Q$ a pair of parabolic subgroups in $\msf^G(M_0)$, the restriction $X(M_Q)_F\hookrightarrow X(M_P)_F$ induces a pair of dual maps $\fa_P\twoheadrightarrow\fa_Q$ and $\fa_Q^*\hookrightarrow \fa_P^*$. Let $\fa_P^Q$ be the kernel of the former map $\fa_P\twoheadrightarrow\fa_Q$. Set $\Delta_P^Q$ to be the set of simple roots for the action of $A_P$ on $P\cap M_Q$. Denote by $(\Delta_P^Q)^\vee$ the set of ``coroots'' as in \cite[p. 26]{MR2192011}. Then $(\Delta_P^Q)^\vee$ is a basis of the $\BR$-linear space $\fa_P^Q$. 

\subsection{Heights}
We fix a height function $\|\cdot\|: G(F)\ra\BR$ as in \cite[\S4]{MR1114210}. It satisfies the following properties: 

(1) $\|x\|\geq 1, \forall x\in G(F)$; 

(2) $\|xy\|\leq \|x\| \|y\|, \forall x, y\in G(F)$; 

(3) there exist constants $c>0$ and $N\in\BN$ such that $\|x^{-1}\|\leq c \|x\|^N, \forall x\in G(F)$. 

If $P\in\msf^G(M_0)$, for any $x\in G(F)$, we can choose $m_P(x)\in M_P(F), n_P(x)\in N_P(F)$ and $k_P(x)\in K$ such that $x=m_P(x)n_P(x)k_P(x)$. Then

(4) there exist constants $c>0$ and $N\in\BN$ such that $\|m_P(x)\| + \|n_P(x)\| \leq c \|x\|^N, \forall x\in G(F)$. 

We also fix a Euclidean norm (still denoted by $\|\cdot\|$) on the $\BR$-linear space $\fa_{M_0}$ which is invariant under the action of $W_0^G$ on $\fa_{M_0}$. Then

(5) there exist $c_1,c_2>0$ such that
$$ c_1(1+\log\|y\|) \leq 1+\|H_{M_0}(y)\| \leq c_2(1+\log\|y\|), \forall y\in M_0(F). $$

In addition, we require that $\|\cdot\|$ is a norm on $G(F)$ in the sense of \cite[\S18.2]{MR2192014}. This is possible. For example, for $G=GL_n$, by writing $(g,g^{-1})=(g_{ij},h_{ij})_{1\leq i,j\leq n}$, one may define $\|g\|:=\sup\limits_{i,j} \{|g_{ij}|_F, |h_{ij}|_F\}$ for $g\in G(F)$. Since $\{g_{ij},h_{ij}\}_{1\leq i,j\leq n}$ is a set of generators for the ring of regular functions of $G$ (viewed as an affine variety over $F$), this defines a norm in the sense of \cite[\S18.2]{MR2192014} on $GL_n(F)$. For general $G$, one can choose an closed embedding $G\ra GL_n$ over $F$ and define the norm on $G(F)$ by the pull-back of the norm on $GL_n(F)$. By \cite[Proposition 18.1.(2)]{MR2192014}, this defines a norm in the sense of \cite[\S18.2]{MR2192014} on $G(F)$. 

\subsection{Functions and distributions}
Let $\fg:=\Lie(G)$. More generally, we shall use a minuscule Fraktur letter to denote the Lie algebra of its corresponding algebraic group. Denote by $\Ad$ the adjoint action of $G$ on itself or $\fg$. The adjoint action of $\fg$ on itself is denoted by $\ad$. 

For a locally compact and totally disconnected topological space $X$ (e.g. $G(F)$ or $\fg(F)$), denote by $\CC_c^\infty(X)$ the space of locally constant, compactly supported, complex-valued functions on $X$. For $f\in\CC_c^\infty(X)$, denote by $\Supp(f)$ its support. Denote by $\CC_c^\infty(X)^\ast$ the space of distributions on $X$, i.e., the linear dual of $\CC_c^\infty(X)$. 

Suppose that there is a left action of $G(F)$ on such an $X$. Then $G(F)$ acts on $\CC_c^\infty(X)$ (or more generally the space of complex functions on $X$) from the left by
$$ g\cdot f(x):=f(g^{-1}\cdot x), \forall g\in G(F), f\in\CC_c^\infty(X), x\in X. $$
Moreover, $G(F)$ acts on $\CC_c^\infty(X)^\ast$ from the left by
$$ g\cdot d(f):=d(g^{-1}\cdot f), \forall g\in G(F), d\in\CC_c^\infty(X)^\ast, f\in\CC_c^\infty(X). $$
Let $\eta: G(F)\ra\BC^\times$ be a locally constant character. We say a function $f\in\CC_c^\infty(X)$ (resp. a distribution $d\in\CC_c^\infty(X)^\ast$) is $\eta$-invariant if $g\cdot f=\eta(g)f$ (resp. $g\cdot d=\eta(g)d$) for all $g\in G(F)$. For trivial $\eta$, we simply say that such a function (resp. distribution) is invariant. 

\subsection{Haar measures}\label{haarmeasure}
Fix the Haar measure on $K$ such that $\vol(K)=1$. Following \cite[\S I.4]{MR1344131}, for all $P\in\msf^G(M_0)$, we fix a Haar measure on $N_P(F)$ such that
$$ \int_{N_P(F)} \exp (2\rho_{\ov{P}}(H_{\ov{P}}(n))) dn =1, $$
where $\rho_{\ov{P}}$ is the half of the sum of roots (with multiplicity) associated to the parabolic subgroup $\ov{P}$ opposite to $P$. From \cite[p. 12]{MR1114210}, for all $M\in\msl^G(M_0)$, there are compatible Haar measures on $G(F)$ and $M(F)$ such that for all $P\in\msp^G(M)$ and $f\in\CC_c^\infty(G(F))$, we have
$$ \int_{G(F)} f(x) dx = \int_{M(F)\times N_P(F)\times K} f(mnk) dkdndm. $$
We shall fix such measures. 

For a $F$-split torus $T$, we choose the Haar measure on $T(F)$ such that the maximal compact subgroup of $T(F)$ is of volume $1$. For a general torus $T$, we choose the Haar measure on $T$ such that the induced measure on $T(F)/A_T(F)$ satisfies $\vol(T(F)/A_T(F))=1$. 

Notice that if $M_0$ is a torus, we have associated to it two measures. However, it will be clear which one should be used according to the context.  

Fix open neighbourhoods $V_\fg$ of $0$ in $\fg(F)$ and $V_G$ of $1$ in $G(F)$ such that the exponential map $\exp: V_\fg\ra V_G$ induces a homeomorphism between them. 
Choose the unique Haar measure on $\fg(F)$ such that $\exp$ preserves the measures. Similarly, we obtain Haar measures on $F$-points of Lie algebras of algebraic subgroups of $G$. 

From the fixed Euclidean norm $\|\cdot\|$ on $\fa_{M_0}$, we deduce measures on $\fa_{M_0}$ and its subspaces.

\subsection{$(G,M)$-families}
Following \cite[p. 15]{MR625344}, we define
$$ \theta_P^Q(\lambda):=\vol(\fa_P^Q/\BZ(\Delta_P^Q)^\vee)^{-1}\prod_{\alpha^\vee\in(\Delta_P^Q)^\vee} \lambda(\alpha^\vee), \forall \lambda\in i\fa_P^\ast, $$
where $\BZ(\Delta_P^Q)^\vee$ denotes the lattice in $\fa_P^Q$ generated by $(\Delta_P^Q)^\vee$. 

Suppose that $M\in\msl^G(M_0)$ and that $Q\in\msf^G(M)$. Let $(c_P)_{P\in\msp^G(M)}$ be a $(G,M)$-family in the sense of \cite[p. 36]{MR625344}. By \cite[Lemma 6.2]{MR625344}, we can define
$$ c_M^Q:=\lim_{\lambda\ra 0} \sum_{\{P\in\msp^G(M): P\subseteq Q\}} c_P(\lambda) \theta_P^Q(\lambda)^{-1}. $$
We sometimes write $c_M:=c_M^G$ if $Q=G$. 

An important example is following. According to \cite[p. 40-41]{MR625344}, for $x\in G(F)$, 
$$ v_P(\lambda,x):=e^{-\lambda(H_P(x))}, \forall \lambda\in i\fa_M^\ast, P\in\msp^G(M), $$
is a $(G,M)$-family (denoted by $(v_P(x))_{P\in\msp^G(M)}$). Then we obtain a function
\begin{equation}\label{localwtfn}
 v_M^Q(x):=\lim_{\lambda\ra 0} \sum_{\{P\in\msp^G(M): P\subseteq Q\}} v_P(\lambda,x) \theta_P^Q(\lambda)^{-1}, \forall x\in G(F). 
 \end{equation}

For a smooth function $c_P(\lambda)$ on $i\fa_P^\ast$, we can associate to it a smooth function $c'_P(\lambda)$ on $i\fa_P^\ast$ as in \cite[(6.3) in \S6]{MR625344}. Denote by $c'_P$ the value of $c'_P(\lambda)$ at $\lambda=0$. Let $(c_P)_{P\in\msp^G(M)}$ and $(d_P)_{P\in\msp^G(M)}$ be two $(G,M)$-families, we define their product $((cd)_P)_{P\in\msp^G(M)}$ in the obvious way and have the following product formula (see \cite[Lemma 6.3]{MR625344})
\begin{equation}\label{split6.3}
 (cd)_M=\sum_{Q\in\msf^G(M)} c'_Q d_M^Q. 
\end{equation}

\subsection{The maps $d_M^G$ and $s$}\label{gmmaps}
Suppose that $M\in\msl^G(M_0)$. As in \cite[p. 356]{MR928262}, we define a map
$$ d_M^G: \msl^G(M)\times\msl^G(M)\ra\BR_{\geq 0} $$
such that for all $(L_1,L_2)\in\msl^G(M)\times\msl^G(M)$, 

(1) $d_M^G(G,M)=d_M^G(M,G)=1$; 

(2) $d_M^G(L_1,L_2)=d_M^G(L_2,L_1)$; 

(3) $d_M^G(L_1,L_2)\neq 0$ if and only if $\fa_M^G=\fa_M^{L_1}\oplus\fa_M^{L_2}$. 

Following \cite[\S II.4]{MR1344131}, we also choose a map (not unique)
$$ s: \msl^G(M)\times\msl^G(M)\ra\msf^G(M)\times\msf^G(M) $$
such that for all $(L_1,L_2)\in\msl^G(M)\times\msl^G(M)$, 

(4) $s(L_1, L_2)\in\msp^G(L_1)\times\msp^G(L_2)$; 

(5) if $s(L_1,L_2)=(Q_1,Q_2)$, then $s(L_2,L_1)=(\ov{Q_2},\ov{Q_1})$; 

(6) (splitting formula) if $(c_P)_{P\in\msp^G(M)}$ and $(d_P)_{P\in\msp^G(M)}$ are $(G,M)$-families, we have the equality
$$ (cd)_M=\sum_{L_1,L_2\in\msl^G(M)} d_M^G(L_1,L_2)c_M^{Q_1}c_M^{Q_2}, $$
where $(Q_1,Q_2):=s(L_1,L_2)$; 

(7) (descent formula) if $(c_P)_{P\in\msp^G(M)}$ is a $(G,M)$-family and $L\in\msl^G(M)$, we have the equality
$$ c_L=\sum_{L'\in\msl^G(M)} d_M^G(L,L')c_M^{Q'}, $$
where $Q'$ denotes the second component of $s(L,L')$.


\section{\textbf{Symmetric pairs}}\label{secsympai}

\subsection{General cases}\label{generalcases}

Following \cite[Definition 7.1.1]{MR2553879}, by a symmetric pair, we mean a triple $(G, H, \theta)$ where $H\subseteq G$ are a pair of reductive groups, and $\theta$ is an involution on $G$ such that $H$ is the subgroup of fixed points of $\theta$. 

Suppose that $(G, H, \theta)$ is a symmetric pair. Let $\fg:=\Lie(G)$ and $\fh:=\Lie(H)$. Write $d\theta$ for the differential of $\theta$. Then 
$$ \fh=\{X\in\fg:(d\theta)(X)=X\}. $$ 
Define 
$$ S:=\{x\in G: x\theta(x)=1\}, $$
on which $G$ acts by 
$$ x\cdot s:=xs\theta(x)^{-1}, \forall x\in G, s\in S. $$
We have a symmetrization map $G/H\ra S$ defined by $x\mapsto x\theta(x)^{-1}$ which identifies the symmetric space $G/H$ with the $G$-orbit of the neutral element on $S$. The induced action of $H$ on $S$ is the conjugation by $H$. 
Let $\fs$ be the tangent space at the neutral element of $S$. Then 
$$ \fs=\{X\in\fg: (d\theta)(X)=-X\}, $$ 
on which $H$ acts by restriction of the adjoint action. Notice that 
$$ \fg=\fh\oplus\fs. $$ 
For $X\in\fs$, denote by $H_X$ (resp. $\fg_X$, $\fh_X$, $\fs_X$) the centraliser of $X$ in $H$ (resp. $\fg$, $\fh$, $\fs$). 

We say an element $X\in\fs$ is semi-simple if the orbit $\Ad(H)(X)$ is Zariski closed in $\fs$. 
By \cite[Fact A, p. 108-109]{MR1375304}, $X\in\fs(F)$ is semi-simple if and only if $\Ad(H(F))(X)$ is closed in $\fs(F)$ in the analytic topology. We say an element $X\in\fs$ is regular if $H_X$ has minimal dimension. It is known that $X\in\fs$ is regular if and only if $\fg_X$ (or $\fs_X$) has minimal dimension (cf. \cite[Proposition 7]{MR311837}). Denote by $\fs_\rs$ the subset of $\fs$ consisting of regular semi-simple elements in $\fs$. 
For $X\in\fs$, consider the characteristic polynomial $\det(\lambda-\ad(X)|_\fg)$. Denote by $\fD^\fs(X)$ the coefficient of the least power of $\lambda$ appearing nontrivially in this polynomial. Then $\fD^\fs(X)$ is an $H$-invariant polynomial on $\fs$. From \cite[end of p. 107]{MR1375304}, we know that $X\in\fs$ is regular semi-simple if and only if $\fD^\fs(X)\neq 0$. Thus $\fs_\rs$ is a principal Zariski open subset of $\fs$. 

By a Cartan subspace of $\fs$, we mean a maximal abelian subspace for the Lie bracket $\fc\subseteq\fs$ (defined over $F$) consisting of semi-simple elements. For such $\fc$, let $\fc_\reg:=\fc\cap\fs_\rs$. By \cite[Theorem 5.1 and Corollary of Theorem 5.2]{MR404361}, if $\fc\subseteq\fs$ is a Cartan subspace, then the centraliser of $\fc$ in $\fs$ is $\fc$ itself. By \cite[Corollary 2 of Theorem 3.2]{MR404361}, if $X\in\fs_\rs(F)$, then $\fs_X$ is a Cartan subspace of $\fs$. Hence, if $X\in\fc_\reg(F)$ with $\fc\subseteq\fs$ being a Cartan subspace, then $\fs_X=\fc$ because of the maximality. By \cite[Corollary 2 of Theorem 4.1]{MR404361}, all Cartan subspaces of $\fs$ are of the form $\fs_X$ where $X\in\fs_\rs(F)$. By \cite[p. 323]{MR366941}, all Cartan subspaces of $\fs$ are conjugate by $H$ over an algebraic closure of $F$. In particular, all Cartan subspaces of $\fs$ have the same dimension. 

For a Cartan subspace $\fc\subseteq\fs$, denote by $T_\fc$ the centraliser of $\fc$ in $H$ and set $\ft_\fc:=\Lie(T_\fc)$. It is known that if $X\in\fc_\reg(F)$ with $\fc\subseteq\fs$ being a Cartan subspace, then $H_X=T_\fc$ (see \cite[p. 112]{MR1375304}). Note that $T_\fc$ is not necessarily a torus in general (though it is always a torus in cases of main interest in this paper thanks to Corollary \ref{corsselemts1}.3)). 

For $X\in\fs(F)$, we define the Weyl discriminant factor 
$$ |D^\fs(X)|_F:=|\fD^\fs(X)|_F^{1/2}. $$
From the above discussion, we see that if $X\in\fc(F)$ with $\fc\subseteq\fs$ being a Cartan subspace, then 
\begin{equation}\label{weyldisc}
 |D^\fs(X)|_F=|\det(\ad(X)|_{\fh/\ft_\fc\oplus\fs/\fc})|_F^{1/2}.  
 \end{equation}

For a Cartan subspace $\fc\subseteq\fs$, set $W(H, \fc):=\Norm_{H(F)}(\fc)/T_\fc(F)$ to be its Weyl group, where $\Norm_{H(F)}(\fc)$ denotes the normaliser of $\fc$ in $H(F)$. Fix a set $\mst(\fs)$ of representatives for $H(F)$-conjugacy classes of Cartan subspaces in $\fs$, which is a finite set by \cite[p. 105]{MR1375304}. Then we have the Weyl integration formula (see \cite[p. 106]{MR1375304})
\begin{equation}\label{wif}
 \int_{\fs(F)} f(X) dX=\sum_{\fc\in\mst(\fs)} |W(H, \fc)|^{-1} \int_{\fc_\reg(F)} |D^\fs(X)|_F \int_{T_\fc(F)\bs H(F)} f(\Ad(x^{-1})(X)) dx dX 
\end{equation}
for all $f\in\CC_c^\infty(\fs(F))$. Recall that the adjoint action induces a local isomorphism $\beta: (T_\fc(F)\bs H(F))\times\fc_\reg(F)\ra\fs_\rs(F)$ of $F$-analytic manifolds, whose image is open in $\fs_\rs(F)$. Here we should use compatible Haar measures on $\fs(F)$ and $\fc_\reg(F)$, i.e., we require that $\beta$ should preserve the measures. For particular cases to be considered, we shall fix Haar measures on $\fs(F)$ in the following sections. Notice that we shall not use the Haar measure on $\fc_\reg(F)$ obtained via the exponential map. 

The lemma below makes the definition of Fourier transform on $\fs(F)$ possible. 

\begin{lem}
Let $(G,H,\theta)$ be a symmetric pair. Then there exists a $G$-invariant $\theta$-invariant non-degenerate symmetric bilinear form $\langle\cdot,\cdot\rangle$ on $\fg$. In particular, $\fg=\fh\oplus\fs$ is an orthogonal direct sum with respect to $\langle\cdot,\cdot\rangle$, and the restriction of $\langle\cdot,\cdot\rangle$ to $\fh$ or $\fs$ is non-degenerate. 
\end{lem}

\begin{proof}
This is \cite[Lemma 7.1.9]{MR2553879}. 
\end{proof}

The following lemma is a special case of \cite[Lemma 3.10]{MR3245011}, which is an analogue of Harish-Chandra's compactness lemma \cite[Lemma 25]{MR0414797}. 

\begin{lem}\label{lem25}
Let $\sigma_\fs$ be a compact subset of $\fs(F)$. Suppose that $\fc$ is a Cartan subspace of $\fs$. Let $\sigma_\fc$ be a compact subset of $\fc_\reg(F)$. Then
$$ \{x\in T_\fc(F)\bs H(F): \Ad(x^{-1})(\sigma_\fc)\cap\sigma_\fs\neq\emptyset\} $$
is relatively compact in $T_\fc(F)\bs H(F)$. 
\end{lem}

\begin{proof}
Choose an arbitrary $X\in\sigma_\fc$. We have $H_X=T_\fc$. Let $N_{\Ad(H)(X), X}^\fs$ be the normal space (see \cite[Notation 2.3.3]{MR2553879}) to $\Ad(H)(X)$ in $\fs$ at the point $X$. 
Note that $N_{\Ad(H)(X),X}^\fs\simeq\fs/\ad(\fh)(X)$. From 
$$ \fg=\fg_X\oplus\ad(\fg)(X)=\fh_X\oplus\fs_X\oplus\ad(\fh)(X)\oplus\ad(\fs)(X), $$
we see that $\fs=\fs_X\oplus\ad(\fh)(X)$. Then $N_{\Ad(H)(X),X}^\fs\simeq\fs_X$ as $H_X$-spaces. Since $X\in\fc_\reg$, $\fc_\reg\subseteq\fs_X$ is an \etale Luna slice at $X$ in the sense of \cite[Theorem A.2.3]{MR2553879}. Thus we can apply \cite[Lemma 3.10]{MR3245011}. 
\end{proof}

Denote by $F[\fs]^H$ the $F$-algebra of $H$-invariant polynomial functions on $\fs$. For $\fc\subseteq\fs$ a Cartan subspace, denote by $F[\fc]^{W(H, \fc)}$ the $F$-algebra of $W(H, \fc)$-invariant polynomial functions on $\fc$. An analogue of Chevalley restriction theorem holds for symmetric spaces. 

\begin{lem}\label{chevalley}
	1) Let $\fc$ be a Cartan subspace of $\fs$. The restriction 
	$$ F[\fs]^H\ra F[\fc]^{W(H, \fc)} $$
	is an isomorphism. 
	
	2) The $F$-algebra $F[\fs]^H$ has a finite system of algebraically independent homogeneous generators. 
\end{lem}

\begin{proof}
	See \cite[Theorem 7 and Corollary of Theorem 8]{zbMATH03577478} and \cite[Proposition 2.3]{MR1375304}. 
\end{proof}

The next lemma is an analogue of \cite[Lemma 28]{MR0414797}. 

\begin{lem}\label{lem28}
Let $\sigma\subseteq\fs(F)$ be a compact subset. Let $\fc$ be a Cartan subspace of $\fs$. Then $\fc(F)\cap\Cl(\Ad(H(F))(\sigma))$ is relatively compact in $\fc(F)$, where $\Cl$ denotes the closure of a subset in $\fs(F)$. 
\end{lem}

\begin{proof}
This is \cite[Lemma 6.12]{MR3414387}, whose proof relying on Lemma \ref{chevalley} applies to an arbitrary symmetric pair. 
\end{proof}

The following lemma is an analogue of \cite[Lemme III.4]{MR1344131}. 

\begin{lem}\label{lemIII.4}
Let $\sigma\subseteq \fs(F)$ be a compact subset. Let $\fc$ be a Cartan subspace of $\fs$ and $T_\fc$ the centraliser of $\fc$ in $H$. Then there exists $c_\sigma>0$ such that for all $x\in H(F)$ and $X\in\fc_\reg(F)$ satisfying $\Ad(x^{-1})(X)\in\sigma$, we have
$$ \inf_{\tau\in T_\fc(F)} \log\|\tau x\| \leq c_\sigma\sup\{1, -\log |D^\fs(X)|_F\}. $$
\end{lem}

\begin{proof}
Let $\|\cdot\|_{T_\fc\bs H}$ be any norm on $(T_\fc\bs H)(F)$ in the sense of \cite[\S18.2]{MR2192014}. Applying the argument of \cite[Lemma 20.3]{MR2192014} to the finite morphism
$$ \beta: (T_\fc\bs H)\times\fc_\reg\ra\fs_\rs $$
of affine algebraic varieties defined by $\beta(x,X):=\Ad(x^{-1})(X)$, we show the inequality
$$ \log\|x\|_{T_\fc\bs H}\leq c_\sigma\sup\{1, -\log |D^\fs(X)|_F\}. $$
By \cite[Proposition 18.3]{MR2192014}, the quotient $H\ra T_\fc\bs H$ has the norm descent property in the sense of \cite[\S18.6]{MR2192014}. That is to say, the restriction of $\|\cdot\|_{T_\fc\bs H}$ to $T_\fc(F)\bs H(F)$ is equivalent to the abstract norm $\inf_{\tau\in T_\fc(F)} \|\tau \cdot\|$ on $T_\fc(F)\bs H(F)$. 
\end{proof} 

The lemma below is an analogue of \cite[Lemma 44]{MR0414797}. 

\begin{lem}\label{lem44}
There exists $\varepsilon>0$ such that the function $|D^\fs(X)|_F^{-\varepsilon}$ is locally integrable on $\fc(F)$ for any Cartan subspace $\fc$ of $\fs$. 
\end{lem}

\begin{proof}
See \cite[Lemma 4.3]{MR3245011}. 
\end{proof}

\begin{coro}\label{cor20.2}
For any $r\geq0$, the function $\sup\{1, -\log |D^\fs(X)|_F\}^r$ is locally integrable on $\fc(F)$ for any Cartan subspace $\fc$ of $\fs$.
\end{coro}

\begin{proof}
We have the elementary fact (cf. the proof of \cite[Corollary 20.2]{MR2192014}): for $\varepsilon>0$ and $r\geq0$, there exists $c>0$ such that
$$ \sup\{1, \log y\}^r\leq cy^\varepsilon+1,\forall y>0. $$
Then it suffices to apply Lemma \ref{lem44}. 
\end{proof}

We say an element $X\in\fs(F)$ is nilpotent if $0\in\Cl(\Ad(H(F))(X))$, where $\Cl$ denotes the closure of a subset in $\fs(F)$. From \cite[Lemmas 2.3.12 and 7.3.8]{MR2553879}, we know that $X\in\fs(F)$ is nilpotent if and only if it is a nilpotent element in $\fg$. Denote by $\CN^\fs$ the set of nilpotent elements in $\fs(F)$, which is a cone. The following lemma is an analogue of Jacobson-Morozov theorem. 

\begin{lem}\label{jm}
Let $(G,H,\theta)$ be a symmetric pair and $X\in\CN^\fs$. Then there exists a group homomorphism $\varphi: SL_2(F)\ra G(F)$ such that
$$ d\varphi\mat(0,1,0,0)=X, d\varphi\mat(0,0,1,0)\in\fs(F) \text{ and } \varphi\mat(t,,,t^{-1})\in H(F), \forall t\in F^\times. $$
\end{lem}

\begin{proof}
This is \cite[Lemma 7.1.11]{MR2553879}. 
\end{proof}

Let $X\in\fs_\rs(F)$. The orbital integral of $X$ is the distribution $I_X$ on $\fs(F)$ defined by
\begin{equation}\label{oiwochar}
 \forall f\in\CC_c^\infty(\fs(F)), I_X(f):=|D^\fs(X)|_F^{1/2}\int_{H_X(F)\bs H(F)} f(\Ad(x^{-1})(X)) dx. 
\end{equation}
The next lemma is an analogue of Harish-Chandra's submersion principle \cite[Theorem 11]{MR0414797}. 

\begin{lem}\label{thm11}
Let $I: \fs_\rs(F)\ra\BC$ be a function. The following conditions are equivalent: 

(1) $I$ is locally constant, invariant by the adjoint action of $H(F)$ and of support included in $\Ad(H(F))(\sigma)$ with $\sigma\subseteq\fs_\rs(F)$ a compact subset;

(2) there exists $f\in\CC_c^\infty(\fs_\rs(F))$ such that
$$ \forall X\in\fs_\rs(F), I(X)=I_X(f). $$
\end{lem}

\begin{proof}
For $\fc\in\mst(\fs)$, apply the argument of \cite[Lemme 6.1]{MR2856369} to the morphism
$$ \Phi_\fc: (T_\fc(F)\bs H(F))\times\fc_\reg(F)\ra\fs_\rs(F) $$
defined by $\Phi_\fc(x,X):=\Ad(x^{-1})(X)$. Then glue the results for all $\fc\in\mst(\fs)$ together. 
\end{proof}

\subsection{The case of $(G,H)$}\label{sectsympar1}
Let $D$ be a central division algebra over $F$. Let $\fg$ be the associative algebra of $2n\times 2n$ matrices with entries in $D$. Denote by $G:=\fg^\times=GL_{2n,D}$ the group of invertible elements in $\fg$. Equipped with the Lie bracket, $\fg$ is identified with the Lie algebra of $G$. Let $\epsilon:=\mat(1_n,,,-1_n)$. Denote by $H$ the subgroup of fixed points of the involution $\Ad(\epsilon)$ on $G$. Then $H=GL_{n,D}\times GL_{n,D}$ is diagonally embedded into $G$. Define 
$$ S:=\{x\in G: x\Ad(\epsilon)(x)=1\} $$
and its tangent space at the neutral element 
$$ \fs:=\{X\in\fg: \Ad(\epsilon)(X)=-X\}. $$
Then we have 
$$ \fg=\fh\oplus\fs $$
and 
$$ \dim\fs=\dim\fh=\frac{1}{2}\dim\fg. $$
Let $M_0$ be the group of diagonal elements in $G$, which is a common minimal Levi subgroup of $G$ and $H$. For a linear subspace $\fv\subseteq\fg$, we write $\fv^\times:=\fv\cap G$. 

Note that the $F$-rank of $G$ is $2n$. Denote $G_n:=G$, $H_n:=H$ and $\fs_n:=\fs$. Recall the following description of semi-simple elements and descendants (see \cite[Definition 7.2.2]{MR2553879}). 

\begin{prop}\label{sselemts1}
The results below hold for $F$ being any field of characteristic zero. 

1) An element $X$ of $\fs(F)$ is semi-simple if and only if it is $H(F)$-conjugate to an element of the form
$$ X(A):=
\left( \begin{array}{cccc}
0 & 0 & 1_{m} & 0 \\
0 & 0 & 0 & 0 \\
A & 0 & 0 & 0 \\
0 & 0 & 0 & 0 \\
\end{array} \right) $$
with $A\in GL_{m}(D)$ being semi-simple in the usual sense. More precisely, the set of $H(F)$-conjugacy classes of semi-simple elements of $\fs(F)$ is bijective to the set of pairs $(m, [A])$ where $0\leq m\leq n$ is an integer and $[A]$ is a semi-simple conjugacy class in $GL_{m}(D)$. Moreover, $X(A)$ is regular semi-simple if and only if $m=n$ and $A$ is regular semi-simple in $GL_{n}(D)$ in the usual sense. 

2) Let $X=X(A)\in\fs(F)$ be semi-simple. Then the descendant $(H_X, \fs_X)$ (as a representation) is isomorphic to 
$$ ((GL_{m,D})_A, (\fg\fl_{m,D})_A) \times (H_{n-m},\fs_{n-m}), $$
where $(GL_{m,D})_A$ (resp. $(\fg\fl_{m,D})_A$) denotes the centraliser of $A$ in $GL_{m,D}$ (resp. $\fg\fl_{m,D}$). 
\end{prop}

\begin{proof}
  This is stated in \cite[Proposition 5.2]{MR3299843} for a non-archimedean local field of characteristic zero, but its proof included in \cite[Propositions 2.1 and 2.2]{MR1394521} is valid for any field of characteristic zero. 
\end{proof}

\begin{coro}\label{corsselemts1}
We have 

	1) $\fs_\rs\subseteq\fs^\times$; 
	
	2) for $X\in\fs_\rs$, $H_X$ is a torus of dimension $\frac{1}{2} \sqrt{\dim\fg}$; 
	
	3) for $\fc\subseteq\fs$ a Cartan subspace, $T_\fc$ is a torus of dimension $\frac{1}{2} \sqrt{\dim\fg}$. 
\end{coro}

\begin{proof}
	1) and 2) are clear from Proposition \ref{sselemts1} applied to an algebraic closure of $F$. 3) results from 2) since for $\fc\subseteq\fs$ a Cartan subspace, $T_\fc=H_X$ for any $X\in\fc_\reg\neq\emptyset$. 
\end{proof}

\begin{lem}\label{intpar1}
Let $P$ be a parabolic subgroup of $G$. Then $P\cap H$ is a parabolic subgroup of $H$ if and only if $\epsilon\in P$. Moreover, if $\epsilon$ belongs to a Levi factor $M$ of $P$, then $M\cap H$ is a Levi factor of $P\cap H$. 
\end{lem}

\begin{proof}
One may consider all the groups over an algebraic closure of $F$. We first suppose that $P\cap H$ is a parabolic subgroup of $H$. Then $\epsilon\in\Cent(H)\subseteq P\cap H$, where $\Cent(H)$ denotes the centre of $H$. This establishes one direction. 

We now suppose that $\epsilon\in P$. Denote by $N$ the unipotent radical of $P$ and let $M$ be a Levi factor of $P$. By the argument in the last paragraph of the proof of \cite[Lemma 4.4]{zbMATH07499568}, we show that $\epsilon$ is $N$-conjugate to an element in $M$ with the help of \cite[Lemma 2.1]{MR518111} (actually we need its variant over a local field for the characteristic function of a singleton here, whose proof is similar). Then replacing $M$ by its $N$-conjugate if necessary, we may assume that $\epsilon\in M$. 

Let $G=GL(V)$ for a vector space $V=\oplus_{1\leq i\leq r} V_i$. Suppose that
$$ P=\{g\in G: g(V_1\oplus...\oplus V_i)\subseteq V_1\oplus...\oplus V_i, \forall 1\leq i\leq r\} $$
and that
$$ M=\{g\in G: g(V_i)\subseteq V_i, \forall 1\leq i\leq r\}. $$
Since $\epsilon\in M$, we have $\epsilon(V_i)\subseteq V_i$ for all $1\leq i\leq r$. Let $V_i^+$ (resp. $V_i^-$) be the $+1$(resp. $-1$)-eigenspace of $V_i$ under the action of $\epsilon$. For $\epsilon^2=1$, we have $V_i=V_i^+\oplus V_i^-$. Let $V^+:=\oplus_{1\leq i\leq r} V_i^+$ and $V^-:=\oplus_{1\leq i\leq r} V_i^-$. Then
$$ H=\{g\in G: g(V^+)\subseteq V^+, g(V^-)\subseteq V^-\}. $$
Hence, 
$$ P\cap H=\{g\in G: g(V_1^+\oplus...\oplus V_i^+)\subseteq V_1^+\oplus...\oplus V_i^+, g(V_1^-\oplus...\oplus V_i^-)\subseteq V_1^-\oplus...\oplus V_i^-, \forall 1\leq i\leq r\}. $$
It means exactly that $P\cap H$ is a parabolic subgroup of $H$ and proves the other direction. Morover, 
$$ M\cap H=\{g\in G: g(V_i^+)\subseteq V_i^+, g(V_i^-)\subseteq V_i^-, \forall 1\leq i\leq r\}. $$
That is to say, $M\cap H$ is a Levi factor of $P\cap H$. 
\end{proof}

Set $\omega:=\mat(0,1_n,1_n,0)$. For $P\in\msf^G(M_0)$, we say that $P$ is ``$\omega$-stable'' if $\omega\in P$. Denote by $\msf^{G,\omega}(M_0)$ the subset of $\omega$-stable parabolic subgroups in $\msf^G(M_0)$. For $M\in\msl^G(M_0)$, we say that $M$ is ``$\omega$-stable'' if $\omega\in M$. Notice that $M$ is $\omega$-stable if and only if $M=M_P$ for some $P\in\msf^{G,\omega}(M_0)$. 
Denote by $\msl^{G,\omega}(M_0)$ the subset of $\omega$-stable Levi subgroups in $\msl^G(M_0)$. Let $A_n$ be the group of diagonal matrices in $GL_n$, which is a minimal Levi subgroup of $GL_n$. For $M_n\in\msl^{GL_n}(A_n)$, denote by $M_{n,D}$ the reductive group over $F$ whose $F$-points are $M_n(D)$ and let $\fm_{n,D}:=\Lie(M_{n,D})$; in particular, the notion $A_{n,D}$ makes sense. There is a bijection $M_n\mapsto M=\mat(\fm_{n,D},\fm_{n,D},\fm_{n,D},\fm_{n,D})^\times$ between $\msl^{GL_n}(A_n)$ and  $\msl^{G,\omega}(M_0)$ (cf. \cite[Proposition 5.2]{zbMATH07499568}). We shall always use the notation $M_n$ to denote the preimage of $M$ under this bijection. Notice that if $M\in\msl^{G,\omega}(M_0)$ and $Q\in\msf^G(M)$, then $Q\in\msf^{G,\omega}(M_0)$. 

Suppose that $M\in\msl^{G,\omega}(M_0)$. If $X\in\fm\cap\fs$, then $A_M\subseteq H_X$. We say an element $X\in(\fm\cap\fs_\rs)(F)$ is $M$-elliptic if $A_M$ is the maximal $F$-split torus in $H_X$ (which is a torus by Corollary \ref{corsselemts1}.2)). Denote by $(\fm\cap\fs_\rs)(F)_\el$ the set of $M$-elliptic elements in $(\fm\cap\fs_\rs)(F)$. Write $M_H:=M\cap H$. For $X\in(\fm\cap\fs_\rs)(F)$, let $[X]_{M_H}$ be the $M_H(F)$-conjugacy class of $X$. Denote by $\Gamma_\el((\fm\cap\fs_\rs)(F))$ the set of $M_H(F)$-conjugacy classes in $(\fm\cap\fs_\rs)(F)_\el$. If $\fc\subseteq\fm\cap\fs$ is a Cartan subspace, then $A_M\subseteq T_\fc$. We say a Cartan subspace $\fc\subseteq\fm\cap\fs$ is $M$-elliptic if $A_M$ is the maximal $F$-split torus in $T_\fc$ (which is a torus by Corollary \ref{corsselemts1}.3)). Since $(M_H, \fm\cap\fs)$ appears as the product of some copies of the form $(H,\fs)$ in lower dimensions (cf. \cite[\S5.1]{zbMATH07499568}), we define $W(M_H,\fc)$ and $\mst(\fm\cap\fs)$ as in Section \ref{generalcases}. In our notation, we write $(\fm\cap\fs)_\rs$ for the Zariski open subset of $\fm\cap\fs$ consisting of regular semi-simple elements with respect to the $M_H$-action. By Proposition \ref{sselemts1} and its product form, we see that $\fm\cap\fs_\rs\subseteq(\fm\cap\fs)_\rs$. Since a Cartan subspace is generally the same as the centraliser of a regular semi-simple element, we also deduce that a Cartan subspace of $\fm\cap\fs$ is nothing but a Cartan subspace of $\fs$ contained in $\fm\cap\fs$. Denote by $\mst_\el(\fm\cap\fs)$ the subset of representatives that are $M$-elliptic Cartan subspaces in $\mst(\fm\cap\fs)$. 

\begin{lem}\label{countlemforwif2.1}
Let $M\in\msl^{G,\omega}(M_0)$ and $[X]_{M_H}\in\Gamma_\el((\fm\cap\fs_\rs)(F))$. 

1) Let $M'\in\msl^{G,\omega}(M_0)$ and $[X']_{M'_H}\in\Gamma_\el((\fm'\cap\fs_\rs)(F))$ be such that $X'$ is $H(F)$-conjugate to $X$. Then there exists 
$$ w\in\left\{\mat(w_n,,,w_n): w_n\in W_0^{GL_{n,D}}\right\}, $$
where $W_0^{GL_{n,D}}$ denotes the Weyl group of $(GL_{n,D}, A_{n,D})$, such that
$$ (\Ad(w)(M), [\Ad(w)(X)]_{\Ad(w)(M)\cap H})=(M', [X']_{M'_H}). $$ 

2) The cardinality of
$$ \{(M', [X']_{M'_H}): M'\in\msl^{G,\omega}(M_0), [X']_{M'_H}\in\Gamma_\el((\fm'\cap\fs_\rs)(F)), \text{$X'$ is $H(F)$-conjugate to $X$}\} $$
is
$$ |W_0^{GL_n}||W_0^{M_n}|^{-1}, $$
where $W_0^{GL_n}$ (resp. $W_0^{M_n}$) denotes the Weyl group of $(GL_n,A_n)$ (resp. $(M_n, A_n)$). 
\end{lem}

\begin{proof}
1) Let $x\in H(F)$ be such that $\Ad(x)(X)=X'$. Then $\Ad(x)(H_X)=H_{X'}$. Since $X\in(\fm\cap\fs_\rs)(F)_\el$ and $X'\in(\fm'\cap\fs_\rs)(F)_\el$, we have $\Ad(x)(A_M)=A_{M'}$ and thus $\Ad(x)(M)=M'$. As $x\in H(F)$, we have $\Ad(x)(M_H)=M'_H$. We see that $\Ad(x)(A_{M_0})\subseteq M'_H$ is a maximal $F$-split torus, so there exists $m'\in M'_H(F)$ such that $\Ad({m'}^{-1}x)(A_{M_0})=A_{M_0}$. That is to say, $w':={m'}^{-1}x\in\Norm_{H(F)}(A_{M_0})=\Norm_{H(F)}(M_0)$, where $\Norm_{H(F)}(A_{M_0})$ denotes the normaliser of $A_{M_0}$ in $H(F)$. Now $\Ad(x)(A_M)=A_{M'}$ implies that $\Ad(w')(A_M)=A_{M'}$. Because $M, M'\in\msl^{G,\omega}(M_0)$, it is shown in \cite[\S9.1]{zbMATH07499568} that any isomorphism $A_M\ra A_{M'}$ induced by $W_0^H$ can be given by $\left\{\mat(w_n,,,w_n): w_n\in W_0^{GL_{n,D}}\right\}$. Hence, there exists $w\in\left\{\mat(w_n,,,w_n): w_n\in W_0^{GL_{n,D}}\right\}$ such that $w^{-1}w'\in\Cent_{W_0^H}(A_M)=W_0^{M_H}$, where $\Cent_{W_0^H}(A_M)$ denotes the centraliser of $A_M$ in $W_0^H$. We can check that such a $w$ satisfies the condition in the lemma. 

2) By 1), the group $\left\{\mat(w_n,,,w_n): w_n\in W_0^{GL_{n,D}}\right\}$ acts transitively on this set. Let
$$ w\in\left\{\mat(w_n,,,w_n): w_n\in W_0^{GL_{n,D}}\right\} $$
be such that
$$ (\Ad(w)(M), [\Ad(w)(X)]_{\Ad(w)(M)\cap H})=(M, [X]_{M_H}). $$
Then $w\in M_H(F)$. Thus the condition on $w$ is equivalent to 
$$ w\in\left\{\mat(w_n,,,w_n): w_n\in W_0^{M_{n,D}}\right\}, $$
where $W_0^{M_{n,D}}$ denotes the Weyl group of $(M_{n,D},A_{n,D})$. We see that the cardinality of the set in the lemma is $|W_0^{GL_{n,D}}||W_0^{M_{n,D}}|^{-1}$ or $|W_0^{GL_{n}}||W_0^{M_{n}}|^{-1}$. 
\end{proof}

\begin{prop}\label{wif2.1}
For $f\in\CC_c^\infty(\fs(F))$, we have the equality 
\[\begin{split}
 \int_{\fs(F)} f(X) dX=&\sum_{M\in\msl^{G,\omega}(M_0)} |W_0^{M_n}| |W_0^{GL_n}|^{-1} \sum_{\fc\in\mst_\el(\fm\cap\fs)} |W(M_H, \fc)|^{-1} \int_{\fc_\reg(F)} |D^\fs(X)|_F \int_{A_M(F)\bs H(F)} \\
 &f(\Ad(x^{-1})(X)) dx dX. 
\end{split}\]
\end{prop}

\begin{proof}
From \cite[beginning of \S9.2]{zbMATH07499568}, any $H(F)$-conjugacy class in $\fs_\rs(F)$ is the image of a class $[X]_{M_H}\in\Gamma_\el((\fm\cap\fs_\rs)(F))$ for some $M\in\msl^{G,\omega}(M_0)$. By Lemma \ref{countlemforwif2.1}, the Weyl integration formula (\ref{wif}) can be written as the above equality (cf. \cite[p. 16-17]{MR1114210} and \cite[(3) in \S I.3]{MR1344131}). 
\end{proof}

For $P_n\in\msl^{GL_n}(A_n)$, denote by $P_{n,D}$ the algebraic group over $F$ whose $F$-points are $P_n(D)$ and let $\fp_{n,D}:=\Lie(P_{n,D})$. There is a bijection $P_n\mapsto P=\mat(\fp_{n,D},\fp_{n,D},\fp_{n,D},\fp_{n,D})^\times$ between $\msf^{GL_n}(A_n)$ and  $\msf^{G,\omega}(M_0)$ (cf. \cite[Proposition 5.2]{zbMATH07499568}). We shall always use the notation $P_n$ to denote the preimage of $P$ under this bijection. Following \cite[p. 1846]{MR3414387}, we shall fix the Haar measures on some subspaces of $\fs(F)$ as follows. Let $P\in\msf^{G,\omega}(M_0)$. Then we have $\fm_P=\mat(\fm_{P_{n,D}}, \fm_{P_{n,D}}, \fm_{P_{n,D}}, \fm_{P_{n,D}})$ and $\fn_P=\mat(\fn_{P_{n,D}}, \fn_{P_{n,D}}, \fn_{P_{n,D}}, \fn_{P_{n,D}})$. 
We have fixed the Haar measures on $\fm_{P_{n,D}}(F)$ and $\fn_{P_{n,D}}(F)$ in Section \ref{haarmeasure}. We shall choose the same Haar measure for any of the four copies in $\fm_P(F)$ or $\fn_P(F)$ under these identifications. In particular, we obtain the Haar measures on $(\fm_P\cap\fs)(F)$ and $(\fn_P\cap\fs)(F)$. 

\begin{lem}\label{jac1}
Let $Q\in\msf^{G,\omega}(M_0)$. For $Y\in(\fm_Q\cap\fs_\rs)(F)$, the map
$$ N_{Q_H}(F)\ra(\fn_Q\cap\fs)(F), n\mapsto \Ad(n^{-1})(Y)-Y $$
is an isomorphism of $F$-analytic manifolds whose Jacobian is $|D^\fs(Y)|_F^{1/2}|D^{\fm_Q\cap\fs}(Y)|_F^{-1/2}$. 
\end{lem}

\begin{proof}
See \cite[Lemma 8.1]{zbMATH07499568} for the isomorphism and the proof of \cite[Proposition 6.3.(ii)]{MR3414387} for the Jacobian.
\end{proof}

Fix a continuous and nontrivial unitary character $\Psi: F\ra\BC^\times$. Let $\langle\cdot,\cdot\rangle$ be the non-degenerate symmetric bilinear form on $\fg(F)$ defined by
$$ \langle X,Y\rangle:=\Trd(XY),\forall X,Y\in\fg(F), $$
where $\Trd$ denotes the reduced trace on $\fg(F)$. It is invariant by the adjoint action of $G(F)$ and $\Ad(\epsilon)$. For $f\in\CC_c^\infty(\fs(F))$, define its normalised Fourier transform $\hat{f}\in\CC_c^\infty(\fs(F))$ by
\begin{equation}\label{fouriertransform1}
 \forall X\in\fs(F), \hat{f}(X):=c_\Psi(\fs(F)) \int_{\fs(F)} f(Y) \Psi(\langle X,Y \rangle) dY, 
\end{equation}
where $c_\Psi(\fs(F))$ is the unique constant such that $\hat{\hat{f}}(X)=f(-X)$ for all $f\in\CC_c^\infty(\fs(F))$ and all $X\in\fs(F)$. For any $M\in\msl^{G,\omega}(M_0)$, the restriction of $\langle\cdot,\cdot\rangle$ on $\fm\cap\fs$ is non-degenerate. Then we can define similarly the normalised Fourier transform of $f\in\CC_c^\infty((\fm\cap\fs)(F))$. 

Suppose that $P\in\msf^{G,\omega}(M_0)$. Let $\eta$ be the quadratic character of $F^\times$ attached to a quadratic extension $E/F$. Denote by $\Nrd$ the reduced norm on $G(F)$. For $f\in\CC_c^\infty(\fs(F))$, we define a function (parabolic descent) $f_P^\eta\in\CC_c^\infty((\fm_P\cap\fs)(F))$ by
\begin{equation}\label{pardes1}
 f_P^\eta(Z):=\int_{K_H\times(\fn_P\cap\fs)(F)} f(\Ad(k^{-1})(Z+U)) \eta(\Nrd(k)) dUdk 
\end{equation}
for all $Z\in(\fm_P\cap\fs)(F)$. We show that $(\hat{f})_P^\eta=(f_P^\eta)^\string^$, so we shall denote it by $\hat{f}_P^\eta$ without confusion. In fact, the integral on $(\fn_P\cap\fs)(F)$ and the Fourier transform commute by our choices of Haar measures (see \cite[\S I.7]{MR1344131}); the commutativity of the integral on $K_H$ and the Fourier transform results from the $H(F)$-invariance of $\langle\cdot,\cdot\rangle$. 

The following result is an analogue of \cite[Theorem 13]{MR0414797}. 

\begin{prop}\label{bdoi1}
Let $f\in\CC_c^\infty(\fs(F))$. Then
$$ \sup_{X\in\fs_\rs(F)}|D^\fs(X)|_F^{1/2} \int_{H_X(F)\bs H(F)} |f(\Ad(x^{-1})(X))| dx<+\infty. $$
\end{prop}

\begin{proof}
It is proved in \cite[Theorem 6.11]{MR3414387} (see also \cite[p. 77]{MR3299843}) that for any fixed Cartan subspace $\fc$ of $\fs$, 
$$ \sup_{X\in\fc_\reg(F)}|D^\fs(X)|_F^{1/2} \int_{H_X(F)\bs H(F)} |f(\Ad(x^{-1})(X))| dx<+\infty. $$
Since $\mst(\fs)$ is a finite set and the orbital integral is constant on any $H(F)$-orbit, we obtain a uniform bound for all $X\in\fs_\rs(F)$. 
\end{proof}

The lemma below is an analogue of \cite[Theorem 15]{MR0414797}. 

\begin{lem}\label{thm15}
There exists $\varepsilon>0$ such that the function $|D^\fs(X)|_F^{-\frac{1}{2}-\varepsilon}$ is locally integrable on $\fs(F)$. 
\end{lem}

\begin{proof}
Choose $\varepsilon>0$ verifying the condition of Lemma \ref{lem44}. Let $f\in\CC_c^\infty(\fs(F))$ with $f\geq0$. By the Weyl integration formula (\ref{wif}), we have
$$ \int_{\fs(F)} |D^\fs(X)|_F^{-\frac{1}{2}-\varepsilon}f(X) dX=\sum_{\fc\in\mst(\fs)} |W(H, \fc)|^{-1} \int_{\fc_\reg(F)} |D^\fs(X)|_F^{\frac{1}{2}-\varepsilon} \int_{A_{T_\fc}(F)\bs H(F)} f(\Ad(x^{-1})(X)) dx dX. $$
The convergence of the right hand side results from Proposition \ref{bdoi1} and Lemmas \ref{lem28} and \ref{lem44}. 
\end{proof}

\begin{coro}\label{corII.1}
For any $r\geq0$, the function $|D^\fs(X)|_F^{-\frac{1}{2}}\sup\{1, -\log |D^\fs(X)|_F\}^r$ is locally integrable on $\fs(F)$. 
\end{coro}

\begin{proof}
It is the same as the proof of Corollary \ref{cor20.2}. 
\end{proof}

\subsection{The case of $(G',H')$}\label{sectsympar2}

Let $E$ be a quadratic extension of $F$. Let $\fg'$ be a central simple algebra over $F$ with a fixed embedding of $F$-algebras $E\hookrightarrow\fg'$. Let $\fh':=\Cent_{\fg'}(E)$ be the centraliser of $E$ in $\fg'$. By the double centraliser theorem (see \cite[Theorem 3.1 in Chapter IV]{milneCFT} for example), we know that $\fh'$ is a central simple algebra over $E$ and that $\dim\fh'=\frac{1}{2}\dim\fg'$. Denote by $G':={\fg'}^\times$ (resp. $H':={\fh'}^\times$) the group of invertible elements in $\fg'$ (resp. $\fh'$), which is considered as a reductive group over $F$ with Lie algebra $\fg'$ (resp. $\fh'$). Let $\alpha\in E\bs F$ be such that $\alpha^2\in F$. Then $E=F(\alpha)$ and $\Ad(\alpha)$ is an involution on $G'$ independant of the choice of $\alpha$. We see that $H'$ is the subgroup of fixed points of $\Ad(\alpha)$. Define 
$$ S':=\{x\in G': x\Ad(\alpha)(x)=1\} $$
and its tangent space at the neutral element 
$$ \fs':=\{Y\in\fg': \Ad(\alpha)(Y)=-Y\}. $$
Then we have 
$$ \fg'=\fh'\oplus\fs' $$ 
and 
$$ \dim\fs'=\dim\fh'=\frac{1}{2}\dim\fg'. $$ 
For a linear subspace $\fv'\subseteq\fg'$, we write ${\fv'}^\times:=\fv'\cap G'$. 

Notice that if $\dim\fg=\dim\fg'$, then after a base change to an algebraic closure of $F$ containing $E$, the symmetric pair $(G',H')$ is the same as the one $(G,H)$ defined in Section \ref{sectsympar1}. By Corollary \ref{corsselemts1}, we see that $\fs'_\rs\subseteq{\fs'}^\times$ and that for $Y\in\fs'_\rs$ (resp. $\fc'\subseteq\fs'$ a Cartan subspace), $H'_Y$ (resp. $T_{\fc'}$) is a torus of dimension $\frac{1}{2}\sqrt{\dim\fg'}$. 

\begin{lem}\label{intpar2}
Let $\wt{P'}$ be a parabolic subgroup of $G'$. Then $\wt{P'}\cap H'$ is a parabolic subgroup of $H'$ if and only if $\alpha\in \wt{P'}$. Moreover, if $\alpha$ belongs to a Levi factor $\wt{M'}$ of $\wt{P'}$, then $\wt{M'}\cap H'$ is a Levi factor of $\wt{P'}\cap H'$. 
\end{lem}

\begin{proof}
It is the same as the proof of Lemma \ref{intpar1}, except that one replaces $\epsilon$ with $\alpha$. 
\end{proof}

By the Wedderburn-Artin theorem, $G'$ is isomorphic to $GL_{n,D}$ for some positive integer $n$ and some central division algebra $D$ over $F$. Since $E$ embeds into $\fg'(F)$, we see that $n\deg(D)$ is even, where $\deg(D)$ denotes the degree of $D$. From the Noether-Skolem theorem (see \cite[Theorem 2.10 of Chapter IV]{milneCFT} for example), up to conjugation by $G'(F)$, the emdedding $H'\hookrightarrow G'$ is reduced to one of the two cases below (see \cite[\S2.1 and \S3.1]{MR3958071} and \cite[\S3.4]{MR4350885}). 

\textbf{Case I}: if $\deg(D)$ is even, then $(G',H')=(GL_{n,D},\Res_{E/F} GL_{n,D'})$, where $D':=\Cent_D(E)$ denoting the centraliser of $E$ in $D$ is a central division algebra over $E$ of degree $\frac{\deg(D)}{2}$. 
Let $M'_0\simeq (\Res_{E/F} \BG_{m,D'})^n$ (resp. $M'_{\wt{0}}\simeq (\BG_{m,D})^n$) be the subgroup of diagonal elements in $H'$ (resp. $G'$), which is a minimal Levi subgroup of $H'$ (resp. $G'$). For $M'\in\msl^{H'}(M'_0)$, denote  by $\wt{M'}$ the unique element in $\msl^{G'}(M'_{\wt{0}})$ such that $\wt{M'}\cap H=M'$. Then the map $M'\mapsto\wt{M'}$ induces a bijection between $\msl^{H'}(M'_0)$ and $\msl^{G'}(M'_{\wt{0}})$ (cf. \cite[Case I in \S3.4]{MR4350885}). Notice that we can identify $A_{M'}$ with $A_{\wt{M'}}$. 

\textbf{Case II}: if $\deg(D)$ is odd, then $(G',H')=(GL_{n,D},\Res_{E/F} GL_{\frac{n}{2},D\otimes_F E})$, where $D\otimes_F E$ is a central division algebra over $E$ of degree $\deg(D)$. 
Let $M'_0\simeq (\Res_{E/F} \BG_{m,D\otimes_F E})^\frac{n}{2}$ (resp. $M'_{\wt{0}}\simeq (\BG_{m,D})^n$) be the subgroup of diagonal elements in $H'$ (resp. $G'$), which is a minimal Levi subgroup of $H'$ (resp. $G'$). Denote by $\msl^{G'}(M'_{\wt{0}}, M'_0)$ the subset of elements in $\msl^{G'}(M'_{\wt{0}})$ containing $M'_0$. For $M'\in\msl^{H'}(M'_0)$, denote  by $\wt{M'}$ the unique element in $\msl^{G'}(M'_{\wt{0}}, M'_0)$ such that $\wt{M'}\cap H=M'$. Then the map $M'\mapsto\wt{M'}$ induces a bijection between $\msl^{H'}(M'_0)$ and $\msl^{G'}(M'_{\wt{0}}, M'_0)$ (cf. \cite[Case II in \S3.4]{MR4350885}). Notice that we can identify $A_{M'}$ with $A_{\wt{M'}}$. 

Suppose that $M'\in\msl^{H'}(M'_0)$. If $Y\in\wt{\fm'}\cap\fs'$, then $A_{M'}\subseteq H'_Y$. We say an element $Y\in(\wt{\fm'}\cap\fs'_\rs)(F)$ is $M'$-elliptic if $A_{M'}$ is the maximal $F$-split torus in $H'_Y$ (which is a torus by Corollary \ref{corsselemts1}.2)). Denote by $(\wt{\fm'}\cap\fs'_\rs)(F)_\el$ the set of $M'$-elliptic elements in $(\wt{\fm'}\cap\fs'_\rs)(F)$. For $Y\in(\wt{\fm'}\cap\fs'_\rs)(F)$, let $[Y]_{M'}$ be the $M'(F)$-conjugacy class of $Y$. Denote by $\Gamma_\el((\wt{\fm'}\cap\fs'_\rs)(F))$ the set of $M'(F)$-conjugacy classes in $(\wt{\fm'}\cap\fs'_\rs)(F)_\el$. If $\fc'\subseteq\wt{\fm'}\cap\fs'$ is a Cartan subspace, then $A_{M'}\subseteq T_{\fc'}$. We say a Cartan subspace $\fc'\subseteq\wt{\fm'}\cap\fs'$ is $M'$-elliptic if $A_{M'}$ is the maximal $F$-split torus in $T_{\fc'}$ (which is a torus by Corollary \ref{corsselemts1}.3)). Since $(M', \wt{\fm'}\cap\fs')$ appears as the product of some copies of the form $(H',\fs')$ in lower dimensions (cf. \cite[beginning of \S5]{MR4350885}), we define $W(M',\fc')$ and $\mst(\wt{\fm'}\cap\fs')$ as in Section \ref{generalcases}. In our notation, we write $(\wt{\fm'}\cap\fs')_\rs$ for the Zariski open subset of $\wt{\fm'}\cap\fs'$ consisting of regular semi-simple elements with respect to the $M'$-action. By a base change to an algebraic closure of $F$ containing $E$ and the relevant discussion in Section \ref{sectsympar1}, we see that $\wt{\fm'}\cap\fs'_\rs\subseteq(\wt{\fm'}\cap\fs')_\rs$ and that a Cartan subspace of $\wt{\fm'}\cap\fs'$ is nothing but a Cartan subspace of $\fs'$ contained in $\wt{\fm'}\cap\fs'$. Denote by $\mst_\el(\wt{\fm'}\cap\fs')$ the subset of representatives that are $M'$-elliptic Cartan subspaces in $\mst(\wt{\fm'}\cap\fs')$. 

\begin{lem}[cf. Lemma \ref{countlemforwif2.1}]\label{countlemforwif2.2}
Let $M'\in\msl^{H'}(M'_0)$ and $[Y']_{M'}\in\Gamma_\el((\wt{\fm'}\cap\fs'_\rs)(F))$. 

1) Let $M\in\msl^{H'}(M'_0)$ and $[Y]_{M}\in\Gamma_\el((\wt{\fm}\cap\fs'_\rs)(F))$ be such that $Y$ is $H'(F)$-conjugate to $Y'$. Then there exists $w\in W_0^{H'}$ such that
$$ (\Ad(w)(M'), [\Ad(w)(Y')]_{\Ad(w)(M')})=(M, [Y]_{M}). $$ 

2) The cardinality of
$$ \{(M, [Y]_{M}): M\in\msl^{H'}(M'_0), [Y]_{M}\in\Gamma_\el((\wt{\fm}\cap\fs'_\rs)(F)), \text{$Y$ is $H'(F)$-conjugate to $Y'$}\} $$
is 
$$ |W_0^{H'}||W_0^{M'}|^{-1}. $$
\end{lem}

\begin{prop}\label{wif2.2}
For $f'\in\CC_c^\infty(\fs'(F))$, we have the equality 
\[\begin{split}
 \int_{\fs'(F)} f'(Y) dY=&\sum_{M'\in\msl^{H'}(M'_0)} |W_0^{M'}||W_0^{H'}|^{-1} \sum_{\fc'\in\mst_\el(\wt{\fm'}\cap\fs')} |W(M', \fc')|^{-1} \int_{\fc'_\reg(F)} |D^{\fs'}(Y)|_F \int_{A_{M'}(F)\bs H'(F)} \\
 &f'(\Ad(x^{-1})(Y)) dx dY. 
\end{split}\]
\end{prop}

\begin{proof}
From \cite[Lemma 9.1]{MR4350885}, any $H'(F)$-conjugacy class in $\fs'_\rs(F)$ is the image of a class $[Y]_{M'}\in\Gamma_\el((\wt{\fm'}\cap\fs'_\rs)(F))$ for some $M'\in\msl^{H'}(M'_0)$ in our case. By Lemma \ref{countlemforwif2.2}, the Weyl integration formula (\ref{wif}) can be written as the above equality (cf. \cite[p. 16-17]{MR1114210} and \cite[(3) in \S I.3]{MR1344131}). 
\end{proof}

In both of \textbf{Case I} and \textbf{Case II}, for $P'\in\msf^{H'}(M'_0)$, denote by $\wt{P'}$ the unique element in $\msf^{G'}(\wt{M'_0})$ such that $\wt{P'}\cap H=P'$. Then the map $P'\mapsto\wt{P'}$ induces a bijection between $\msf^{H'}(M'_0)$ and $\msf^{G'}(\wt{M'_0})$ (cf. \cite[\S3.4]{MR4350885}). 
Let $\tau\in D^\times$ in \textbf{Case I} (resp. $\tau\in GL_2(D)$ in \textbf{Case II}) be an element such that $\Ad(\alpha)(\tau)=-\tau$. Let $P'\in\msf^{H'}(M'_0)$. Then we have $\fm_{\wt{P'}}\cap\fs'=\fm_{P'}\tau=\tau\fm_{P'}$ and $\fn_{\wt{P'}}\cap\fs'=\fn_{P'}\tau=\tau\fn_{P'}$ by \cite[Proposition 3.12]{MR4350885}. We have fixed the Haar measures on $\fm_{P'}(F)$ and $\fn_{P'}(F)$ in Section \ref{haarmeasure}. We shall choose the same Haar measures on $(\fm_{\wt{P'}}\cap\fs')(F)$ and $(\fn_{\wt{P'}}\cap\fs')(F)$ using above identifications induced by $\tau$. Such Haar measures depend on the choice of $\tau$. 

\begin{lem}\label{jac2}
Let $Q'\in\msf^{H'}(M'_0)$. For $X\in(\fm_{\wt{Q'}}\cap\fs'_\rs)(F)$, the map
$$ N_{Q'}(F)\ra(\fn_{\wt{Q'}}\cap\fs')(F), n\mapsto \Ad(n^{-1})(X)-X $$
is an isomorphism of $F$-analytic manifolds whose Jacobian is $|D^{\fs'}(X)|_F^{1/2}|D^{\fm_{\wt{Q'}}\cap\fs'}(X)|_F^{-1/2}$. 
\end{lem}

\begin{proof}
See \cite[Lemma 8.1]{MR4350885} for the isomorphism. The computation of its Jacobian is close to the proof of \cite[Proposition 6.3.(ii)]{MR3414387}. 
\end{proof}

Fix a continuous and nontrivial unitary character $\Psi: F\ra\BC^\times$. Let $\langle\cdot,\cdot\rangle$ be the symmetric bilinear form on $\fg'(F)$ defined by
$$ \langle Y,X\rangle:=\Trd(YX),\forall Y,X\in\fg'(F), $$
where $\Trd$ denotes the reduced trace on $\fg'(F)$. It is non-degenerate, which can be seen after the base change to an algebraic closure of $F$. It is also invariant by the adjoint action of $G'(F)$ and $\Ad(\alpha)$. For $f'\in\CC_c^\infty(\fs'(F))$, define its normalised Fourier transform $\hat{f'}\in\CC_c^\infty(\fs'(F))$ by
\begin{equation}\label{fouriertransform2}
 \forall Y\in\fs'(F), \hat{f'}(Y):=c_\Psi(\fs'(F)) \int_{\fs'(F)} f'(X) \Psi(\langle Y,X \rangle) dX, 
\end{equation}
where $c_\Psi(\fs'(F))$ is the unique constant such that $\hat{\hat{f'}}(Y)=f'(-Y)$ for all $f'\in\CC_c^\infty(\fs'(F))$ and all $Y\in\fs'(F)$. For any $M'\in\msl^{H'}(M'_0)$, the restriction of $\langle\cdot,\cdot\rangle$ on $\wt{\fm'}\cap\fs'$ is non-degenerate. Then we can define similarly the normalised Fourier transform of $f'\in\CC_c^\infty((\wt{\fm'}\cap\fs')(F))$. 

Suppose that $P'\in\msf^{H'}(M'_0)$. For $f'\in\CC_c^\infty(\fs'(F))$, we define a function (parabolic descent) $f'_{P'}\in\CC_c^\infty((\fm_{\wt{P'}}\cap\fs')(F))$ by
\begin{equation}\label{pardes2}
 f'_{P'}(Z):=\int_{K_{H'}\times(\fn_{\wt{P'}}\cap\fs')(F)} f'(\Ad(k^{-1})(Z+U)) dUdk 
\end{equation}
for all $Z\in(\fm_{\wt{P'}}\cap\fs')(F)$. By our choices of Haar measures (see \cite[\S I.7]{MR1344131}) 
and the $H'(F)$-invariance of $\langle\cdot,\cdot\rangle$, we show that $(\hat{f'})_{P'}=(f'_{P'})^\string^$, which will be denoted by $\hat{f'}_{P'}$ without confusion. 

\begin{prop}\label{bdoi2}
Let $f'\in\CC_c^\infty(\fs'(F))$. Then
$$ \sup_{Y\in\fs'_\rs(F)}|D^{\fs'}(Y)|_F^{1/2} \int_{H'_Y(F)\bs H'(F)} |f'(\Ad(x^{-1})(Y))| dx<+\infty. $$
\end{prop}

\begin{coro}[cf. Corollary \ref{corII.1}]\label{corII.2}
For any $r\geq0$, the function $|D^{\fs'}(Y)|_F^{-\frac{1}{2}}\sup\{1, -\log |D^{\fs'}(Y)|_F\}^r$ is locally integrable on $\fs'(F)$. 
\end{coro}

The rest of this section is devoted to the proof of Proposition \ref{bdoi2}. We shall follow the main steps in \cite[\S6.3]{MR3414387}, which is similar to the proof of \cite[Theorem 13]{MR0414797}, and only point out some additional ingredients. Let $n$ be the $F$-rank of $G'$. Denote $G'_n:=G', H'_n:=H'$ and $\fs'_n:=\fs'$. Recall that the $F$-rank of $H'_n$ is $n$ in \textbf{Case I} (resp. $\frac{n}{2}$ in \textbf{Case II}). We shall use induction on $n$. For $n=1$ in \textbf{Case I} (resp. $n=2$ in \textbf{Case II}), the proposition is evident since $H'_Y(F)\bs H'(F)$ is compact in our case. 

The following description of semi-simple elements and descendants (see \cite[Definition 7.2.2]{MR2553879}) is a generalisation of \cite[Lemma 2.1]{MR1487565} (see also \cite[Proposition 4.7]{MR3414387}). 

\begin{prop}\label{sselemts2}
1) An element $Y$ of $\fs'(F)$ is semi-simple if and only if it is $H'(F)$-conjugate to an element of the form
$$ Y(B):=\mat(B,0,0,0), $$
with $B\in{\fs'_m}^\times(F)$ being semi-simple with respect to the $H'_m$-action. More precisely, the set of $H'(F)$-conjugacy classes of semi-simple elements in $\fs'(F)$ is bijective to the set of pairs $(m, [B])$ where $0\leq m\leq n$ is an integer in \textbf{Case I} (resp. an even number in \textbf{Case II}) and $[B]$ is a semi-simple $H'_m(F)$-conjugacy class in ${\fs'_m}^\times(F)$. Moreover, $Y(B)$ is regular semi-simple if and only if $m=n$ and $B$ is regular semi-simple in ${\fs'}^\times(F)$. 

2) Let $Y=Y(B)\in\fs'(F)$ be semi-simple. Then the descendant $(H'_Y, \fs'_Y)$ (as a representation) is isomorphic to
$$ (H'_{m,B}, \fs'_{m,B})\times(H'_{n-m},\fs'_{n-m}), $$
where $H'_{m,B}$ (resp. $\fs'_{m,B}$) denotes the centraliser of $B$ in $H'_m$ (resp. $\fs'_m$). 
\end{prop}

\begin{proof}
1) By the base change to an algebraic closure of $F$ containing $E$, we see from \cite[Proposition 2.1]{MR1394521} that an element $Y\in\fs'(F)$ which is $H'(F)$-conjugate to $Y(B)$ in the proposition is semi-simple. Now we suppose that $Y\in\fs'(F)$ is semi-simple. Since $Y^2\in\fh'(F)$, up to $H'(F)$ conjugation, we may suppose that $Y^2=\mat(A,0,0,0)$ with $A\in{\fh'_m}^\times(F)$ being semi-simple in the usual sense. 
From $\mat(A,0,0,0)Y=Y\mat(A,0,0,0)$, we deduce that $Y=\mat(B,0,0,C)$ for some $B\in\fs'_m(F)$ such that $AB=BA$ and some $C\in\fs'_{n-m}(F)$. As $Y^2=\mat(A,0,0,0)$, we have $B\in{\fs'_m}^\times(F)$. Because $Y$ is semi-simple, it is shown in \cite[p. 71]{MR1394521} that $Y$ and $Y^2$ have the same rank over an algebraic closure of $F$ containing $E$. Then $C=0$. We can also see from \cite[Proposition 2.1]{MR1394521} that $B$ is semi-simple with respect to the $H'_m$-action after the base change. We have established the first statement. 

For the second statement, it suffices to notice that two such elements $Y(B_1)$ with $B_1\in{\fs'_{m_1}}^\times(F)$ and $Y(B_2)$ with $B_2\in{\fs'_{m_2}}^\times(F)$ in the proposition are $H'(F)$-conjugate if and only if $m_1=m_2$ (denoted by $m$) and $B_1$ and $B_2$ are $H'_m(F)$-conjugate. 

The third statement follows from the base change or 2). 

2) It can be shown by direct calculation. 
\end{proof}

Fix a Cartan subspace $\fc'$ of $\fs'$. The lemma below is an analogue of \cite[Lemma 29]{MR0414797}. 

\begin{lem}\label{lem29}
Let $f'\in\CC_c^\infty(\fs'(F)-\CN^{\fs'})$. Then
$$ \sup_{Y\in\fc'_\reg(F)}|D^{\fs'}(Y)|_F^{1/2} \int_{H'_Y(F)\bs H'(F)} |f'(\Ad(x^{-1})(Y))| dx<+\infty. $$
\end{lem}

\begin{proof}
We may apply the argument of \cite[Lemma 6.14]{MR3414387} relying on Lemma \ref{lem28} and \cite[Proposition 3.11]{MR3245011}, which is an analogue of Harish-Chandra's semi-simple descent for orbital integrals \cite[Lemma 16.1]{MR2192014}. By Proposition \ref{sselemts2}.2), it suffices to prove the boundedness of orbital integrals for
$$ (H'_{m,B}, \fs'_{m,B})\times(H'_{n-m},\fs'_{n-m}) $$
with $B\in{\fs'_m}^\times(F)$ being semi-simple with respect to the $H'_m$-action and $0<m\leq n$. Since there exists an $H'_{m,B}$-equivariant linear isomorphism $\fs'_{m,B}\ra \fh'_{m,B}$ induced by $Z\mapsto ZB$, the first factor $(H'_{m,B}, \fs'_{m,B})$ is covered by Harish-Chandra's work \cite[Theorem 13]{MR0414797} on classical orbital integrals on Lie algebras. Then we conclude by applying the induction hypothesis to the second factor $(H'_{n-m},\fs'_{n-m}) $. 
\end{proof}

Consider $X_0\in\CN^{\fs'}$. 
By the Jacobson-Morozov theorem for symmetric spaces (Lemma \ref{jm}), there exists a group homomorphism $\varphi: SL_2(F)\ra G'(F)$ such that
$$ X_0=d\varphi\mat(0,1,0,0), Y_0:=d\varphi\mat(0,0,1,0)\in\fs'(F) \text{ and } D_0:=d\varphi\mat(1,0,0,-1)\in \fh'(F). $$
Write $r':=\dim\fs'_{Y_0}$ and $m':=\frac{1}{2}\tr(\ad(-D_0)|_{\fs'_{Y_0}})$. 

\begin{lem}\label{prop4.4}
We have 

1) $r'\geq\frac{1}{2}\sqrt{\dim\fg'}$; 

2) $r'+m'>\frac{1}{4}\dim\fg'+\frac{1}{4}\sqrt{\dim\fg'}$. 
\end{lem}

\begin{proof}
It suffices to check these relations after a base change to an algebraic closure of $F$ containing $E$. Then the lemma is exactly \cite[Proposition 4.4]{MR3414387}. 
\end{proof}

Let $\fs'_\bdd$ be the set of $X\in\fs'(F)$ such that there exists an open neighbourhood $\sigma'$ of $X$ in $\fs'(F)$ satisfying
$$ \sup_{Y\in\fc'_\reg(F)}|D^{\fs'}(Y)|_F^{1/2} \int_{H'_Y(F)\bs H'(F)} |f'(\Ad(x^{-1})(Y))| dx<+\infty $$
for all $f'\in\CC_c^\infty(\fs'(F))$ with $\Supp(f')\subseteq\sigma'$. 
The next lemma is an analogue of \cite[Lemma 38]{MR0414797}. 

\begin{lem}\label{lem38}
We have $\CN^{\fs'}-\{0\}\subseteq \fs'_\bdd$. 
\end{lem}

\begin{proof}
We may apply the argument of \cite[Lemma 6.16]{MR3414387} thanks to Lemma \ref{prop4.4}. We only sketch some steps for convenience and omit the details. In fact, we use induction on $r'$, which is also $\dim\fh'_{X_0}$ by \cite[Lemma 3.2]{MR3414387} and $\dim\fs'=\dim\fh'=\frac{1}{2}\dim\fg'$ (cf. \cite[Lemma 4.3]{MR3414387}), for $X_0\in\CN^{\fs'}-\{0\}$. Suppose that for all $X'_0\in\CN^{\fs'}-\{0\}$ satisfying $\fh'_{X'_0}<r'$, we have  $X'_0\in\fs'_\bdd$. We need to show that $X_0\in\fs'_\bdd$. Denote by $\CN^{\fs'}_q$ the union of all $H'(F)$-orbits in $\CN^{\fs'}$ of dimension $\leq q$. Then $X_0\in\CN^{\fs'}_{\frac{1}{2}\dim\fg'-r'}-\CN^{\fs'}_{\frac{1}{2}\dim\fg'-r'-1}$, and the induction hypothesis means that $\CN^{\fs'}-\CN^{\fs'}_{\frac{1}{2}\dim\fg'-r'}\subseteq\fs'_\bdd$. Fix $t\in F^\times$ with $|t|_F>1$ and then construct an open $H'(F)$-invariant neighbourhood $\omega_0$ of $X_0$ in $\fs'(F)$ as in the proof of \cite[Lemma 6.16]{MR3414387}. Let $Y\in\fc'_\reg(F)$ and $f'\in\CC_c^\infty(\omega_0)$. As in \cite[(1) in p. 1853]{MR3414387} and the first paragraph of the proof of \cite[Lemma 6.16]{MR3414387}, there exists $f'_{\alpha_0}\in\CC_c^\infty(\fs'(F))$ such that $\Supp(f'_{\alpha_0})\cap\CN^{\fs'}_{\frac{1}{2}\dim\fg'-r'}=\emptyset$ and that (see \eqref{oiwochar} for the notation)
\begin{equation}\label{formlem38}
 I_{t^{-1}Y}(f')=|t|_F^{\frac{1}{4}\dim\fg'+\frac{1}{4}\sqrt{\dim\fg'}-r'-m'}(I_Y(f')+I_Y(f'_{\alpha_0})). 
\end{equation}
By Lemma \ref{lem29} and the induction hypothesis, we have 
$$ a:=\sup_{Y\in\fc'_\reg(F)} I_Y(|f'_{\alpha_0}|)<+\infty. $$
Denote $c:=|t|_F^{\frac{1}{4}\dim\fg'+\frac{1}{4}\sqrt{\dim\fg'}-r'-m'}<1$ (by Lemma \ref{prop4.4} and $|t|_F>1$). By iteration, we deduce from \eqref{formlem38} that 
$$ I_{t^{-d}Y}(f')=c^dI_Y(f')+\sum_{1\leq k\leq d}c^k I_{t^{k-d}Y}(f'_{\alpha_0}). $$
Replacing $Y$ with $t^d Y$, we get 
$$ I_{Y}(f')=c^d I_{t^d Y}(f')+\sum_{1\leq k\leq d}c^k I_{t^k Y}(f'_{\alpha_0}). $$
By Lemma \ref{lem28}, we have $\lim\limits_{d\ra+\infty} I_{t^d Y}(f')=0$. Then 
$$ |I_{Y}(f')|\leq a\sum_{k=1}^{+\infty} c^k=\frac{ac}{1-c}. $$
As we may replace $f'$ by $|f'|$ in the above argument, we have completed the proof. 
\end{proof}

\begin{proof}[Proof of Proposition \ref{bdoi2}]
We may use the argument in \cite[\S VI.7]{MR0414797} to show that $0\in\fs'_\bdd$. Then the proposition follows from Lemmas \ref{lem29} and \ref{lem38}. 
\end{proof}


\section{\textbf{Weighted orbital integrals}}\label{sectwoi}

\subsection{The case of $(G,H)$}

Let $E/F$ be a quadratic field extension and $\eta$ the quadratic character of $F^\times/NE^\times$ attached to it, where $N$ is the norm map $E^\times\ra F^\times$. For $x\in H(F)$, which is viewed as an element in $G(F)$, we denote by $\Nrd(x)$ its reduced norm. Suppose that $M\in\msl^{G,\omega}(M_0)$ and that $Q\in\msf^G(M)$. For all $f\in \CC_c^\infty(\fs(F))$ and $X\in (\fm\cap\fs_\rs)(F)$, we define the weighted orbital integral
\begin{equation}\label{defwoi1}
 J_M^Q(\eta, X, f):=|D^\fs(X)|_F^{1/2} \int_{H_X(F)\bs H(F)} f(\Ad(x^{-1})(X)) \eta(\Nrd(x)) v_M^Q(x) dx. 
\end{equation}
Recall that this definition depends implicitly on the choice of a maximal compact subgroup of $G(F)$, which has been fixed in Section \ref{gpandmap}. Since $v_M^Q(x)$ is left-invariant by $M(F)$ and we have $H_X\subseteq M_H$ for $X\in\fm\cap\fs_\rs$, we see that $v_M^Q(x)$ is left-invariant by $H_X(F)$. This integral is absolutely convergent since the orbit $\Ad(H(F))(X)$ is closed in $\fs(F)$, which ensures that the integrand is a compactly supported (and locally constant) function on the homogeneous space. 

Notice that for $x\in M_H(F)$, we have $J_M^Q(\eta, \Ad(x^{-1})(X), f)=\eta(\Nrd(x)) J_M^Q(\eta, X, f)$. 
Sometimes it is convenient to introduce a transfer factor as in \cite[Definition 5.7]{MR3414387}: for $X=\mat(0,A,B,0)\in\fs_\rs(F)$, define 
\begin{equation}\label{eqdeftf} 
\kappa(X):=\eta(\Nrd(A)), 
\end{equation}
where $\Nrd(A)$ denotes the reduced norm of $A\in GL_n(D)$. Then we have $\kappa(\Ad(x^{-1})(X))=\eta(\Nrd(x))\kappa(X)$, and thus the function $\kappa(\cdot)J_M^Q(\eta, \cdot, f)$ is constant on $\Ad(M_H(F))(X)$. 

Though we mainly consider $M\in\msl^{G,\omega}(M_0)$, it is unharmful to extend our definition by \eqref{defwoi1} to all Levi subgroups of the form $M=\Ad(w)(L)$, where $L\in\msl^{G,\omega}(M_0)$ and $w\in W_0^H$. 

One may also extend in the obvious way the definition (\ref{defwoi1}) of weighted orbital integrals to the symmetric pair $(L, L_H, \Ad(\epsilon))$, where $L\in\msl^{G,\omega}(M_0)$, since it appears as the product of some copies of the form $(G, H, \Ad(\epsilon))$ in lower dimensions. 

\begin{prop}\label{propwoi1}
Suppose that $M\in\msl^{G,\omega}(M_0)$ and that $Q\in\msf^G(M)$. 

1) For $X\in(\fm\cap\fs_\rs)(F)$ fixed, the support of the distribution $J_M^Q(\eta, X, \cdot)$ is contained in the closed orbit $\Ad(H(F))(X)$. 

2) For $f\in \CC_c^\infty(\fs(F))$ fixed, the function $J_M^Q(\eta, \cdot, f)$ is locally constant  on $(\fm\cap\fs_\rs)(F)$. If $\fc\subseteq\fm\cap\fs$ is a Cartan subspace, the restriction of this function to $\fc_\reg(F)$ vanishes outside a compact subset of $\fc(F)$. 

3) If $w\in \Norm_{H(F)}(M_0)$, $x\in M_H(F)$ and $k\in K_H$, we have the equality
$$ J_{\Ad(w)(M)}^G(\eta, \Ad(wx)(X), \Ad(k)(f))=\eta(\Nrd(wxk)) J_M^G(\eta, X, f) $$
for all $X\in(\fm\cap\fs_\rs)(F)$ and $f\in\CC_c^\infty(\fs(F))$. 

4) For $X\in(\fm\cap\fs_\rs)(F)$ and $f\in \CC_c^\infty(\fs(F))$, we have the equality
$$ J_M^Q(\eta, X, f)=J_M^{M_Q}(\eta, X, f_Q^\eta), $$
where $f_Q^\eta\in\CC_c^\infty((\fm_Q\cap\fs)(F))$ is defined by (\ref{pardes1}). 

5) (Descent formula) If $L\in\msl^{G,\omega}(M_0)$, $L\subseteq M$ and $X\in(\fl\cap\fs_\rs)(F)$, we have
$$ J_M^G(\eta, X, f)=\sum_{L'\in\msl^G(L)} d_L^G(M, L') J_L^{L'}(\eta, X, f_{Q'}^\eta) $$
for all $f\in\CC_c^\infty(\fs(F))$, where $Q'$ denotes the second component of $s(M,L')$ (see Section \ref{gmmaps}). 

6) (Non-equivariance) For $X\in(\fm\cap\fs_\rs)(F)$, $y\in H(F)$ and $f\in\CC_c^\infty(\fs(F))$, we have the equality
$$ J_M^G(\eta, X, \Ad(y^{-1})(f))=\eta(\Nrd(y)) \sum_{Q\in\msf^G(M)} J_M^{M_Q}(\eta, X, f_{Q,y}^\eta), $$
where $f_{Q,y}^\eta\in\CC_c^\infty((\fm_Q\cap\fs)(F))$ is defined by
\begin{equation}\label{nonequcst1}
 f_{Q,y}^\eta(Z):=\int_{K_H\times(\fn_Q\cap\fs)(F)} f(\Ad(k^{-1})(Z+U)) \eta(\Nrd(k)) v'_Q(ky) dUdk, \forall Z\in(\fm_Q\cap\fs)(F). 
\end{equation}
\end{prop}

\begin{proof}
1) This is obvious from the definition. 

2) Let $Y\in(\fm\cap\fs_\rs)(F)$. Then $\fc:=\fs_Y\subseteq\fm\cap\fs$ is a Cartan subspace and $Y\in\fc_\reg(F)$. 
Since $\Ad(M_H(F))(\fc_\reg(F))$ is an open subset of $(\fm\cap\fs_\rs)(F)$ (see \cite[p. 105]{MR1375304}), in order to prove the first statement, it suffices to find a neighbourhood $U$ of $Y$ in $\fc_\reg(F)$ on which the function $\kappa(\cdot)J_M^Q(\eta, \cdot, f)$ is constant. We shall follow the proof of \cite[Theorem 17.11]{MR2192014}. Consider the function $\phi$ on $\fc_\reg(F)\times(T_\fc(F)\bs H(F))$ defined by $\phi(X,x):=(\kappa f)(\Ad(x^{-1})(X))$. Then $\phi$ is locally constant but usually not compactly supported. However, now choosing a compact neighbourhood $\sigma_\fc$ of $Y$ in $\fc_\reg(F)$, we see from Harish-Chandra's compactness lemma for symmetric spaces (Lemma \ref{lem25}) applied to $\sigma_\fs:=\Supp(f)$ that the restriction of $\phi$ to $\sigma_\fc\times(T_\fc(F)\bs H(F))$ is compactly supported. By \cite[Lemma 2.1]{MR2192014}, there exists an open neighbourhood $U$ of $Y$ in $\sigma_\fc$ such that $\phi(X,x)=\phi(Y,x)$ for all $X\in U$ and $x\in T_\fc(F)\bs H(F)$. It follows that the function $\kappa(\cdot)J_M^Q(\eta, \cdot, f)$ is constant on $U$. 

The second statement is a corollary of Lemma \ref{lem28}. 

3) The effect of $\Ad(w)$ is a consequence of our choice of Haar measures. The effect of $\Ad(x)$ results from the left-invariance of $v_M^G(x)$ by $M_H(F)$. The effect of $\Ad(k)$ comes from the right-invariance of $v_M^G(x)$ by $K_H$. One should keep in mind the effect of $\eta(\Nrd(x))$ in every step. 

4) Write $Q_H:=Q\cap H\in\msf^H(M_0)$. One sees that $M_{Q_H}=M_Q\cap H$ and that $N_{Q_H}=N_Q\cap H$. Applying the change of variables $x=mnk$ with $m\in M_{Q_H}(F)$, $n\in N_{Q_H}(F)$ and $k\in K_H$ in (\ref{defwoi1}), since $v_M^Q(x)=v_M^{M_Q}(m)$, we have
$$ J_M^Q(\eta, X, f)=|D^\fs(X)|_F^{1/2} \int_{(M_{Q_H,X}(F)\bs M_{Q_H}(F))\times N_{Q_H}(F)\times K_H} f(\Ad(mnk)^{-1}(X)) \eta(\Nrd(mk)) v_M^{M_Q}(m) dk dn dm. $$

Applying Lemma \ref{jac1} to $Y=\Ad(m^{-1})(X)$, we deduce
\[\begin{split} 
J_M^Q(\eta, X, f)=&|D^{\fm_Q\cap\fs}(X)|_F^{1/2} \int_{(M_{Q_H,X}(F)\bs M_{Q_H}(F))\times K_H\times (\fn_Q\cap\fs)(F)} f(\Ad(k^{-1})(\Ad(m^{-1})(X)+U)) \\ 
&\eta(\Nrd(mk)) v_M^{M_Q}(m) dU dk dm \\
=&|D^{\fm_Q\cap\fs}(X)|_F^{1/2} \int_{M_{Q_H,X}(F)\bs M_{Q_H}(F)} f_Q^\eta(\Ad(m^{-1})(X)) \eta(\Nrd(m)) v_M^{M_Q}(m) dm \\
=&J_M^{M_Q}(\eta, X, f_Q^\eta). 
\end{split}\]

5) It follows from (7) in Section \ref{gmmaps} and 4). 

6) By the change of variables, we see that
$$ J_M^G(\eta, X, \Ad(y^{-1})(f))=|D^\fs(X)|_F^{1/2} \int_{H_X(F)\bs H(F)} f(\Ad(x^{-1})(X)) \eta(\Nrd(xy)) v_M(xy) dx.  $$
For $x\in H(F)$ and $Q\in\msf^{G,\omega}(M_0)$, let $k_Q(x)$ be an element in $K_H$ such that $xk_Q(x)^{-1}\in Q_H(F)$. It follows from the product formula (\ref{split6.3}) that (see the proof of \cite[Lemma 8.2]{MR625344})
$$ v_M(xy)=\sum_{Q\in\msf^G(M)} v_M^Q(x) v'_Q(k_Q(x)y).  $$
As in 4), we write
\[\begin{split}
J_M^G(\eta, X, \Ad(y^{-1})(f))=&\eta(\Nrd(y)) \sum_{Q\in\msf^G(M)} |D^\fs(X)|_F^{1/2} \int_{H_X(F)\bs H(F)} f(\Ad(x^{-1})(X)) \eta(\Nrd(x)) v_M^Q(x) \\ 
&v'_Q(k_Q(x)y) dx \\
=&\eta(\Nrd(y)) \sum_{Q\in\msf^G(M)} |D^\fs(X)|_F^{1/2} \int_{(M_{Q_H,X}(F)\bs M_{Q_H}(F))\times N_{Q_H}(F)\times K_H} \\ 
&f(\Ad(mnk)^{-1}(X)) \eta(\Nrd(mk)) v_M^{M_Q}(m) v'_Q(ky) dk dn dm. 
\end{split}\]

Applying again Lemma \ref{jac1} to $Y=\Ad(m^{-1})(X)$, we obtain
\[\begin{split}
J_M^G(\eta, X, \Ad(y^{-1})(f))=&\eta(\Nrd(y)) \sum_{Q\in\msf^G(M)} |D^{\fm_Q\cap\fs}(X)|_F^{1/2} \int_{(M_{Q_H,X}(F)\bs M_{Q_H}(F))\times K_H\times (\fn_Q\cap\fs)(F)} \\ 
&f(\Ad(k^{-1})(\Ad(m^{-1})(X)+U)) \eta(\Nrd(mk)) v_M^{M_Q}(m) v'_Q(ky) dU dk dm \\ 
=&\eta(\Nrd(y)) \sum_{Q\in\msf^G(M)} |D^{\fm_Q\cap\fs}(X)|_F^{1/2} \int_{M_{Q_H,X}(F)\bs M_{Q_H}(F)} f_{Q,y}^\eta(\Ad(m^{-1})(X)) \\ 
&\eta(\Nrd(m)) v_M^{M_Q}(m) dm \\ 
=&\eta(\Nrd(y)) \sum_{Q\in\msf^G(M)} J_M^{M_Q}(\eta, X, f_{Q,y}^\eta).
\end{split}\]
\end{proof}

\begin{lem}\label{lemIII.5.1}
Suppose that $M\in\msl^{G,\omega}(M_0)$ and that $Q\in\msf^G(M)$. Let $\sigma\subseteq\fs(F)$ be a compact subset. There exist constants $c>0$ and $N\in\BN$ such that if $x\in H(F)$ and $X\in(\fm\cap\fs_\rs)(F)$ satisfy $\Ad(x^{-1})(X)\in\sigma$, then
$$ |v_M^Q(x)|\leq c\sup\{1, -\log |D^\fs(X)|_F\}^N. $$
\end{lem}

\begin{proof}
It is shown in the proof of \cite[Lemme III.5]{MR1344131} that there exist constants $c_1>0$ and $N\in\BN$ such that for all $x\in G(F)$, 
$$ |v_M^Q(x)|\leq c_1(1+\log \|x\|)^N. $$

Suppose that $x\in H(F)$ and $X\in(\fm\cap\fs_\rs)(F)$ satisfy $\Ad(x^{-1})(X)\in\sigma$. If we replace $x$ by $yx$ and $X$ by $\Ad(y)(X)$, where $y\in M_H(F)$, the two sides in the inequality to be proved remain unchanged. Since $\mst(\fm\cap\fs)$ is a finite set, we may fix a Cartan subspace $\fc\subseteq\fm\cap\fs$ and suppose that $X\in\fc_\reg(F)$. Let $\tau\in T_\fc(F)$ be such that
$$ \|\tau x\|=\inf_{\tau'\in T_\fc(F)} \|\tau' x\|. $$
Then
$$ |v_M^Q(x)|= |v_M^Q(\tau x)|\leq c_1(1+\log \|\tau x\|)^N=c_1(1+\inf_{\tau'\in T_\fc(F)}\log\|\tau' x\|)^N. $$
Now it suffices to apply Lemma \ref{lemIII.4}. 
\end{proof}

\begin{coro}\label{corIII.6.1}
Suppose that $M\in\msl^{G,\omega}(M_0)$ and that $Q\in\msf^G(M)$. Let $f\in\CC_c^\infty(\fs(F))$. There exist constants $c>0$ and $N\in\BN$ such that for all $X\in(\fm\cap\fs_\rs)(F)$, we have
$$ |J_M^{M_Q}(\eta, X, f_Q^\eta)| \leq c \sup\{1, -\log |D^\fs(X)|_F\}^N. $$
\end{coro}

\begin{proof}
By Proposition \ref{propwoi1}.4) and Lemma \ref{lemIII.5.1} applied to $\sigma=\Supp(f)$, we see that
\[\begin{split}
|J_M^{M_Q}(\eta, X, f_Q^\eta)|\leq&|D^\fs(X)|_F^{1/2} \int_{H_X(F)\bs H(F)} |f(\Ad(x^{-1})(X)) v_M^Q(x)| dx \\
\leq&c\sup\{1, -\log |D^\fs(X)|_F\}^N |D^\fs(X)|_F^{1/2} \int_{H_X(F)\bs H(F)} |f(\Ad(x^{-1})(X))| dx. 
\end{split}\]
Now we draw our conclusion by Proposition \ref{bdoi1}. 
\end{proof}

\subsection{The case of $(G',H')$}

Suppose that $M'\in\msl^{H'}(M'_0)$ and that $Q'\in\msf^{H'}(M')$. For all $f'\in \CC_c^\infty(\fs'(F))$ and $Y\in(\wt{\fm'}\cap\fs'_\rs)(F)$, we define the weighted orbital integral
\begin{equation}\label{defwoi2}
 J_{M'}^{Q'}(Y, f'):=|D^{\fs'}(Y)|_F^{1/2} \int_{H'_Y(F)\bs H'(F)} f'(\Ad(x^{-1})(Y)) v_{M'}^{Q'}(x) dx. 
\end{equation}
Recall that this definition depends implicitly on the choice of a maximal compact subgroup of $H'(F)$, which has been fixed in Section \ref{gpandmap}. By the base change to an algebraic closure of $F$ containing $E$, we see that $H'_Y\subseteq H'_{Y^2}\subseteq M'$ for $Y\in\wt{\fm'}\cap\fs'_\rs$. Then $v_{M'}^{Q'}(x)$ is left-invariant by $H'_Y(F)$. This integral is absolutely convergent since the orbit $\Ad(H'(F))(Y)$ is closed in $\fs'(F)$. 
Notice that for $x\in M'(F)$, we have $J_{M'}^{Q'}(\Ad(x^{-1})(Y), f')=J_{M'}^{Q'}(Y, f')$, i.e., the function $J_{M'}^{Q'}(\cdot, f')$ is constant on $\Ad(M'(F))(Y)$. 
One may extend in the obvious way the definition (\ref{defwoi2}) to the symmetric pair $(\wt{L'}, L', \Ad(\alpha))$, where $L'\in\msl^{H'}(M'_0)$, since it appears as the product of some copies of the form $(G', H', \Ad(\alpha))$ in lower dimensions. 

\begin{prop}\label{propwoi2}
Suppose that $M'\in\msl^{H'}(M'_0)$ and that $Q'\in\msf^{H'}(M')$. 

1) For $Y\in(\wt{\fm'}\cap\fs'_\rs)(F)$ fixed, the support of the distribution $J_{M'}^{Q'}(Y, \cdot)$ is contained in the closed orbit $\Ad(H'(F))(Y)$. 

2) For $f'\in \CC_c^\infty(\fs'(F))$ fixed, the function $J_{M'}^{Q'}(\cdot, f')$ is locally constant  on $(\wt{\fm'}\cap\fs'_\rs)(F)$. If $\fc'\subseteq\wt{\fm'}\cap\fs'$ is a Cartan subspace, the restriction of this function to $\fc'_\reg(F)$ vanishes outside a compact subset of $\fc'(F)$. 

3) If $w\in \Norm_{H'(F)}(M'_0)$, $x\in M'(F)$ and $k\in K_{H'}$, we have the equality
$$ J_{\Ad(w)(M')}^{H'}(\Ad(wx)(Y), \Ad(k)(f'))=J_{M'}^{H'}(Y, f') $$
for all $Y\in(\wt{\fm'}\cap\fs'_\rs)(F)$ and $f'\in\CC_c^\infty(\fs'(F))$. 

4) For $Y\in(\wt{\fm'}\cap\fs'_\rs)(F)$ and $f'\in\CC_c^\infty(\fs'(F))$, we have the equality
$$ J_{M'}^{Q'}(Y, f')=J_{M'}^{M_{Q'}}(Y, f'_{Q'}), $$
where $f'_{Q'}\in\CC_c^\infty((\fm_{\wt{Q'}}\cap\fs')(F))$ is defined by (\ref{pardes2}). 

5) (Descent formula) If $L'\in\msl^{H'}(M'_0)$, $L'\subseteq M'$ and $Y\in(\wt{\fl'}\cap\fs'_\rs)(F)$, we have
$$ J_{M'}^{H'}(Y, f')=\sum_{L\in\msl^{H'}(L')} d_{L'}^{H'}(M', L) J_{L'}^L(Y, f'_Q) $$
for all $f'\in\CC_c^\infty(\fs'(F))$, where $Q$ denotes the second component of $s(M',L)$ (see Section \ref{gmmaps}). 

6) (Noninvariance) For $Y\in(\wt{\fm'}\cap\fs'_\rs)(F)$, $y\in H'(F)$ and $f'\in\CC_c^\infty(\fs'(F))$, we have the equality
$$ J_{M'}^{H'}(Y, \Ad(y^{-1})(f'))=\sum_{Q'\in\msf^{H'}(M')} J_{M'}^{M_{Q'}}(Y, f'_{Q',y}), $$
where $f'_{Q',y}\in\CC_c^\infty((\fm_{\wt{Q'}}\cap\fs')(F))$ is defined by
\begin{equation}\label{nonequcst2}
 f'_{Q',y}(Z):=\int_{K_{H'}\times(\fn_{\wt{Q'}}\cap\fs')(F)} f'(\Ad(k^{-1})(Z+U)) v'_{Q'}(ky) dUdk, \forall Z\in(\fm_{\wt{Q'}}\cap\fs')(F). 
\end{equation}
\end{prop}

\begin{proof}
It is almost the same as the proof of Proposition \ref{propwoi1}, except that one needs to use Lemma \ref{jac2} to show 4) and 6). 
\end{proof}

\begin{lem}[cf. Lemma \ref{lemIII.5.1}]\label{lemIII.5.2}
Suppose that $M'\in\msl^{H'}(M'_0)$ and that $Q'\in\msf^{H'}(M')$. Let $\sigma'\subseteq\fs'(F)$ be a compact subset. There exist constants $c>0$ and $N\in\BN$ such that if $x\in H'(F)$ and $Y\in(\wt{\fm'}\cap\fs'_\rs)(F)$ satisfy $\Ad(x^{-1})(Y)\in\sigma'$, then
$$ |v_{M'}^{Q'}(x)|\leq c\sup\{1, -\log |D^{\fs'}(Y)|_F\}^N. $$
\end{lem}

\begin{coro}\label{corIII.6.2}
Suppose that $M'\in\msl^{H'}(M'_0)$ and that $Q'\in\msf^{H'}(M')$. Let $f'\in\CC_c^\infty(\fs'(F))$. There exist constants $c>0$ and $N\in\BN$ such that for all $Y\in(\wt{\fm'}\cap\fs'_\rs)(F)$, we have
$$ |J_{M'}^{M_{Q'}}(Y, f'_{Q'})| \leq c \sup\{1, -\log |D^{\fs'}(Y)|_F\}^N. $$
\end{coro}

\begin{proof}
We may apply the argument of Corollary \ref{corIII.6.1} with the help of Proposition \ref{propwoi2}.4), Lemma \ref{lemIII.5.2} and Proposition \ref{bdoi2}. 
\end{proof}


\section{\textbf{The noninvariant trace formula}}\label{secnontf}

\subsection{The case of $(G,H)$}

Suppose that $M\in\msl^{G,\omega}(M_0)$. For $x,y\in G(F)$, we define a $(G,M)$-family $(v_P(x,y))_{P\in\msp^G(M)}$ as in \cite[(12.1) in \S12]{MR1114210} by
$$ v_P(\lambda,x,y):=e^{-\lambda(H_P(y)-H_{\ov{P}}(x))}, \forall \lambda\in i\fa_M^\ast, P\in\msp^G(M), $$
where $\ov{P}\in\msp^G(M)$ is the parabolic subgroup opposite to $P$. Let $E/F$ be a quadratic field extension and $\eta$ the quadratic character of $F^\times/NE^\times$ attached to it. For $f, f'\in\CC_c^\infty(\fs(F))$ and $X\in(\fm\cap\fs_\rs)(F)_\el$, we define
\begin{equation}\label{defV.1.1}
 J_M^G(\eta, X, f, f'):=|D^\fs(X)|_F \int_{(A_M(F)\bs H(F))^2} f(\Ad(x^{-1})(X)) f'(\Ad(y^{-1})(X)) \eta(\Nrd(x^{-1}y)) v_M(x, y) dx dy. 
\end{equation}

\begin{prop}\label{propV.1.1}
Suppose that $M\in\msl^{G,\omega}(M_0)$ and that $f,f'\in\CC_c^\infty(\fs(F))$. 

1) The integral (\ref{defV.1.1}) is absolutely convergent. 

2) The function $J_M^G(\eta, \cdot, f, f')$ is locally constant on $(\fm\cap\fs_\rs)(F)_\el$. 

3) If $\fc\subseteq\fm\cap\fs$ is an $M$-elliptic Cartan subspace, the restriction of $J_M^G(\eta, \cdot, f, f')$ to $\fc_\reg(F)$ vanishes outside a compact subset of $\fc(F)$. 

4) If $w\in \Norm_{H(F)}(M_0)$, $x\in M_H(F)$ and $k, k'\in K_H$, we have the equality
$$ J_{\Ad(w)(M)}^G(\eta, \Ad(wx)(X), \Ad(k)(f), \Ad(k')(f'))=\eta(\Nrd(kk')) J_M^G(\eta, X, f, f') $$
for all $X\in(\fm\cap\fs_\rs)(F)_\el$. 

5) There exist constants $c>0$ and $N\in\BN$ such that for all $X\in(\fm\cap\fs_\rs)(F)_\el$, we have
$$ |J_M^G(\eta, X, f, f')| \leq c \sup\{1, -\log |D^\fs(X)|_F\}^N. $$

6) For all $X\in(\fm\cap\fs_\rs)(F)_\el$, we have
$$ J_M^G(\eta, X, f, f')=\sum_{L_1,L_2\in\msl^G(M)} d_M^G(L_1,L_2) J_M^{L_1}(\eta,X,f_{\ov{Q_1}}^\eta) J_M^{L_2}(\eta,X,{f'}_{Q_2}^\eta), $$
where $(Q_1,Q_2):=s(L_1,L_2)$ (see Section \ref{gmmaps}). 
\end{prop}

\begin{proof}
The statements 1)-4) can be proved in the same way as the proof of analogous properties for (\ref{defwoi1}) in Section \ref{sectwoi}. Notice that the $\eta(\Nrd(\cdot))$-invariant effects coming from $x$ and $y$ may sometimes cancel. 

For $x\in G(F)$, we define a $(G,M)$-family $(\ov{v}_P(x))_{P\in\msp^G(M)}$ by
$$ \ov{v}_P(\lambda, x):=e^{\lambda(H_{\ov{P}}(x))},\forall\lambda\in i\fa_M^\ast,P\in\msp^G(M). $$
Then $v_{P}(x,y)=\ov{v}_P(x)v_P(y)$ as the product of $(G,M)$-families. Notice that for all $Q\in\msp^G(M)$ and $x\in G(F)$, we have
$$ \ov{v}_M^Q(x)=v_M^{\ov{Q}}(x). $$

The statement 6) is a consequence of the splitting formula of $(G,M)$-families ((6) in Section \ref{gmmaps}) and Proposition \ref{propwoi1}.4). It together with Corollary \ref{corIII.6.1} implies the statement 5). 
\end{proof}

For $f,f'\in\CC_c^\infty(\fs(F))$, we define
\begin{equation}\label{geomnoninvltf1}
\begin{split}
 J^G(\eta, f, f'):=&\sum_{M\in\msl^{G,\omega}(M_0)} |W_0^{M_n}| |W_0^{GL_n}|^{-1} (-1)^{\dim(A_M/A_G)} \sum_{\fc\in\mst_\el(\fm\cap\fs)} |W(M_H, \fc)|^{-1} \int_{\fc_\reg(F)} \\ 
&J_M^G(\eta, X, f, f') dX. 
\end{split}
\end{equation}
This expression is absolutely convergent by Proposition \ref{propV.1.1}.5) and Corollary \ref{cor20.2}. 

\begin{remark}\label{rmkV.1.1}
We have the equality
$$ J^G(\eta, f, f')=J^G(\eta, f', f). $$
It results from the fact that for all $M\in\msl^{G,\omega}(M_0)$ and all $x,y\in G(F)$, we have $v_M(x,y)=v_M(y,x)$. 
\end{remark}

Again, one may extend in the obvious way the definitions (\ref{defV.1.1}) and (\ref{geomnoninvltf1}) to the symmetric pair $(L, L_H, \Ad(\epsilon))$, where $L\in\msl^{G,\omega}(M_0)$, since it appears as the product of some copies of the form $(G, H, \Ad(\epsilon))$ in lower dimensions. 

\begin{thm}[Noninvariant trace formula]\label{noninvltf1}
For all $f, f'\in\CC_c^\infty(\fs(F))$, we have the equality
$$ J^G(\eta, f, \hat{f'})=J^G(\eta, \hat{f}, f'). $$
\end{thm}

The rest of this section is devoted to the proof of Theorem \ref{noninvltf1}. 

Fix $P_0\in\msp^G(M_0)$. Denote
$$ \ov{\fa_{P_0}^+}:=\{T\in\fa_{M_0}: \alpha(T)\geq0, \forall\alpha\in\Delta_{P_0}^G\}. $$
For $T\in\ov{\fa_{P_0}^+}$, write
$$ d(T):=\inf_{\alpha\in\Delta_{P_0}^G}\alpha(T), $$
which is invariant under the translation by $\fa_G$. Set $R_0:=(\fa_{M_0,F}+\fa_G)/\fa_G$, which is a lattice in $\fa_{M_0}/\fa_G$. For $T\in R_0\cap(\ov{\fa_{P_0}^+}/\fa_G)$, we define a function $u(\cdot,T)$ on $A_G(F)\bs G(F)$ as in \cite[p. 21]{MR1114210}, which is the characteristic function of certain compact subset. To be precise, let $C_{M_0}(T)$ be the convex hull in $\fa_{M_0}/\fa_G$ of 
$$ \{T_B: B\in\msp^G(M_0)\}, $$
where $T_B$ denotes the unique $W_0^G$-translate of $T$ which lies in $\ov{\fa_B^+}$. Then $u(x,T)$ is defined as the characteristic function of the set of points 
$$ x=k_1mk_2, m\in A_G(F)\bs M_0(F), k_1,k_2\in K_G, $$
in $A_G(F)\bs G(F)$ such that $H_{M_0}(m)$ lies in $C_{M_0}(T)$. 

Let $f, f'\in\CC_c^\infty(\fs(F))$. For $x\in H(F)$, we define
$$ k(x, f, f'):=\int_{\fs(F)} f(X) f'(\Ad(x^{-1})(X)) dX. $$
For $T\in R_0\cap(\ov{\fa_{P_0}^+}/\fa_G)$, we define 
$$ K^T(\eta, f, f'):=\int_{A_G(F)\bs H(F)} k(x, f, f') \eta(\Nrd(x)) u(x, T) dx. $$
Since $A_G(F)\bs H(F)$ is a closed subgroup of $A_G(F)\bs G(F)$, the restriction of $u(x, T)$ to $A_G(F)\bs H(F)$ is also compactly supported, and the above integral is absolutely convergent. 

By the Weyl integration formula (Proposition \ref{wif2.1}), we obtain the geometric expansion
$$ K^T(\eta, f, f')=\sum_{M\in\msl^{G,\omega}(M_0)} |W_0^{M_n}| |W_0^{GL_n}|^{-1} \sum_{\fc\in\mst_\el(\fm\cap\fs)} |W(M_H, \fc)|^{-1} \int_{\fc_\reg(F)} K^T(\eta, X, f, f') dX, $$
where
\[\begin{split}
 K^T(\eta, X, f, f'):=&|D^\fs(X)|_F \int_{A_G(F)\bs H(F)} \int_{A_M(F)\bs H(F)} f(\Ad(x^{-1})(X)) f'(\Ad(xy)^{-1} (X)) \eta(\Nrd(y)) \\ 
&u(y, T) dx dy. 
\end{split}\]
By the change of variables $xy\mapsto y$, we can write
\begin{equation}\label{firstwoi}
 K^T(\eta, X, f, f')=|D^\fs(X)|_F \int_{(A_M(F)\bs H(F))^2} f(\Ad(x^{-1})(X)) f'(\Ad(y^{-1})(X)) \eta(\Nrd(x^{-1}y)) u_M(x, y, T) dx dy, 
\end{equation}
where
$$ u_M(x, y, T):=\int_{A_G(F)\bs A_M(F)} u(x^{-1}ay, T) da $$
is defined as in \cite[p. 21]{MR1114210}. 

For $x,y\in G(F)$ and $T\in R_0\cap(\ov{\fa_{P_0}^+}/\fa_G)$, we define the second weight function $v_M(x, y, T)$ as in \cite[p. 30]{MR1114210}, which is left-invariant under the multiplication of $A_M(F)$ on $x$ or $y$. To be precise, let $\lambda\in\fa_M^\ast\otimes_\BR\BC$ be a point whose real part $\lambda_\BR\in\fa_M^\ast$ is in general position. For $P\in\msp^G(M)$, set 
$$ \Delta_P^\lambda:=\{\alpha\in\Delta_P^G: \lambda_\BR(\alpha^\vee)<0\}, $$
where $\alpha^\vee$ is the ``coroot'' associated to $\alpha\in\Delta_P^G$ (see \cite[p. 26]{MR2192011}). Denote by $\varphi_P^\lambda$ the characteristic function of the set of $T'\in\fa_M$ such that $\varpi_\alpha(T')>0$ for each $\alpha\in\Delta_P^\lambda$ and that $\varpi_\alpha(T')\leq0$ for each $\alpha\in\Delta_P^G-\Delta_P^\lambda$, where $\{\varpi_\alpha: \alpha\in\Delta_P^G\}$ is the basis of $(\fa_P^G)^\ast$ which is dual to $\{\alpha^\vee: \alpha\in\Delta_P^G\}$. Let 
$$ Y_P(x,y,T):=T_P+H_P(x)-H_{\ov{P}}(y), \forall P\in\msp^G(M). $$
The set $\msy_M(x,y,T)=\{Y_P(x,y,T): P\in\msp^G(M)\}$ is a $(G,M)$-orthogonal set in the sense of \cite[p. 19]{MR1114210}. Define 
$$ \sigma_M(T', \msy_M(x,y,T)):=\sum_{P\in\msp^G(M)} (-1)^{|\Delta_P^\lambda|} \varphi_P^\lambda(T'-Y_P(x,y,T)), \forall T'\in\fa_M/\fa_G. $$
The function $\sigma_M(\cdot, \msy_M(x,y,T))$ is known to be compactly supported (see \cite[p. 22]{MR1114210}). Then $v_M(x, y, T)$ is defined as the integral 
$$ v_M(x, y, T):=\int_{A_G(F)\bs A_M(F)} \sigma_M(H_M(a), \msy_M(x,y,T))da. $$
Now, we define the corresponding weighted orbital integral 
\begin{equation}\label{secondwoi}
 J^T(\eta, X, f, f'):=|D^\fs(X)|_F \int_{(A_M(F)\bs H(F))^2} f(\Ad(x^{-1})(X)) f'(\Ad(y^{-1})(X)) \eta(\Nrd(x^{-1}y)) v_M(x, y, T) dx dy. 
\end{equation}

Let $\fc\subseteq\fm\cap\fs$ be an $M$-elliptic Cartan subspace. For $\varepsilon>0$ and $T\in\fa_{M_0,F}\cap\ov{\fa_{P_0}^+}$ with large $\|T\|$, consider the domain near the singular set
$$ \fc(\varepsilon, T):=\{X\in\fc_\reg(F): |D^\fs(X)|\leq e^{-\varepsilon\| T\|}\}. $$

\begin{lem}\label{lem43}
Fix an arbitary constant $\varepsilon_0>0$. Fix a constant $\varepsilon'>0$ satisfying the condition of Lemma \ref{lem44}. Let $\fc\subseteq\fm\cap\fs$ be an $M$-elliptic Cartan subspace. Given $\varepsilon>0$, there exists $c>0$ such that for any $T\in\fa_{M_0,F}\cap\ov{\fa_{P_0}^+}$ with $\|T\|\geq\varepsilon_0$, 
$$ \int_{\fc(\varepsilon, T)} (|K^T(\eta, X, f, f')|+|J^T(\eta, X, f, f')|) dX \leq c e^{-\frac{\varepsilon'\varepsilon\|T\|}{2}}. $$
\end{lem}

\begin{proof}
It is shown in \cite[(4.8) in p. 31]{MR1114210} that there exist $c_1,d_1>0$ such that for all $x,y\in G(F)$ and $T\in\fa_{M_0,F}\cap\ov{\fa_{P_0}^+}$ with $\|T\|\geq\varepsilon_0$, 
$$ u_M(x, y, T) \leq c_1 (\| T\| + \log \| x\| + \log \| y\|)^{d_1}. $$

For any $a_1,a_2\in A_M(F)$, we deduce that
$$ u_M(x, y, T)=u_M(a_1 x, a_2 y, T)\leq c_1 (\|T\| + \log\|a_1 x\| + \log\|a_2 y\|)^{d_1}. $$
Since $T_\fc(F)/A_M(F)$ is compact, there exists $c_2>0$ such that
$$ u_M(x, y, T)\leq c_2 (\|T\| + \inf_{\tau_1\in T_\fc(F)} \log\|\tau_1 x\| + \inf_{\tau_2\in T_\fc(F)} \log\|\tau_2 y\|)^{d_1}. $$

Now let $x,y\in H(F)$ and $X\in\fc(\varepsilon,T)$, and assume that
$$ f(\Ad(x^{-1})(X))f'(\Ad(y^{-1})(X))\neq0. $$
Let $\sigma\subseteq\fs(F)$ be a compact subset containing $\Supp(f)\cup\Supp(f')$. From Lemma \ref{lemIII.4}, there exists $c_\sigma>0$ such that
$$ \inf_{\tau_1\in T_\fc(F)} \log\|\tau_1 x\|, \inf_{\tau_2\in T_\fc(F)} \log\|\tau_2 y\|\leq c_\sigma\sup\{1,-\log|D^\fs(X)|_F\}. $$
Therefore, there exists $c'_\sigma>0$ such that
\begin{equation}\label{est1}
 u_M(x, y, T) \leq c'_\sigma (\| T\| - \log |D^\fs(X)|_F)^{d_1}. 
\end{equation}

By Proposition \ref{bdoi1}, there exists $c_3>0$ such that
\begin{equation}\label{est2}
 |D^\fs(X)|_F^{1/2}\int_{A_M(F)\bs H(F)} |f(x^{-1}Xx)| dx \leq c_3 
\end{equation}
and 
\begin{equation}\label{est3}
 |D^\fs(X)|_F^{1/2}\int_{A_M(F)\bs H(F)} |f'(y^{-1}Xy)| dy \leq c_3 
\end{equation}
for all $X\in\fc_\reg(F)$. 

Putting the estimates (\ref{est1}), (\ref{est2}) and (\ref{est3}) into the definition (\ref{firstwoi}) of $K^T(\eta, X, f, f')$, we obtain the inequality
$$ |K^T(\eta, X, f, f')| \leq c'_\sigma c_3^2 (\| T\| - \log |D^\fs(X)|_F)^{d_1}. $$

By Lemma \ref{lem44}, for any subset $B$ of $\fc_\reg(F)$ which is relatively compact in $\fc(F)$, there exists $c_B>0$ such that 
$$ \int_{B} |D^\fs(X)|_F^{-\varepsilon'} dX\leq c_B. $$ 
We deduce that for $m\in\BZ$, 
\begin{equation}\label{strengthenlem44}
 \vol(B\cap\{X\in\fc_\reg(F): |D^\fs(X)|_F=q^{-\frac{m}{2}}\})\leq c_B q^{-\frac{\varepsilon'm}{2}}. 
\end{equation}
 
We claim that for any $B$ as above, there exists $c'_B>0$ such that
$$ \int_{B\cap\fc(\varepsilon, T)} (\| T\| - \log |D^\fs(X)|_F)^{d_1} dX \leq c'_B e^{-\frac{\varepsilon'\varepsilon\| T\|}{2}}. $$
This is an analogue of the exercise in \cite[p. 32]{MR1114210} and we include here a proof for completeness. For $X\in\fc(\varepsilon,T)$, we have
$$ \|T\|\leq-\frac{1}{\varepsilon}\log|D^\fs(X)|_F. $$
Therefore, 
$$ \int_{B\cap\fc(\varepsilon, T)} (\| T\| - \log |D^\fs(X)|_F)^{d_1} dX \leq \left(1+\frac{1}{\varepsilon}\right)^{d_1}\int_{B\cap\fc(\varepsilon, T)} (- \log |D^\fs(X)|_F)^{d_1} dX. $$
Since
$$ B\cap\fc(\varepsilon, T)=\coprod_{m\geq\frac{2\varepsilon\|T\|}{\log q}} (B\cap\{X\in\fc_\reg(F): |D^\fs(X)|_F=q^{-\frac{m}{2}}\}), $$
we have
$$ \int_{B\cap\fc(\varepsilon, T)} (- \log |D^\fs(X)|_F)^{d_1} dX=\sum_{m\geq\frac{2\varepsilon\|T\|}{\log q}} \left(\frac{m\log q}{2}\right)^{d_1} \vol(B\cap\{X\in\fc_\reg(F): |D^\fs(X)|_F=q^{-\frac{m}{2}}\}). $$
Applying (\ref{strengthenlem44}), we obtain
$$ \int_{B\cap\fc(\varepsilon, T)} (- \log |D^\fs(X)|_F)^{d_1} dX \leq \sum_{m\geq\frac{2\varepsilon\|T\|}{\log q}} \left(\frac{m\log q}{2}\right)^{d_1} c_B q^{-\frac{\varepsilon'm}{2}}. $$
Now we can confirm our claim by noting the basic fact: for $d>0$ and $a>1$, there exists $c_{d,a}>0$ such that
$$ \sum_{m\geq x} \frac{m^d}{a^m}\leq c_{d,a}a^{-\frac{x}{2}}, \forall x\geq 0. $$

Taking
$$ B=\{X\in\fc_\reg(F): K^T(\eta, X, f, f')\neq 0\}, $$
we see that
$$ \int_{\fc(\varepsilon, T)} |K^T(\eta, X, f, f')| dX \leq c'_\sigma c_3^2 c'_B e^{-\frac{\varepsilon'\varepsilon\| T\|}{2}}. $$
This is half of the lemma. 

It is proved in \cite[p. 32]{MR1114210} that there exist $c_4,d_2>0$ such that for all $x,y\in G(F)$ and $T\in\fa_{M_0,F}\cap\ov{\fa_{P_0}^+}$ with $\|T\|\geq\varepsilon_0$, 
$$ v_M(x, y, T) \leq c_4 (\| T\| + \log \| x\| + \log \| y\|)^{d_2}. $$
By the same argument as before, we obtain
$$ \int_{\fc(\varepsilon, T)} |J^T(\eta, X, f, f')| dX \leq c_5 e^{-\frac{\varepsilon'\varepsilon\| T\|}{2}} $$
for some $c_5>0$. This establishes the other half of the lemma. 
\end{proof}

\begin{lem}\label{lem4.4}
Suppose that $\delta>0$. Then there exist $c, \varepsilon_1, \varepsilon_2>0$ such that
$$ |u_M(x,y,T)-v_M(x,y,T)|\leq ce^{-\varepsilon_1\|T\|} $$
for all $T\in\fa_{M_0,F}\cap\ov{\fa_{P_0}^+}$ with $d(T)\geq\delta\|T\|$, and all $x,y\in\{x\in G(F): \|x\|\leq e^{\varepsilon_2\|T\|}\}$. 
\end{lem}

\begin{proof}
This is Arthur's main geometric lemma \cite[Lemma 4.4]{MR1114210}. 
\end{proof}

\begin{lem}\label{lemofprop45}
Suppose that $\delta>0$. Let $\fc\subseteq\fm\cap\fs$ be an $M$-elliptic Cartan subspace. Then there exist $c, \varepsilon>0$ such that
$$ \int_{\fc_\reg(F)} |K^T(\eta, X, f, f') - J^T(\eta, X, f, f')| dX \leq c e^{-\varepsilon\| T\|} $$
for all $T\in\fa_{M_0,F}\cap\ov{\fa_{P_0}^+}$ with sufficiently large $\|T\|$ and $d(T)\geq\delta\| T\|$. 
\end{lem}

\begin{proof}
Fix $\varepsilon_2>0$ to be the constant given by Lemma \ref{lem4.4}. Let $x,y\in H(F)$ and $X\in\fc_\reg(F)-\fc(\frac{\varepsilon_2}{2},T)$, and assume that
$$ f(\Ad(x^{-1})(X))f'(\Ad(y^{-1})(X))\neq0. $$
Let $\sigma\subseteq\fs(F)$ be a compact subset containing $\Supp(f)\cup\Supp(f')$. From Lemma \ref{lemIII.4}, there exists $c_\sigma>0$ such that
$$ \inf_{\tau_1\in T_\fc(F)} \|\tau_1 x\|, \inf_{\tau_2\in T_\fc(F)} \|\tau_2 y\|\leq c_\sigma\sup\{1,|D^\fs(X)|_F^{-1}\}. $$
Since $X\in\fc_\reg(F)-\fc(\frac{\varepsilon_2}{2},T)$, we have
$$ \sup\{1,|D^\fs(X)|_F^{-1}\}\leq\sup\{1,e^{\frac{\varepsilon_2\|T\|}{2}}\}=e^{\frac{\varepsilon_2\|T\|}{2}}. $$
Then, multiplying $x$ and $y$ by elements in $T_\fc(F)$ if necessary, and taking $\|T\|\geq\frac{2\log c_\sigma}{\varepsilon_2}$, we can assume that
$$ \|x\|,\|y\|\leq e^{\varepsilon_2\|T\|}. $$
It follows from Lemma \ref{lem4.4} that
$$ |u_M(x,y,T)-v_M(x,y,T)|\leq ce^{-\varepsilon_1\|T\|}. $$
By the definitions (see (\ref{firstwoi}) and (\ref{secondwoi})) of $K^T(\eta, X, f, f')$ and $J^T(\eta, X, f, f')$, we obtain that
$$ \int_{\fc_\reg(F)-\fc(\frac{\varepsilon_2}{2},T)} |K^T(\eta, X, f, f')-J^T(\eta, X, f, f')| dX\leq c_1 e^{-\varepsilon_1\| T\|}, $$
where
$$ c_1:=c \int_{\fc_\reg(F)} |D^\fs(X)|_F \int_{(A_M(F)\bs H(F))^2} |f(\Ad(x^{-1})(X)) f'(\Ad(y^{-1})(X))| dx dy dX $$
is finite by Proposition \ref{bdoi1} and Lemma \ref{lem28}. 

One can draw the conclusion by combining this with Lemma \ref{lem43}. 
\end{proof}

Define
$$ J^T(\eta, f, f'):=\sum_{M\in\msl^{G,\omega}(M_0)} |W_0^{M_n}| |W_0^{GL_n}|^{-1} \sum_{\fc\in\mst_\el(\fm\cap\fs)} |W(M_H, \fc)|^{-1} \int_{\fc_\reg(F)} J^T(\eta, X, f, f') dX, $$
where $J^T(\eta, X, f, f')$ is defined by (\ref{secondwoi}). 

\begin{prop}\label{prop45}
Suppose that $\delta>0$. Then there exist $c,\varepsilon>0$ such that
$$ |K^T(\eta, f, f')-J^T(\eta, f, f')| \leq c e^{-\varepsilon\| T\|} $$
for all $T\in\fa_{M_0,F}\cap\ov{\fa_{P_0}^+}$ with sufficiently large $\|T\|$ and $d(T)\geq\delta\| T\|$. 
\end{prop}

\begin{proof}
Apply Lemma \ref{lemofprop45}. 
\end{proof}

It is proved in \cite[(6.5) in p. 46]{MR1114210} that $v_M(x,y,T)$ is an exponential polynomial in $T\in R_0\cap(\ov{\fa_{P_0}^+}/\fa_G)$. 
Denote by $\tilde{v}_M(x, y)$ the constant term of $v_M(x, y, T)$ as in \cite[(6.6) in p. 46]{MR1114210}. Then for $f,f'\in\CC_c^\infty(\fs(F))$ and $X\in(\fm\cap\fs_\rs)(F)_\el$, $J^T(\eta, X, f, f')$ is also an exponential polynomial in $T\in R_0\cap(\ov{\fa_{P_0}^+}/\fa_G)$ whose constant term is given by
\begin{equation}\label{cstsecondwoi}
 \tilde{J}_M(\eta, X, f, f'):=|D^\fs(X)|_F \int_{(A_M(F)\bs H(F))^2} f(\Ad(x^{-1})(X)) f'(\Ad(y^{-1})(X)) \eta(\Nrd(x^{-1}y)) \tilde{v}_M(x, y) dx dy. 
\end{equation}
Thus for $f,f'\in\CC_c^\infty(\fs(F))$, $J^T(\eta, f, f')$ is still an exponential polynomial in $T\in R_0\cap(\ov{\fa_{P_0}^+}/\fa_G)$ whose constant term is given by
\begin{equation}\label{cstsecondgeom}
 \tilde{J}(\eta, f, f'):=\sum_{M\in\msl^{G,\omega}(M_0)} |W_0^{M_n}| |W_0^{GL_n}|^{-1} \sum_{\fc\in\mst_\el(\fm\cap\fs)} |W(M_H, \fc)|^{-1} \int_{\fc_\reg(F)} \tilde{J}_M(\eta, X, f, f') dX. 
\end{equation}

\begin{coro}\label{corV.4}
For $f,f'\in\CC_c^\infty(\fs(F))$, we have the equality
$$ \tilde{J}(\eta, f, \hat{f'})=\tilde{J}(\eta, \hat{f}, f'). $$
\end{coro}

\begin{proof}
By the Plancherel formula, for $x\in H(F)$, we have
$$ k(x, f, \hat{f'})=k(x, \hat{f}, f'). $$
Then for all $T\in R_0\cap(\ov{\fa_{P_0}^+}/\fa_G)$, 
$$ K^T(\eta, f, \hat{f'})=K^T(\eta, \hat{f}, f'). $$
Finally, apply Proposition \ref{prop45} to conclude. 
\end{proof}

\begin{lem}\label{lemV.5}
For all $Q\in\msf^{G,\omega}(M_0)$, there exists a constant $c'_Q$ such that for all $f,f'\in\CC_c^\infty(\fs(F))$, we have the equality
$$ \tilde{J}(\eta, f, f')=\sum_{Q\in\msf^{G,\omega}(M_0)} |W_0^{M_{Q_n}}| |W_0^{GL_n}|^{-1} (-1)^{\dim(A_Q/A_G)}J^{M_Q}(\eta,f_{\ov{Q}}^\eta,{f'}_Q^\eta)c'_Q, $$
where $J^{M_Q}(\eta,f_{\ov{Q}}^\eta,{f'}_Q^\eta)$ is defined by (\ref{geomnoninvltf1}). 
\end{lem}

\begin{proof}
Suppose that $M\in\msl^{G,\omega}(M_0)$. It is shown in \cite[p. 92]{MR1114210} that 
$$ \tilde{v}_M(x,y)=(-1)^{\dim(A_M/A_G)}\sum_{Q\in\msf^G(M)} v_M^Q(x,y) c'_Q, $$
where $c'_Q$ is a constant for each $Q\in\msf^{G,\omega}(M_0)$. 

Now substitude this in the definition (\ref{cstsecondwoi}) of $\tilde{J}_M(\eta, X, f, f')$. Note that
$$ v_M^Q(m_1 n_1 k_1, m_2 n_2 k_2)=v_M^{M_Q}(m_1, m_2) $$
for $m_1\in M_{\ov{Q}_H}(F), n_1\in N_{\ov{Q}_H}(F), m_2\in M_{Q_H}(F), n_2\in N_{Q_H}(F), k_1,k_2\in K_H$. By the same argument as the proof of Proposition \ref{propwoi1}.4), one shows that
$$ |D^\fs(X)|_F \int_{(A_M(F)\bs H(F))^2} f(\Ad(x^{-1})(X)) f'(\Ad(y^{-1})(X)) \eta(\Nrd(x^{-1}y)) v_M^Q(x, y) dx dy=J_M^{M_Q}(\eta,X,f_{\ov{Q}}^\eta,{f'}_Q^\eta), $$
where $J_M^{M_Q}(\eta,X,f_{\ov{Q}}^\eta,{f'}_Q^\eta)$ is defined by (\ref{defV.1.1}). Therefore, we have
$$ \tilde{J}_M(\eta, X, f, f')=(-1)^{\dim(A_M/A_G)} \sum_{Q\in\msf^G(M)} J_M^{M_Q}(\eta,X,f_{\ov{Q}}^\eta,{f'}_Q^\eta) c'_Q. $$
Then the lemma follows from the definition (\ref{cstsecondgeom}) of $\tilde{J}(\eta, f, f')$. 
\end{proof}

\begin{proof}[Proof of Theorem \ref{noninvltf1}]
Using Lemma \ref{lemV.5} and Corollary \ref{corV.4}, we can prove the theorem by induction on the dimension of $G$. 
\end{proof}

\subsection{The case of $(G',H')$}

Suppose that $M'\in\msl^{H'}(M'_0)$. For $x,y\in H'(F)$, we define an $(H',M')$-family $(v_{P'}(x,y))_{P'\in\msp^{H'}(M')}$ by
$$ v_{P'}(\lambda,x,y):=e^{-\lambda(H_{P'}(y)-H_{\ov{P'}}(x))}, \forall \lambda\in i\fa_{M'}^\ast, P'\in\msp^{H'}(M'), $$
where $\ov{P'}\in\msp^{H'}(M')$ is the parabolic subgroup opposite to $P'$. For $f, f'\in\CC_c^\infty(\fs'(F))$ and $Y\in(\wt{\fm'}\cap\fs'_\rs)(F)_\el$, we define
\begin{equation}\label{defV.1.2}
 J_{M'}^{H'}(Y, f, f'):=|D^{\fs'}(Y)|_F \int_{(A_{M'}(F)\bs H'(F))^2} f(\Ad(x^{-1})(Y)) f'(\Ad(y^{-1})(Y)) v_{M'}(x, y) dx dy. 
\end{equation}

\begin{prop}\label{propV.1.2}
Suppose that $M'\in\msl^{H'}(M'_0)$ and that $f,f'\in\CC_c^\infty(\fs'(F))$. 

1) The integral (\ref{defV.1.2}) is absolutely convergent. 

2) The function $J_{M'}^{H'}(\cdot, f, f')$ is locally constant on $(\wt{\fm'}\cap\fs'_\rs)(F)_\el$. 

3) If $\fc'\subseteq\wt{\fm'}\cap\fs'$ is an $M'$-elliptic Cartan subspace, the restriction of $J_{M'}^{H'}(\cdot, f, f')$ to $\fc'_\reg(F)$ vanishes outside a compact subset of $\fc'(F)$. 

4) If $w\in \Norm_{H'(F)}(M'_0)$, $x\in M'(F)$ and $k, k'\in K_{H'}$, we have the equality
$$ J_{\Ad(w)(M')}^{H'}(\Ad(wx)(Y), \Ad(k)(f), \Ad(k')(f'))=J_{M'}^{H'}(Y, f, f') $$
for all $Y\in(\wt{\fm'}\cap\fs'_\rs)(F)_\el$. 

5) There exist constants $c>0$ and $N\in\BN$ such that for all $Y\in(\wt{\fm'}\cap\fs'_\rs)(F)_\el$, we have
$$ |J_{M'}^{H'}(Y, f, f')| \leq c \sup\{1, -\log |D^{\fs'}(Y)|_F\}^N. $$

6) For all $Y\in(\wt{\fm'}\cap\fs'_\rs)(F)_\el$, we have
$$ J_{M'}^{H'}(Y, f, f')=\sum_{L'_1,L'_2\in\msl^{H'}(M')} d_{M'}^{H'}(L'_1,L'_2) J_{M'}^{L'_1}(Y,f_{\ov{Q'_1}}) J_{M'}^{L'_2}(Y,f'_{Q'_2}), $$
where $(Q'_1,Q'_2):=s(L'_1,L'_2)$ (see Section \ref{gmmaps}). 
\end{prop}

\begin{proof}
It is almost the same as the proof of Proposition \ref{propV.1.1}, except that one needs to use Proposition \ref{propwoi2}.4) and Corollary \ref{corIII.6.2} to show 6) and 5). 
\end{proof}

For $f,f'\in\CC_c^\infty(\fs'(F))$, we define
\begin{equation}\label{geomnoninvltf2}
\begin{split}
 J^{H'}(f, f'):=&\sum_{M'\in\msl^{H'}(M'_0)} |W_0^{H'}||W_0^{M'}|^{-1} (-1)^{\dim(A_{M'}/A_{H'})} \sum_{\fc'\in\mst_\el(\wt{\fm'}\cap\fs')} |W(M', \fc')|^{-1} \int_{\fc'_\reg(F)} \\
&J_{M'}^{H'}(Y, f, f') dY. 
\end{split}
\end{equation}
This expression is absolutely convergent by Proposition \ref{propV.1.2}.5) and Corollary \ref{cor20.2}. One may extend in the obvious way the definitions (\ref{defV.1.2}) and (\ref{geomnoninvltf2}) to the symmetric pair $(\wt{L'}, L', \Ad(\alpha))$, where $L'\in\msl^{H'}(M'_0)$. 

\begin{remark}\label{rmkV.1.2}
We have the equality
$$ J^{H'}(f, f')=J^{H'}(f', f). $$
It results from the fact that for all $M'\in\msl^{H'}(M'_0)$ and all $x,y\in H'(F)$, we have $v_{M'}(x,y)=v_{M'}(y,x)$. 
\end{remark}

\begin{thm}[Noninvariant trace formula]\label{noninvltf2}
For all $f, f'\in\CC_c^\infty(\fs'(F))$, we have the equality
$$ J^{H'}(f, \hat{f'})=J^{H'}(\hat{f}, f'). $$
\end{thm}

\begin{proof}
We may simply copy the proof of Theorem \ref{noninvltf1} here with obvious modifications. Especially, one needs to use Proposition \ref{bdoi2} to show analogues of Lemmas \ref{lem43} and \ref{lemofprop45} for the case of $(G',H')$. 
\end{proof}


\section{\textbf{Howe's finiteness for weighted orbital integrals}}\label{sechowfin}

\subsection{The case of $(G,H)$}

For an open compact subgroup $r$ of $\fs(F)$, denote by $\CC_c^\infty(\fs(F)/r)$ the subspace of $\CC_c^\infty(\fs(F))$ consisting of the functions invariant by translation of $r$. Let ${i_r^\fs}^*: \CC_c^\infty(\fs(F))^*\ra\CC_c^\infty(\fs(F)/r)^*$ be the dual map of the natural injection $\CC_c^\infty(\fs(F)/r)\hookrightarrow\CC_c^\infty(\fs(F))$. 

For any set $\sigma$, denote by $\BC[\sigma]$ the $\BC$-linear space of maps from $\sigma$ to $\BC$ with finite support. For $M\in\msl^{G,\omega}(M_0)$ and $\sigma\subseteq (\fm\cap\fs_\rs)(F)$, we define the linear map
\[\begin{split}
\delta_M^G: \BC[\sigma] &\ra \CC_c^\infty(\fs(F))^*, \\
(z_X)_{X\in\sigma} &\mapsto \sum_{X\in\sigma} z_X J_M^G(\eta, X, \cdot), \\
\end{split}\]
where $z_X\in\BC$ is the coordinate at $X\in\sigma$. 

For $L\in\msl^{G,\omega}(M_0)$ and $r_L\subseteq(\fl\cap\fs)(F)$ an open compact subgroup, we define ${i_{r_L}^{\fl\cap\fs}}^*$ and $\delta_M^L$ in the same way by using $(L, L_H, \Ad(\epsilon))$ in place of $(G, H, \Ad(\epsilon))$. 

\begin{prop}[Howe's finiteness]\label{howe1}
Let $r$ be an open compact subgroup of $\fs(F)$, $M\in\msl^{G,\omega}(M_0)$ and $\sigma\subseteq (\fm\cap\fs_\rs)(F)$. Suppose that there exists a compact subset $\sigma_0\subseteq(\fm\cap\fs)(F)$ such that $\sigma\subseteq \Ad(M_H(F))(\sigma_0)$. Then the image of the linear map
$$ {i_r^\fs}^*\circ\delta_M^G: \BC[\sigma]\ra\CC_c^\infty(\fs(F)/r)^* $$
is of finite dimension. 
\end{prop}

\begin{remark}
For $M=G$, Proposition \ref{howe1} is essentially included in a more general result \cite[Theorem 6.1]{MR1375304} in the context of $\theta$-groups (in the sense of \cite[p. 467]{zbMATH03577478}). 
\end{remark}

The rest of this section is devoted to the proof of Proposition \ref{howe1}. We shall follow the main steps in \cite[\S IV.2-6]{MR1344131}. Since the function $\kappa(\cdot)J_M^G(\eta,\cdot,f)$ is constant on $\Ad(M_H(F))(X)$ for $f\in\CC_c^\infty(\fs(F))$ and $X\in(\fm\cap\fs_\rs)(F)$, and $\mst(\fm\cap\fs)$ is a finite set, by Lemma \ref{lem28}, we may and shall suppose that $\sigma$ is relatively compact in $(\fm\cap\fs)(F)$. The proposition will be proved by induction on the dimension of $G$. 

Recall that we have chosen the standard maximal compact subgroup $K=GL_{2n}(\CO_D)$ of $G(F)=GL_{2n}(D)$. Let $k=\fg(\CO_F):=\fg\fl_{2n}(\CO_D)$, which is an $\CO_F$-lattice in $\fg(F)=\fg\fl_{2n}(D)$ and is stable under the adjoint action of $K$. Since $H\in\msl^G(M_0)$, we have set $K_H=K\cap H(F)=GL_n(\CO_D)\times GL_n(\CO_D)$. Let $\fh(\CO_F):=k\cap\fh(F)$ and $\fs(\CO_F):=k\cap\fs(F)$. Then we see that $k=\fh(\CO_F)\oplus\fs(\CO_F)$ and that $\fs(\CO_F)$ is stable under the adjoint action of $K_H$. For all $P\in\msf^G(M_0)$, we fix $a_P\in A_P(F)$ such that $|\alpha(a_P)|_F<1, \forall \alpha\in\Delta_P^G$. 

Recall that we denote by $\CN^\fs$ the set of nilpotent elements in $\fs(F)$ and fix a uniformiser $\varpi$ of $\CO_F$. Let $X\in\CN^\fs\cap(k-\varpi k)$. By the Jacobson-Morozov theorem for symmetric spaces (Lemma \ref{jm}), there exists a group homomorphism $\varphi: SL_2(F)\ra G(F)$ such that
$$ d\varphi\mat(0,1,0,0)=X \text{ and } a^\varphi:=\varphi\mat(\varpi,,,\varpi^{-1})\in H(F). $$
We define the parabolic subgroup $P^X$ of $G$ as in \cite[\S IV.3]{MR1344131}. More concretely, set
$$ \fg[i]:=\{Y\in\fg: \Ad(a^\varphi)(Y)=\varpi^i Y, \forall i\in\BZ\} $$
and
$$ \fp^X:=\bigoplus_{i\geq0} \fg[i]; $$
then let
\begin{equation}\label{defP^X}
 P^X:=\{x\in G: \Ad(x)(\fp^X)=\fp^X\}. 
\end{equation}
Note that $P^X$ is independent of the choice of $\varphi$ by \cite[Proposition 5.7.1]{MR794307}. Since $a^\varphi$ commutes with $\epsilon$, one has $\epsilon\in P^X$. By Lemma \ref{intpar1}, $P^X\cap H$ is a parabolic subgroup of $H$. Then there exists an element $x\in K_H$ such that $P':=\Ad(x)(P^X)\in\msf^G(M_0)$. We shall fix such an $x$. Let $a^X:=\Ad(x^{-1})(a_{P'})\in H(F)$. Note that $a^X$ depends on the choice of $x$, but this is unimportant. It is proved in \cite[(3) in \S IV.3]{MR1344131} that 
$$ \Ad(a^X)(X)\in\varpi k. $$
But $\Ad(a^X)(X)\in\fs(F)$. Thus we have
\begin{equation}\label{equIV.3(3)1}
 \Ad(a^X)(X)\in\varpi\fs(\CO_F). 
\end{equation}

\begin{lem}\label{lemIV.3}
There exists an integer $h\in\BN$ such that for all $Y\in\varpi^h \fs(\CO_F)$, all integers $l\geq h$ and all $Z\in\Ad(a^X)^{-1}(\varpi^l \fs(\CO_F))$, there exists $\gamma\in K_H$ with $\eta(\Nrd(\gamma))=1$ such that
$$ \Ad(\gamma)(X+Y+Z)\in X+Y+\varpi^l \fs(\CO_F). $$
\end{lem}

\begin{proof}
It is analogous to the proof of \cite[Lemme IV.3]{MR1344131}. The point is that we require $\gamma\in H(F)$ and $\eta(\Nrd(\gamma))=1$ here. 

Set
$$ \fn^X:=\bigoplus_{i\geq1} \fg[i]. $$
Since $a^\varphi$ commutes with $\epsilon$, one has $\fg[i]=(\fg[i]\cap\fh)\oplus(\fg[i]\cap\fs)$. 

It is proved in \cite[(4) in \S IV.3]{MR1344131} that 
$$ \Ad(a^X)^{-1}(k)\subseteq k+\Ad(a^X)^{-1}(k)\cap\fn^X.  $$
By the intersection of $\fs(F)$ with both sides in the last expression, we have
\begin{equation}\label{lemIV.3(4)}
 \Ad(a^X)^{-1}(\fs(\CO_F))\subseteq \fs(\CO_F)+\Ad(a^X)^{-1}(\fs(\CO_F))\cap\fn^X. 
\end{equation}

It is shown in \cite[(2) in \S IV.3]{MR1344131} and the proof of \cite[Lemme IV.3]{MR1344131} that there exists $c\in\BN$ such that 
$$ \fn^X\cap\Ad(a^X)^{-1}(k)\subseteq \ad(X)(\varpi^{-c}k).  $$
By the intersection of $\fs(F)$ with both sides in the last expression, we have 
\begin{equation}\label{lemIV.3(2)}
 \fn^X\cap\Ad(a^X)^{-1}(\fs(\CO_F))\subseteq \ad(X)(\varpi^{-c}\fh(\CO_F)). 
\end{equation}
Fix such an integer $c$. 
There exists $c'\in\BN$ such that for all $l\geq c'$, we have 
\begin{itemize}
	\item $\varpi^l k\subseteq V_\fg$; 
	\item $\exp(\varpi^l k)\subseteq \{x\in K: \eta(\Nrd(x))=1\}$; 
	\item $\Ad(\exp(Y'))(Y)\equiv Y+\ad(Y')(Y) \mod \varpi^{l+c}k, \forall Y\in k, Y'\in\varpi^l k$. 
\end{itemize}
Fix such an integer $c'$. 
There exists $c''\in\BN$ such that $\Ad(a^X)^{-1}(\varpi^{c''}k)\subseteq k$. Fix such an integer $c''$. 

Let $h\in\BN$ be such that $h\geq\max\{c+c', c+c''\}$. Let $Y,l$ and $Z$ be as in the statement. Thanks to (\ref{lemIV.3(4)}), we can write $Z=Z_1+Z_2$, where $Z_1\in\varpi^l\fs(\CO_F)$ and $Z_2\in\Ad(a^X)^{-1}(\varpi^l\fs(\CO_F))\cap\fn^X$. Because of (\ref{lemIV.3(2)}), we can choose $Z'\in\varpi^{l-c}\fh(\CO_F)$ such that $Z_2=\ad(X)(Z')$. Since $Z'\in\varpi^{c'}\fh(\CO_F)$ from the hypothesis $l\geq h\geq c+c'$, we can define $\gamma:=\exp(Z')$. Then we have $\gamma\in K_H$ and $\eta(\Nrd(\gamma))=1$. 
It is shown in the proof of \cite[Lemme IV.3]{MR1344131} that $\gamma$ verifies 
$$ \Ad(\gamma)(X+Y+Z)\in X+Y+\varpi^l k. $$
But both of $\Ad(\gamma)(X+Y+Z)$ and $X+Y$ belong to $\fs(F)$. Hence, we deduce that 
$$ \Ad(\gamma)(X+Y+Z)\in X+Y+\varpi^l \fs(\CO_F). $$
\end{proof}

For $X\in\CN^\fs\cap(k-\varpi k)$, we shall fix an integer $h^X$ such that

(1) $h^X$ verifies the condition of Lemma \ref{lemIV.3}; 

(2) $h^X\geq 1$; 

(3) $\Ad(a^X)(\varpi^{h^X}\fs(\CO_F))\subseteq\varpi \fs(\CO_F)$. 

Denote by $\CN^G$ the set of nilpotent elements in $\fg(F)$. Let $\BP\fg(F)$ (resp. $\BP\fs(F)$) be the projective space associated to $\fg(F)$ (resp. $\fs(F)$) and $\pi: \fg(F)-\{0\}\ra \BP\fg(F)$ the natural projection. Since $\pi(\CN^G-\{0\})$ is compact and $\CN^\fs-\{0\}$ is a closed subset of $\CN^G-\{0\}$, we know that $\pi(\CN^\fs-\{0\})$ is compact. One also sees that
$$ \pi(\CN^\fs\cap(k-\varpi k))=\pi(\CN^\fs-\{0\}). $$
We can and shall choose a finite set $\CN_0\subseteq\CN^\fs\cap(k-\varpi k)$ such that
$$ \bigcup_{X\in\CN_0} \pi(X+\varpi^{h^X}\fs(\CO_F)) $$
is an open neighbourhood of $\pi(\CN^\fs-\{0\})\subseteq\BP\fs(F)$. 

\begin{lem}\label{lemIV.4}
There exists an integer $c\in\BN$ such that for all $d\in\delta_M^G(\BC[\sigma])$ and all $f\in\CC_c^\infty(\fs(F))$ satisfying $d(f)\neq 0$, we have
$$ \Supp(f)\cap[\varpi^{-c}\fs(\CO_F)\cup\bigcup_{X\in\CN_0} F^\times(X+\varpi^{h^X}\fs(\CO_F))]\neq\emptyset. $$
\end{lem}

\begin{proof}
Recall that $\sigma$ is assumed to be relatively compact in $(\fm\cap\fs)(F)$. Fix an open compact neighbourhood $\sigma'$ of $\sigma$ in $\fs(F)$. Fix $P_0\in\msp^H(M_0)$ and set
$$ A_{P_0}^+:=\{a\in A_{M_0}(F): |\alpha(a)|_F\geq 1, \forall \alpha\in\Delta_{P_0}^H\}. $$
Similarly, for all $B\in\msp^G(M_0)$, set
$$ A_B^+:=\{a\in A_{M_0}(F): |\alpha(a)|_F\geq 1, \forall \alpha\in\Delta_B^G\}. $$
We see from the argument of \cite[Lemma 4.11]{zbMATH07499568} that
$$ A_{P_0}^+\subseteq \bigcup_{\{B\in\msp^G(M_0): P_0\subseteq B\}} A_B^+. $$

By the Cartan decomposition, there exists a compact subset $\Gamma\subseteq H(F)$ such that $H(F)=K_H A_{P_0}^+\Gamma$. 
Fix such a $\Gamma$. Then
$$ H(F)\subseteq \bigcup_{\{B\in\msp^G(M_0): P_0\subseteq B\}} K_H A_B^+ \Gamma. $$
Fix $c'\in\BN$ such that $\Ad(\Gamma)(\sigma')\subseteq\varpi^{-c'}\fs(\CO_F)$. Since $\Ad(A_B^+)(\fs(\CO_F))\subseteq \fs(\CO_F)+(\fn_B\cap\fs)(F)$ and $(\fn_B\cap\fs)(F)\subseteq\CN^\fs$, we obtain
$$ \Ad(H(F))(\sigma')\subseteq\varpi^{-c'}\fs(\CO_F)+\CN^\fs. $$
Choose an integer $c$ such that
$$ c\geq c'+\sup_{X\in\CN_0} h^X. $$

Let $d$ and $f$ be as in the statement. It is evident that
$$ \Supp(f)\cap\Ad(H(F))(\sigma')\neq\emptyset. $$
Suppose that $X_1\in\varpi^{-c'}\fs(\CO_F), X_2\in\CN^\fs$ satisfy $X_1+X_2\in\Supp(f)$. If $X_2\in\varpi^{-c}\fs(\CO_F)$, we are done. If not, let $X\in\CN_0$ be such that $\pi(X_2)\in\pi(X+\varpi^{h^X}\fs(\CO_F))$. Then there exists $\lambda\in F^\times$ such that $X_2\in\lambda(X+\varpi^{h^X}\fs(\CO_F))$. Since $X_2\notin\varpi^{-c}\fs(\CO_F)$, we have $v_F(\lambda)<-c$ and thus $v_F(\lambda)+h^X<-c'$. Therefore, $X_1\in\lambda\varpi^{h^X}\fs(\CO_F)$ and then $X_1+X_2\in\lambda(X+\varpi^{h^X}\fs(\CO_F))\subseteq F^\times(X+\varpi^{h^X}\fs(\CO_F))$. This is exactly what we want to prove. 
\end{proof}

We shall fix an integer $c$ verifying the condition in Lemma \ref{lemIV.4}. We shall also fix another integer $h$ such that
\begin{equation}\label{condonh}
 h\geq\sup_{X\in\CN_0} h^X. 
\end{equation}
Denote
$$ C:=\{f\in\CC_c^\infty(\fs(F)/\varpi^h\fs(\CO_F)): \Supp(f)\subseteq\varpi^{-c}\fs(\CO_F)\}, $$
which is a $\BC$-linear space of finite dimension. For all $L\in\msl^G(M), L\neq G$, let $r(L):=\varpi^h(\fl\cap\fs)(\CO_F)$. 

\begin{lem}\label{lemIV.5}
Let $\fz=(z_X)_{X\in\sigma}\in\BC[\sigma]$. Suppose that

1) ${i_{r(L)}^{\fl\cap\fs}}^*\circ\delta_M^L(\fz)=0$ for all $L\in\msl^G(M), L\neq G$;

2) $\delta_M^G(\fz)(C)=0$. 

Then ${i_{\varpi^h\fs(\CO_F)}^\fs}^*\circ\delta_M^G(\fz)=0$. 
\end{lem}

\begin{proof}
Write $d:=\delta_M^G(\fz)$. It suffices to prove by induction on the integer $e\geq c$ the assertion 

$(A)_e$: for all $f\in\CC_c^\infty(\fs(F)/\varpi^h\fs(\CO_F))$ with $\Supp(f)\subseteq\varpi^{-e}\fs(\CO_F)$, we have $d(f)=0$. 

If $e=c$, this is true by the hypothesis 2). Fix an $e>c$ and suppose that $(A)_{e-1}$ is true. For any open compact subset $s\subseteq\fg(F)$, denote by $1_s$ its characteristic function. It suffices to prove that for all $Y\in\varpi^{-e}\fs(\CO_F)-\varpi^{-e+1}\fs(\CO_F)$, we have $d(1_{Y+\varpi^h\fs(\CO_F)})=0$. 

Suppose that $Y\in\varpi^{-e}\fs(\CO_F)-\varpi^{-e+1}\fs(\CO_F)$. The hypothesis $e>c$ implies $\varpi^{-c}\fs(\CO_F)\subseteq\varpi^{-e+1}\fs(\CO_F)$ and thus $Y\notin\varpi^{-c}\fs(\CO_F)$. Suppose on the contary that $d(1_{Y+\varpi^h\fs(\CO_F)})\neq0$. By Lemma \ref{lemIV.4}, there exists $X\in\CN_0$ and $\lambda\in F^\times$ such that $(Y+\varpi^h\fs(\CO_F))\cap\lambda(X+\varpi^{h^X}\fs(\CO_F))\neq\emptyset$. Fix such $X$ and $\lambda$. Since $v_F(Y)=-e$ and $v_F(X)=0$, we have $v_F(\lambda)=-e$. As $h\geq h^X\geq h^X-e$ (see (\ref{condonh})), we obtain $\varpi^h\fs(\CO_F)\subseteq\lambda\varpi^{h^X}\fs(\CO_F)$ and then $Y\in\lambda(X+\varpi^{h^X}\fs(\CO_F))$. Let $Y'\in\varpi^{h^X}\fs(\CO_F)$ such that $Y=\lambda(X+Y')$. 

Let $Z\in\Ad(a^X)^{-1}(\varpi^h\fs(\CO_F))$. Since $h+e\geq h^X$, we can apply Lemma \ref{lemIV.3} to $X, Y',\lambda^{-1}Z$ and $l:=h+e$. Then there exists $\gamma\in K_H$ with $\eta(\Nrd(\gamma))=1$ such that
$$ \Ad(\gamma)(X+Y'+\lambda^{-1}Z)\in X+Y'+\varpi^{h+e}\fs(\CO_F). $$
From $v_F(\lambda)=-e$, we deduce that
$$ \Ad(\gamma)(\lambda(X+Y')+Z)\in\lambda(X+Y')+\varpi^h\fs(\CO_F), $$
i.e., 
$$ \Ad(\gamma)(Y+Z)\in Y+\varpi^h\fs(\CO_F). $$
Since $\gamma\in K_H$, this is equivalent to
$$ \Ad(\gamma)(Y+Z+\varpi^h\fs(\CO_F))=Y+\varpi^h\fs(\CO_F) $$
or
$$ \Ad(\gamma)(1_{Y+Z+\varpi^h\fs(\CO_F)})=1_{Y+\varpi^h\fs(\CO_F)}. $$
By Proposition \ref{propwoi1}.3), we obtain
$$ d(1_{Y+Z+\varpi^h\fs(\CO_F)})=\eta(\Nrd(\gamma))d(1_{Y+\varpi^h\fs(\CO_F)}). $$
Because $\eta(\Nrd(\gamma))=1$, we have
$$ d(1_{Y+Z+\varpi^h\fs(\CO_F)})=d(1_{Y+\varpi^h\fs(\CO_F)}). $$
Now, by the sum over $Z\in\left(\Ad(a^X)^{-1}(\varpi^h\fs(\CO_F))+\varpi^h\fs(\CO_F)\right)/\varpi^h\fs(\CO_F)$ (a finite set), we get
\begin{equation}\label{equIV.5(3)}
 d(1_{Y+\varpi^h\fs(\CO_F)})=[k':\fs(\CO_F)]^{-1}d(1_{Y+\varpi^h k'}), 
\end{equation}
where
\begin{equation}\label{equk'}
 k':=\Ad(a^X)^{-1}(\fs(\CO_F))+\fs(\CO_F). 
\end{equation}

By Proposition \ref{propwoi1}.6), we have the equality
$$ d\left(\Ad(a^X)(1_{Y+\varpi^h k'})-\eta(\Nrd(a^X))1_{Y+\varpi^h k'}\right)=\eta(\Nrd(a^X)) \sum_{Q\in\msf^G(M),Q\neq G} \delta_M^{M_Q}(\fz) \left((1_{Y+\varpi^h k'})_{Q,(a^X)^{-1}}^\eta\right). $$
From (\ref{nonequcst1}), it is clear that $(1_{Y+\varpi^h k'})_{Q,(a^X)^{-1}}^\eta$ is invariant by translation of $r(M_Q)$. By the hypothesis 1), we have $\delta_M^{M_Q}(\fz) \left((1_{Y+\varpi^h k'})_{Q,(a^X)^{-1}}^\eta\right)=0$ for all $Q\in\msf^G(M),Q\neq G$, so
\begin{equation}\label{equIV.5(4)}
 \eta(\Nrd(a^X))d(1_{Y+\varpi^h k'})=d\left(\Ad(a^X)(1_{Y+\varpi^h k'})\right). 
\end{equation}
We see easily that
$$ \Ad(a^X)(1_{Y+\varpi^h k'})=1_s, $$
where
\begin{equation}\label{equs}
 s:=\Ad(a^X)(Y)+\varpi^h\Ad(a^X)(k'). 
\end{equation}

Recall $Y=\lambda(X+Y')$ above. As $X\in\CN_0$, by (\ref{equIV.3(3)1}), we have $\Ad(a^X)(X)\in\varpi\fs(\CO_F)$. Since $Y'\in\varpi^{h^X}\fs(\CO_F)$, by the hypothesis (3) on $h^X$, we have $\Ad(a^X)(Y')\in\varpi\fs(\CO_F)$. For $v_F(\lambda)=-e$, we obtain
\begin{equation}\label{equs1}
 \Ad(a^X)(Y)\in\varpi^{-e+1}\fs(\CO_F). 
\end{equation}
We see from (\ref{equk'}) that
$$ \Ad(a^X)(k')=\fs(\CO_F)+\Ad(a^X)(\fs(\CO_F)). $$
By (\ref{condonh}) and the hypothesis (3) on $h^X$, we have $\varpi^h\Ad(a^X)(\fs(\CO_F))\subseteq\varpi^{h^X}\Ad(a^X)(\fs(\CO_F))\subseteq\varpi\fs(\CO_F)$. Then by (\ref{condonh}) and the hypothesis (2) on $h^X$, we have
\begin{equation}\label{equs2}
 \varpi^h\Ad(a^X)(k')=\varpi^h\fs(\CO_F)+\varpi^h\Ad(a^X)(\fs(\CO_F))\subseteq\varpi^{h^X}\fs(\CO_F)+\varpi\fs(\CO_F)=\varpi\fs(\CO_F). 
\end{equation}
From (\ref{equs}), (\ref{equs1}) and (\ref{equs2}), we see that
$$ \Supp(1_s)\subseteq\varpi^{-e+1}\fs(\CO_F). $$
Since $\fs(\CO_F)\subseteq\Ad(a^X)(k')$, we know that $1_s$ is invariant by translation of $\varpi^h\fs(\CO_F)$. Using the induction hypothesis $(A)_{e-1}$, we have
$$ d(1_s)=0. $$
Thanks to (\ref{equIV.5(3)}) and (\ref{equIV.5(4)}), we obtain
$$ d(1_{Y+\varpi^h\fs(\CO_F)})=0. $$
This proves $(A)_e$ and thus the lemma. 
\end{proof}

\begin{proof}[Proof of Proposition \ref{howe1}]
We use induction on the dimension of $G$. Suppose that for all $L\in\msl^G(M), L\neq G$, and all open compact subgroups $r_L$ of $(\fl\cap\fs)(F)$, the image of the linear map 
$$ {i_{r_L}^{\fl\cap\fs}}^\ast\circ\delta_M^L:\BC[\sigma]\ra\CC_c^\infty((\fl\cap\fs)(F)/r_L)^\ast $$ 
is of finite dimension. This is actually a product form of the proposition in lower dimensions. Now we would like to prove the proposition. The argument below is also valid for the case $G=M$. 

Enlarge $h$ in (\ref{condonh}) if necessary such that $r\supseteq\varpi^h\fs(\CO_F)$. We shall prove that the image of ${i_{\varpi^h\fs(\CO_F)}^\fs}^*\circ\delta_M^G$ is of finite dimension. Admit this for the moment. Since ${i_r^\fs}^*$ factorises by ${i_{\varpi^h\fs(\CO_F)}^\fs}^*$, the image of ${i_r^\fs}^*\circ\delta_M^G$ is also of finite dimension. Then we finish the proof. 

Let $\CK_1$ be the kernel of the linear map
$$ \bigoplus_{L\in\msl^G(M),L\neq G} {i_{r(L)}^{\fl\cap\fs}}^\ast\circ\delta_M^L: \BC[\sigma]\ra\bigoplus_{L\in\msl^G(M),L\neq G}\CC_c^\infty((\fl\cap\fs)(F)/r(L))^\ast, $$
whose image is of finite dimension by our induction hypothesis applied to $r_L:=r(L)$ for all $L\in\msl^G(M),L\neq G$. Hence, to prove that ${i_{\varpi^h\fs(\CO_F)}^\fs}^*\circ\delta_M^G(\BC[\sigma])$ is of finite dimension, it suffices to prove that ${i_{\varpi^h\fs(\CO_F)}^\fs}^*\circ\delta_M^G(\CK_1)$ is of finite dimension. 

Consider the composition of the linear maps
$$ d_1:={i_{\varpi^h\fs(\CO_F)}^\fs}^*\circ\delta_M^G\big|_{\CK_1}: \CK_1\ra\CC_c^\infty(\fs(F)/\varpi^h\fs(\CO_F))^\ast $$
and
$$ \Res: \CC_c^\infty(\fs(F)/\varpi^h\fs(\CO_F))^\ast\ra C^\ast. $$
The latter map is the natural restriction. Lemma \ref{lemIV.5} says that
$$ \ker(\Res\circ d_1)=\ker(d_1), $$
which is denoted by $\CK_2$. Then
$$ d_1(\CK_1)\simeq\CK_1/\CK_2\simeq\Res\circ d_1(\CK_1)\subseteq C^\ast. $$
Since $C$ is of finite dimension, we see that $d_1(\CK_1)$ is of finite dimension. 
\end{proof}

\begin{coro}\label{corhowe1}
Let $r$ be an open compact subgroup of $\fs(F)$, $M\in\msl^{G,\omega}(M_0)$ and $\sigma\subseteq (\fm\cap\fs_\rs)(F)$. Suppose that there exists a compact subset $\sigma_0\subseteq(\fm\cap\fs)(F)$ such that $\sigma\subseteq \Ad(M_H(F))(\sigma_0)$. Then there exists a finite subset $\{X_i: i\in I\}\subseteq\sigma$ and a finite subset $\{f_i: i\in I\}\subseteq\CC_c^\infty(\fs(F)/r)$ such that for all $X\in\sigma$ and all $f\in\CC_c^\infty(\fs(F)/r)$, we have the equality
$$ J_M^G(\eta, X, f)=\sum_{i\in I} J_M^G(\eta, X_i, f) J_M^G(\eta, X, f_i). $$
\end{coro}

\begin{proof}
By Proposition \ref{howe1}, there exists  a finite subset $\{X_i: i\in I\}\subseteq\sigma$ such that $\{{i_r^\fs}^\ast\circ\delta_M^G(X_i):i\in I\}$ is a basis of ${i_r^\fs}^\ast\circ\delta_M^G(\BC[\sigma])$. By linear algebra, there exists a finite subset $\{f_i: i\in I\}\subseteq\CC_c^\infty(\fs(F)/r)$ such that ${i_r^\fs}^\ast\circ\delta_M^G(X_i)(f_j)=\delta_{ij}, \forall i,j\in I$, where $\delta_{ij}$ denotes the Kronecker delta function. Choose such $\{X_i: i\in I\}$ and $\{f_i: i\in I\}$. 

Then, for all $X\in\sigma$, there exists $\lambda_i, \forall i\in I$ such that
$$ J_M^G(\eta, X, \cdot)=\sum_{i\in I} \lambda_i J_M^G(\eta, X_i, \cdot)\in\CC_c^\infty(\fs(F)/r)^\ast. $$
Hence, for all $i\in I$, 
$$ J_M^G(\eta, X, f_i)=\sum_{j\in I} \lambda_j J_M^G(\eta, X_j, f_i)=\sum_{j\in I} \lambda_j\delta_{ji}=\lambda_i. $$
We have finished the proof. 
\end{proof}

\subsection{The case of $(G',H')$}

For an open compact subgroup $r'$ of $\fs'(F)$, denote by $\CC_c^\infty(\fs'(F)/r')$ the subspace of $\CC_c^\infty(\fs'(F))$ consisting of the functions invariant by translation of $r'$. Let ${i_{r'}^{\fs'}}^*: \CC_c^\infty(\fs'(F))^*\ra\CC_c^\infty(\fs'(F)/r')^*$ be the dual map of the natural injection $\CC_c^\infty(\fs'(F)/r')\hookrightarrow\CC_c^\infty(\fs'(F))$. 

For any set $\sigma'$, denote by $\BC[\sigma']$ the $\BC$-linear space with a basis $\sigma'$. For $M'\in\msl^{H'}(M'_0)$ and $\sigma'\subseteq (\wt{\fm'}\cap\fs'_\rs)(F)$, we define the linear map
\[\begin{split}
\delta_{M'}^{H'}: \BC[\sigma'] &\ra \CC_c^\infty(\fs'(F))^*, \\
(z_Y)_{Y\in\sigma'} &\mapsto \sum_{Y\in\sigma'} z_Y J_{M'}^{H'}(Y, \cdot), \\
\end{split}\]
where $z_Y\in\BC$ is the coordinate at $Y\in\sigma'$. 

For $L'\in\msl^{H'}(M'_0)$ and $r'_{L'}\subseteq(\wt{\fl'}\cap\fs')(F)$ an open compact subgroup, we define ${i_{r'_{L'}}^{\wt{\fl'}\cap\fs'}}^*$ and $\delta_{M'}^{L'}$ in the same way by using $(\wt{L'}, L', \Ad(\alpha))$ in place of $(G', H', \Ad(\alpha))$. 

\begin{prop}[Howe's finiteness]\label{howe2}
Let $r'$ be an open compact subgroup of $\fs'(F)$, $M'\in\msl^{H'}(M'_0)$ and $\sigma'\subseteq (\wt{\fm'}\cap\fs'_\rs)(F)$. Suppose that there exists a compact subset $\sigma'_0\subseteq(\wt{\fm'}\cap\fs')(F)$ such that $\sigma'\subseteq \Ad(M'(F))(\sigma'_0)$. Then the image of the linear map
$$ {i_{r'}^{\fs'}}^*\circ\delta_{M'}^{H'}: \BC[\sigma']\ra\CC_c^\infty(\fs'(F)/r')^* $$
is of finite dimension. 
\end{prop}

The rest of this section is devoted to the proof of Proposition \ref{howe2}. It is similar to the proof of Proposition \ref{howe1} and we only point out some additional argument. Recall that we have chosen the standard maximal compact subgroup $K_{H'}=GL_n(\CO_{D'})$ of $H'(F)=GL_n(D')$ in \textbf{Case I} (resp. $K_{H'}=GL_{\frac{n}{2}}(\CO_{D\otimes_F E})$ of $H'(F)=GL_{\frac{n}{2}}(D\otimes_F E)$ in \textbf{Case II}). Since $\Ad(\tau)(D')=D'$ in \textbf{Case I} (resp. $\Ad(\tau)(D\otimes_F E)=D\otimes_F E$ in \textbf{Case II}), we deduce that $\Ad(\tau)(\CO_{D'})=\CO_{D'}$ in \textbf{Case I} (resp. $\Ad(\tau)(\CO_{D\otimes_F E})=\CO_{D\otimes_F E}$ in \textbf{Case II}). Thus $\Ad(\tau)(K_{H'})=K_{H'}$. Let $\fh'(\CO_F):=\fg\fl_n(\CO_{D'})$ in \textbf{Case I} (resp. $\fh'(\CO_F):=\fg\fl_{\frac{n}{2}}(\CO_{D\otimes_F E})$ in \textbf{Case II}). Let $\fs'(\CO_F):=\fh'(\CO_F)\tau=\tau\fh'(\CO_F)$ be an $\CO_F$-lattice in $\fs'(F)$ (see Section \ref{sectsympar2} for the choice of $\tau$). Let $k':=\fh'(\CO_F)\oplus\fs'(\CO_F)$, whose decomposition is stable under the adjoint action of $K_{H'}$ because $\Ad(\tau)(K_{H'})=K_{H'}$. For all $P'\in\msf^{H'}(M'_0)$, we fix $a_{P'}\in A_{P'}(F)$ such that $|\alpha(a_{P'})|_F<1, \forall \alpha\in\Delta_{P'}^{H'}$. 

Starting from $X\in\CN^{\fs'}\cap(\fs'(\CO_F)-\varpi \fs'(\CO_F))$, we obtain a group homomorphism $\varphi: SL_2(F)\ra G'(F)$ by the Jacobson-Morozov theorem for symmetric spaces (Lemma \ref{jm}). Denote $a^{\varphi}:=\varphi\mat(\varpi,,,\varpi^{-1})\in H'(F)$. Define $P^X$ as in \eqref{defP^X}, which contains the centraliser $\Cent_{G'}(a^\varphi)$ of $a^\varphi$ in $G'$ as a Levi factor. Since $a^\varphi$ commutes with $\alpha$, by Lemma \ref{intpar2}, $P^X\cap H'$ is a parabolic subgroup of $H'$, which contains the centraliser $\Cent_{H'}(a^\varphi)$ of $a^\varphi$ in $H'$ as a Levi factor. We want to show that there exists $x\in K_{H'}$ such that $\Ad(x)(P^X\cap H')\in\msf^{H'}(M'_0)$ and that $\Ad(x)(P^X)\in\msf^{G'}(M'_{\wt{0}})$. 

\begin{lem}\label{jordanform}
For $Y\in\CN^{\fs'}$, there exists $x\in H'(F)$ such that $\Ad(x)(Y)$ is in the Jordan normal form, i.e., diagonal block matrices with entries in $D'\tau$ in \textbf{Case I} (resp. $(D\otimes_F E)\tau$ in \textbf{Case II}) whose blocks are of the form
$$ \left( \begin{array}{cccc}
  0 & \tau        &           &  \\
     & \ddots & \ddots &  \\
     &           & \ddots & \tau \\
     &           &           & 0  \\
  \end{array} \right). $$ 
\end{lem}

\begin{proof}
It can be proved in the same way as \cite[Lemmas 2.2 and 2.3]{MR1487565} by linear algebra over a division ring. 
\end{proof}

Thanks to Lemma \ref{jordanform}, we can construct explicitly the above morphism $\varphi$ (see \cite[p. 184]{MR2838836}). If $X$ is in the Jordan normal form, by {\it{loc. cit.}}, we may choose $a^\varphi\in A_{L'}(F)$ for some $L'\in\msl^{H'}(M'_0)$ such that $\Cent_{H'}(a^\varphi)=L'$ and that $\Cent_{G'}(a^\varphi)=\wt{L'}$. For a general $X$ as above, by Lemma \ref{jordanform}, there exists $y\in H'(F)$ such that $\Ad(y)(a^\varphi)\in A_{L'}(F)$ for some $L'\in\msl^{H'}(M'_0)$ satisfying $\Cent_{H'}(\Ad(y)(a^\varphi))=L'$ and $\Cent_{G'}(\Ad(y)(a^\varphi))=\wt{L'}$. Let $x\in K_{H'}$ be such that $x^{-1}y\in(P^X\cap H')(F)$. Then $\Ad(x)(P^X\cap H')=\Ad(y)(P^X\cap H')$ contains $L'$ as a Levi factor and $\Ad(x)(P^X)=\Ad(y)(P^X)$ contains $\wt{L'}$ as a Levi factor. Furthermore, since $\Ad(x)(P^X)\cap H'=\Ad(x)(P^X\cap H')$, we see that $\Ad(x)(P^X\cap H')$ and $\Ad(x)(P^X)$ are associated under the bijection $P'\mapsto\wt{P'}$ between $\msf^{H'}(M'_0)$ and $\msf^{G'}(M'_{\wt{0}})$. 

Fix $x\in K_{H'}$ as above and denote $P':=\Ad(x)(P^X \cap H')\in\msf^{H'}(M'_0)$. Then $\wt{P'}=\Ad(x)(P^X)$. Put $a^X:=\Ad(x^{-1})(a_{P'})\in H'(F)$. By the argument of \cite[(3) in \S IV.3]{MR1344131}, we show that
\begin{equation}\label{equIV.3(3)2}
 \Ad(a^X)(X)\in\varpi\fs'(\CO_F). 
\end{equation}

\begin{proof}[Proof of Proposition \ref{howe2}]
We may apply the argument of Proposition \ref{howe1} with obvious modifications. Especially, one needs to use Proposition \ref{propwoi2} and \eqref{equIV.3(3)2} to show an analogue of Lemma \ref{lemIV.5} for the case of $(G',H')$. Additionally, to prove an analogue of Lemma \ref{lemIV.3} for this case, one may resort to the argument rather than the consequence of some steps in the proof of \cite[Lemme IV.3]{MR1344131} since our definition of $k'$ is different from $\fg'(\CO_F):=\fg\fl_n(\CO_D)$. However, there is no essential difficulty with our preparation above and we omit details here. 
\end{proof}

\begin{coro}[cf. Corollary \ref{corhowe1}]\label{corhowe2}
Let $r'$ be an open compact subgroup of $\fs'(F)$, $M'\in\msl^{H'}(M'_0)$ and $\sigma'\subseteq (\wt{\fm'}\cap\fs'_\rs)(F)$. Suppose that there exists a compact subset $\sigma'_0\subseteq(\wt{\fm'}\cap\fs')(F)$ such that $\sigma'\subseteq \Ad(M'(F))(\sigma'_0)$. Then there exists a finite subset $\{Y_i: i\in I\}\subseteq\sigma'$ and a finite subset $\{f_i: i\in I\}\subseteq\CC_c^\infty(\fs'(F)/r')$ such that for all $Y\in\sigma'$ and all $f'\in\CC_c^\infty(\fs'(F)/r')$, we have the equality
$$ J_{M'}^{H'}(Y, f')=\sum_{i\in I} J_{M'}^{H'}(Y_i, f') J_{M'}^{H'}(Y, f'_i). $$
\end{coro}


\section{\textbf{Representability of the Fourier transform of weighted orbital integrals}}\label{secrep}

\subsection{The case of $(G,H)$}

Following \cite[\S V.6]{MR1344131}, we denote by $\mse^\fs$ the space of functions $e: \fs_\rs(F)\ra\BC$ such that

(1) $e$ is locally constant; 

(2) for any open compact subset $r$ of $\fs(F)$, there exist constants $c>0$ and $N\in\BN$ such that for all $X\in r\cap\fs_\rs(F)$, one has the inequality
$$ |e(X)|\leq c\sup\{1, -\log |D^\fs(X)|_F\}^N. $$

If $e\in\mse^\fs$, the function $X\mapsto |D^\fs(X)|_F^{-1/2}e(X)$ is locally integrable on $\fs(F)$ thanks to Corollary \ref{corII.1}. It defines then a distribution on $\fs(F)$:
\begin{equation}\label{defV.6.1}
 \forall f\in\CC_c^\infty(\fs(F)), f\mapsto\int_{\fs(F)} f(X)e(X)|D^\fs(X)|_F^{-1/2}dX. 
\end{equation}
Denote by $\msd^\fs$ the space of distributions obtained in this way. The map $\mse^\fs\ra\msd^\fs$ defined above is an isomorphism. For $d\in\msd^\fs$, we shall always denote by $e_d$ its preimage in $\mse^\fs$. 

Notice that the notion $\mse^\fs$ can be defined for any symmetric pair, and that the definition $\msd^\fs$ can at least be extended to the symmetric pair $(M, M_H, \Ad(\epsilon))$, where $M\in\msl^{G,\omega}(M_0)$, since it appears as the product of some copies of the form $(G, H, \Ad(\epsilon))$ in lower dimensions. 

If $d\in\msd^\fs$ is $\eta(\Nrd(\cdot))$-invariant with respect to the adjoint action of $H(F)$, then so is $e_d\in\mse^\fs$ and by the Weyl integration formula (Proposition \ref{wif2.1}), we have the equality
\begin{equation}\label{equV.6inv1}
 d(f)=\sum_{M\in\msl^{G,\omega}(M_0)} |W_0^{M_n}| |W_0^{GL_n}|^{-1} \sum_{\fc\in\mst_\el(\fm\cap\fs)} |W(M_H, \fc)|^{-1} \int_{\fc_\reg(F)} J_G^G(\eta, X, f) e_d(X) dX 
\end{equation}
for all $f\in\CC_c^\infty(\fs(F))$, where $J_G^G(\eta, X, f)$ is defined by (\ref{defwoi1}). 

\begin{remark}[Glueing]\label{gluermk}
Let $d\in\CC_c^\infty(\fs(F))^*$ and $(r_i)_{i\in I}$ be a family of open compact subsets of $\fs(F)$ such that $\bigcup\limits_{i\in I} r_i=\fs(F)$. Suppose that for all $i\in I$, there exists $d_i\in\msd^\fs$ such that $d(f)=d_i(f)$ for all $f\in\CC_c^\infty(\fs(F))$ with $\Supp(f)\subseteq r_i$. Then $d\in\msd^\fs$. Refer to \cite[Remarque V.6]{MR1344131} for the details. 
\end{remark}

Let $M\in\msl^{G,\omega}(M_0)$ and $X\in(\fm\cap\fs_\rs)(F)$. Denote by $\hat{J}_M^G(\eta, X, \cdot)$ the distribution on $\fs(F)$ defined by
$$ \hat{J}_M^G(\eta, X, f):=J_M^G(\eta, X, \hat{f}) $$
for all $f\in\CC_c^\infty(\fs(F))$, where the right hand side is defined by (\ref{defwoi1}). We also have a similar definition for the symmetric pair $(M, M_H, \Ad(\epsilon))$, where $M\in\msl^{G,\omega}(M_0)$. The main result of this section is the following. 

\begin{prop}[Representability]\label{repr1}
Let $M\in\msl^{G,\omega}(M_0)$ and $X\in(\fm\cap\fs_\rs)(F)$. Then the distribution $\hat{J}_M^G(\eta, X, \cdot)\in\msd^\fs$. 
\end{prop}

\begin{remark}
For $M=G$, Proposition \ref{repr1} is essentially \cite[Theorem 6.1.(i)]{MR3414387} (see also \cite[Theorem 6.2]{MR3299843}). 
\end{remark}

The rest of this section is devoted to the proof of Proposition \ref{repr1}. We shall follow the main steps in \cite[\S V.7-10]{MR1344131}. 

Let $\fc$ be a Cartan subspace of $\fs$. Recall that $T_\fc$ denotes the centraliser of $\fc$ in $H$. Suppose that $e_0: (T_\fc(F)\bs H(F))\times \fc_\reg(F)\ra\BC$ is a function such that

(1) $e_0$ is locally constant; 

(2) for any open compact subset $r$ of $\fs(F)$, there exist constants $c>0$ and $N\in\BN$ such that for all $x\in T_\fc(F)\bs H(F)$ and $X\in\fc_\reg(F)$ satisfying $\Ad(x^{-1})(X)\in r$, one has the inequality
$$ e_0(x, X)\leq c\sup\{1, -\log |D^\fs(X)|_F\}^N. $$
Following \cite[\S V.7]{MR1344131}, for $f\in\CC_c^\infty(\fs(F))$, we define
\begin{equation}\label{defV.7.1}
 d_0(f):=\int_{\fc_\reg(F)} |D^\fs(X)|_F^{1/2} \int_{T_\fc(F)\bs H(F)} f(\Ad(x^{-1})(X)) e_0(x, X) dx dX. 
\end{equation}

\begin{lem}\label{lemV.7.1}
Let $\fc$ be a Cartan subspace of $\fs$. Suppose that $e_0$ satisfies  the above hypotheses. Then the integral (\ref{defV.7.1}) is absolutely convergent. Moreover, the distribution $d_0\in\msd^\fs$. 
\end{lem}

\begin{proof}
We define a function $e':\fs_\rs(F)\ra \BC$ by
\begin{equation}\label{defe'}
 e'(X):=\sum_{\{x\in T_\fc(F)\bs H(F):\Ad(x)(X)\in \fc(F)\}} e_0(x,\Ad(x)(X)) 
\end{equation}
for all $X\in\fs_\rs(F)$. If $X\notin\Ad(H(F))(\fc_\reg(F))$, then $e'(X)=0$. If $X\in\Ad(H(F))(\fc_\reg(F))$, then the sum in (\ref{defe'}) is actually over the finite set $W(H, \fc)y$, where $y\in H(F)$ is any element such that $\Ad(y)(X)\in\fc(F)$. Hence, $e'$ is well-defined. 

Additionally, one may check that $e'\in\mse^\fs$ from the hypotheses on $e_0$. Let $d'\in\msd^\fs$ be the distribution associated to $e'$ by (\ref{defV.6.1}). For all $f\in\CC_c^\infty(\fs(F))$, by the Weyl integration formula (\ref{wif}), we have
\begin{equation}\label{equlemV.7}
 d'(f)=|W(H, \fc)|^{-1} \int_{\fc_\reg(F)} |D^\fs(X)|_F^{1/2} \int_{T_\fc(F)\bs H(F)} f(\Ad(x^{-1})(X)) e'(\Ad(x^{-1})(X)) dx dX. 
\end{equation}
Since
$$ e'(\Ad(x^{-1})(X))=\sum_{w\in W(H,\fc)} e_0(wx,\Ad(w)(X)) $$
for all $x\in T_\fc(F)\bs H(F)$ and all $X\in\fc_\reg(F)$, we deduce that
$$ d'(f)=|W(H, \fc)|^{-1} \sum_{w\in W(H,\fc)} \int_{\fc_\reg(F)} |D^\fs(X)|_F^{1/2} \int_{T_\fc(F)\bs H(F)} f(\Ad(x^{-1})(X)) e_0(wx,\Ad(w)(X)) dx dX. $$
Applying the change of variables $X\mapsto \Ad(w^{-1})(X)$ and $x\mapsto w^{-1}x$, which does not modify the Haar measures, we obtain
$$ d'(f)=d_0(f) $$
for all $f\in\CC_c^\infty(\fs(F))$. That is to say, $d_0=d'\in\msd^\fs$. 

Note that in the argument above, we have used the convergence of an analogue of (\ref{equlemV.7}) with $e_0$ and $f$ replaced by their absolute values. It also results in the absolute convergence of (\ref{defV.7.1}). 
\end{proof}

\begin{coro}[Parabolic induction]\label{corV.8.1}
Let $M\in\msl^{G,\omega}(M_0)$, $P\in\msp^G(M)$ and $d\in\msd^{\fm\cap\fs}$. Then the distribution on $\fs(F)$ defined by $\forall f\in\CC_c^\infty(\fs(F)), f\mapsto d(f_P^\eta)$ belongs to $\msd^\fs$, where $f_P^\eta$ is defined by (\ref{pardes1}). 
\end{coro}

\begin{proof}
Applying the Weyl integration formula (\ref{wif}) to $d(f_P^\eta)$, we see that it suffices to fix a Cartan subspace $\fc\subseteq\fm\cap\fs$ and prove that the distribution on $\fs(F)$ defined by
\begin{equation}\label{equ0corV.8}
 \forall f\in\CC_c^\infty(\fs(F)), f\mapsto \int_{\fc_\reg(F)} |D^{\fm\cap\fs}(X)|_F^{1/2} \int_{T_\fc(F)\bs M_H(F)} f_P^\eta(\Ad(x^{-1})(X)) e_d(\Ad(x^{-1})(X)) dx dX 
\end{equation}
belongs to $\msd^\fs$. Recall that $e_d\in\mse^{\fm\cap\fs}$ is associated to $d$ by (\ref{defV.6.1}). 

Define a function $e_1: H(F)\times\fc_\reg(F)\ra\BC$ by
\begin{equation}\label{defe_1}
 e_1(x,X):=\eta(\Nrd(x))\int_{M_H(F)\cap K} e_d(\Ad(m_P(x)k)^{-1}(X)) \eta(\Nrd(m_P(x)k)) dk 
\end{equation}
for all $x\in H(F)$ and all $X\in\fc_\reg(F)$, where $m_P(x)\in M_H(F)$ is any element such that $m_P(x)^{-1}x\in N_{P_H}(F)K_H$. 
Since $K_H$ is a special subgroup of $H(F)$, such $m_P(x)$ is well-defined modulo $M_H(F)\cap K$. Hence, the integral (\ref{defe_1}) is independent of the choice of $m_P(x)$. 
The function $e_1$ is left $T_\fc(F)$-invariant on the first variable, so it induces a function (still denoted by $e_1$) : $T_\fc(F)\bs H(F)\times\fc_\reg(F)\ra\BC$. 

We shall check that $e_1$ verifies the hypotheses of Lemma \ref{lemV.7.1}. Firstly, $e_1$ is locally constant because $e_d$ is locally constant and $e_1$ is right $K_H$-invariant on the first variable. Secondly, suppose that $r$ is an open compact subset of $\fs(F)$. We fix an open compact subset $r_M\subseteq(\fm\cap\fs)(F)$ such that if $X\in(\fm\cap\fs)(F), U\in(\fn_P\cap\fs)(F), k\in K_H$ satisfy $\Ad(k)(X+U)\in r$, then $X\in r_M$; this is possible for it suffices to let $r_M$ contain the projection of $\Ad(K_H)(r)$ to $(\fm\cap\fs)(F)$. Replacing $r_M$ with $\Ad(M_H(F)\cap K)(r_M)$ if necessary, we may additionally assume that
\begin{equation}\label{equ1corV.8}
 \Ad(M_H(F)\cap K)(r_M)=r_M. 
\end{equation}
Since $e_d\in\mse^{\fm\cap\fs}$, there exist constants $c>0$ and $N\in\BN$ such that
$$ |e_d(X)|\leq c\sup\{1, -\log |D^{\fm\cap\fs}(X)|_F\}^N $$
for all $X\in r_M\cap(\fm\cap\fs)_\rs(F)$. 
One sees from (\ref{weyldisc}) that for all $X\in(\fm\cap\fs)_\rs(F)$, 
$$ |D^\fs(X)|_F|D^{\fm\cap\fs}(X)|_F^{-1}=|\det(\ad(X)|_{\fg/\fm})|_F^{1/2}. $$
Hence, $|D^\fs(X)|_F|D^{\fm\cap\fs}(X)|_F^{-1}$ is bounded for $X\in r_M\cap(\fm\cap\fs)_\rs(F)$. 
Recall that $\fm\cap\fs_\rs\subseteq(\fm\cap\fs)_\rs$ (see Section \ref{sectsympar1}) and that $|D^\fs(X)|_F\neq0$ for $X\in\fs_\rs(F)$ (see Section \ref{generalcases}). 
We deduce that there exists $c'>0$ such that
\begin{equation}\label{equ2corV.8}
 |e_d(X)|\leq c'\sup\{1, -\log |D^\fs(X)|_F\}^N 
\end{equation}
for all $X\in r_M\cap\fs_\rs(F)$. Now, suppose that $x\in T_\fc(F)\bs H(F)$ and $X\in\fc_\reg(F)$ satisfy $\Ad(x^{-1})(X)\in r$. Write $x=mnk$ with $m\in M_H(F), n\in N_{P_H}(F)$ and $k\in K_H$. Then
$$ \Ad(x^{-1})(X)=\Ad(k^{-1})(\Ad(m^{-1})(X)+U), $$
where $U:=\Ad(n^{-1}m^{-1})(X)-\Ad(m^{-1})(X)\in(\fn_P\cap\fs)(F)$. Thus $\Ad(m^{-1})(X)\in r_M$ by our assumption on $r_M$. Thanks to (\ref{equ1corV.8}) and (\ref{equ2corV.8}), we obtain
$$ |e_1(x,X)|\leq c'\sup\{1, -\log |D^\fs(X)|_F\}^N. $$
To sum up, $e_1$ verifies the hypotheses of Lemma \ref{lemV.7.1}. 

Applying Lemma \ref{lemV.7.1} to $\fc$ and $e_1$, we know that the distribution $d_1$ on $\fs(F)$ defined by
$$ \forall f\in\CC_c^\infty(\fs(F)), d_1(f):=\int_{\fc_\reg(F)} |D^\fs(X)|_F^{1/2} \int_{T_\fc(F)\bs H(F)} f(\Ad(x^{-1})(X)) e_1(x, X) dx dX $$
belongs to $\msd^\fs$. Note that $e_1(mnk, X)=\eta(\Nrd(k))e_1(m, X)$ for $m\in M_H(F), n\in N_{P_H}(F), k\in K_H$. By the same argument as the proof of Proposition \ref{propwoi1}.4), one shows that
\[\begin{split} 
d_1(f)=&\int_{\fc_\reg(F)} |D^{\fm\cap\fs}(X)|_F^{1/2} \int_{T_\fc(F)\bs M_H(F)} f_P^\eta(\Ad(m^{-1})(X)) e_1(m,X) dm dX \\
=&\int_{\fc_\reg(F)} |D^{\fm\cap\fs}(X)|_F^{1/2} \int_{T_\fc(F)\bs M_H(F)} f_P^\eta(\Ad(m^{-1})(X)) \int_{M_H(F)\cap K} e_d(\Ad(mk)^{-1}(X))\eta(\Nrd(k)) \\ 
&dk dm dX. 
\end{split}\]
Note that for $k\in M_H(F)\cap K$, we have $\Ad(k^{-1})f_P^\eta=\eta(\Nrd(k))f_P^\eta$. By the change of variables $mk\mapsto m$, one can eliminate the integral over $M_H(F)\cap K$ and see that $d_1$ is the same as (\ref{equ0corV.8}). 
\end{proof}

Let $M\in\msl^{G,\omega}(M_0)$ and $d\in\msd^{\fm\cap\fs}$. Suppose that $d$ is $\eta(\Nrd(\cdot))$-invariant with respect to the adjoint action of $M_H(F)$. Following \cite[\S V.9]{MR1344131}, we define a distribution $\Ind_M^{G, w}(d)$ on $\fs(F)$ by
\begin{equation}\label{wtind1}
 \Ind_M^{G, w}(d)(f):=\sum_{\{L\in\msl^{G,\omega}(M_0): L\subseteq M\}} |W_0^{L_n}| |W_0^{M_n}|^{-1} \sum_{\fc\in\mst_\el(\fl\cap\fs)} |W(L_H, \fc)|^{-1} \int_{\fc_\reg(F)} J_M^G(\eta, X, f) e_d(X) dX
\end{equation}
for all $f\in\CC_c^\infty(\fs(F))$, where $J_M^G(\eta, X, f)$ is defined by (\ref{defwoi1}). In particular, if $M=G$ and $d\in\msd^\fs$ is $\eta(\Nrd(\cdot))$-invariant with respect to the adjoint action of $H(F)$, we have $\Ind_G^{G, w}(d)=d$ by (\ref{equV.6inv1}). 

\begin{coro}\label{corV.9.1}
Let $M\in\msl^{G,\omega}(M_0)$ and $d\in\msd^{\fm\cap\fs}$. Suppose that $d$ is $\eta(\Nrd(\cdot))$-invariant with respect to the adjoint action of $M_H(F)$. Then the integral (\ref{wtind1}) is absolutely convergent. Moreover, the distribution $\Ind_M^{G, w}(d)\in\msd^\fs$. 
\end{coro}

\begin{remark}
This corollary is unnecessary for the proof of Proposition \ref{repr1} but useful in Section \ref{sectinvwoi1}. 
\end{remark}

\begin{proof}[Proof of Corollary \ref{corV.9.1}]
It suffices to fix a Cartan subspace $\fc\subseteq\fm\cap\fs$ and prove the same assertion for the distribution on $\fs(F)$ defined by
\begin{equation}\label{equcorV.9}
 \forall f\in\CC_c^\infty(\fs(F)), f\mapsto\int_{\fc_\reg(F)} J_M^G(\eta, X, f) e_d(X) dX. 
\end{equation}

Define a function $e_2: (T_\fc(F)\bs H(F))\times\fc_\reg(F)\ra\BC$ by
\begin{equation}\label{equpfcorV.9}
 e_2(x,X):=\eta(\Nrd(x))v_M^G(x)e_d(X) 
\end{equation}
for all $x\in T_\fc(F)\bs H(F)$ and all $X\in\fc_\reg(F)$. It is locally constant. Note that $e_d(\Ad(m^{-1})(X))=\eta(\Nrd(m))e_d(X)$ for $x\in M_H(F)$ and $X\in(\fm\cap\fs_\rs)(F)$ by our assumption on $d$. Thus we may use the same argument as in the proof of Corollary \ref{corV.8.1} to show the inequality
$$ |e_d(X)|=|e_d(\Ad(m^{-1})(X))|\leq c'\sup\{1, -\log |D^\fs(X)|_F\}^N. $$
Thanks to Lemma \ref{lemIII.5.1}, one has a similar bound for $v_M^G(x)$. In sum, $e_2$ verifies the hypotheses of Lemma \ref{lemV.7.1}. 

Applying Lemma \ref{lemV.7.1} to $\fc$ and $e_2$, we know that the integral
$$ \forall f\in\CC_c^\infty(\fs(F)), f\mapsto\int_{\fc_\reg(F)} |D^\fs(X)|_F^{1/2} \int_{T_\fc(F)\bs H(F)} f(\Ad(x^{-1})(X)) e_2(x, X) dx dX $$
is absolutely convergent and defines a distribution in $\msd^\fs$. This distribution is the same as (\ref{equcorV.9}). 
\end{proof}

\begin{remark}
Instead of the hypothesis $e_d\in\mse^{\fm\cap\fs}$, an analogue of Corollary \ref{corV.9.1} holds if one assumes that $d$ is defined via (\ref{defV.6.1}) by a function $e_d:(\fm\cap\fs_\rs)(F)\ra\BC$ locally constant, $\eta(\Nrd(\cdot))$-invariant with respect to the adjoint action of $M_H(F)$, and such that for any open compact subset $r$ of $(\fm\cap\fs)(F)$, there exist constants $c>0$ and $N\in\BN$ such that for all $X\in r\cap\fs_\rs(F)$, one has the inequality
$$ |e_d(X)|\leq c\sup\{1, -\log |D^\fs(X)|_F\}^N. $$
The proof is the same, except that we need not use the boundedness of $|D^\fs(X)|_F|D^{\fm\cap\fs}(X)|_F^{-1}$ for $X\in r_M\cap\fs_\rs(F)$ (see the proof of Corollary \ref{corV.8.1}). 
\end{remark}

\begin{lem}\label{countlemforrepr1}
Let $M\in\msl^{G,\omega}(M_0)$ and $\fc\subseteq\fm\cap\fs$ be an $M$-elliptic Cartan subspace. 

1) Let $M'\in\msl^{G,\omega}(M_0), \fc'\in\mst_\el(\fm'\cap\fs)$ and $x\in H(F)$ be such that $\Ad(x)(\fc)=\fc'$. Then there exist elements $m'\in M'_H(F)$ and $w\in\Norm_{H(F)}(M_0)$ such that $x=m'w$. 

2) The cardinality of
$$ \{(M',\fc'): M'\in\msl^{G,\omega}(M_0), \fc'\in\mst_\el(\fm'\cap\fs), \fc'\text{ is }H(F)\text{-conjugate to }\fc\} $$
is
$$ |W_0^{GL_n}||W_0^{M_n}|^{-1}|W(M_H,\fc)||W(H,\fc)|^{-1}. $$
\end{lem}

\begin{proof}
1) Since $\fc\subseteq\fm\cap\fs$ (resp. $\fc'\subseteq\fm'\cap\fs$) is $M$-elliptic (resp. $M'$-elliptic), we have $A_{T_\fc}=A_M$ (resp. $A_{T_{\fc'}}=A_{M'}$). From $\Ad(x)(\fc)=\fc'$, we obtain $\Ad(x)(A_{T_\fc})=A_{T_{\fc'}}$, $\Ad(x)(A_M)=A_{M'}$ and $\Ad(x)(M)=M'$, which implies that $\Ad(x)(M_H)=M'_H$. Then 1) can be shown by the same argument as in the proof of Lemma \ref{countlemforwif2.1}.1). 

2) We can and shall identify an $M$-elliptic (resp. $M'$-elliptic) Cartan subspace in $\fm\cap\fs$ (resp. $\fm'\cap\fs$) with its $M_H(F)$(resp. $M'_H(F)$)-conjugacy class. Then $\mst_\el(\fm'\cap\fs)$ is identified to the set of $M'_H(F)$-conjugacy classes of $M'$-elliptic Cantan subspaces in $\fm'\cap\fs$. As in the proof of Lemma \ref{countlemforwif2.1}.1), we also see that the group $W_{0}^{H, \omega}:=\left\{\mat(w_n,,,w_n): w_n\in W_0^{GL_{n,D}}\right\}$ acts transitively on the set of pairs in 2). 

Firstly, let us count $M'$ appearing in the pairs (cf. \cite[p. 426]{MR2192014}). 

Since $M\in\msl^{G,\omega}(M_0)$, for $w\in W_0^{H,\omega}$, we see that $w=\mat(w_n,,,w_n)\in\Norm_{W_0^{H,\omega}}(M)$ if and only if $w_n\in\Norm_{W_0^{GL_{n,D}}}(M_{n,D})$, where $\Norm_{W_0^{H,\omega}}(M)$ (resp. $\Norm_{W_0^{GL_{n,D}}}(M_{n,D})$) denotes the normaliser of $M$ (resp. $M_{n,D}$) in $W_0^{H,\omega}$ (resp. $W_0^{GL_{n,D}}$). Hence, the number of $M'$ is
$$ |W_0^{GL_{n,D}}||\Norm_{W_0^{GL_{n,D}}}(M_{n,D})|^{-1}. $$ 

Secondly, for such an $M'$ fixed, we count $\fc'$ such that $(M',\fc')$ belongs to the set of pairs in 2) (cf. \cite[Lemma 7.1]{MR2192014}). 

For $x\in H(F)$, we claim that $\Ad(x)(\fc')\subseteq\fm'\cap\fs$ if and only if $x\in\Norm_{H(F)}(M')$, which denotes the normaliser of $M'$ in $H(F)$. On the one hand, suppose that $\Ad(x)(\fc')\subseteq\fm'\cap\fs$. Then $A_{M'}\subseteq \Cent_{H}(\fm'\cap\fs)\subseteq \Ad(x)(T_{\fc'})$, where $\Cent_{H}(\fm'\cap\fs)$ denotes the centraliser of $\fm'\cap\fs$ in $H$. Since $\Ad(x)(A_{M'})=\Ad(x)(A_{T_{\fc'}})$ is the maximal $F$-split torus in $\Ad(x)(T_{\fc'})$, we have $A_{M'}\subseteq \Ad(x)(A_{M'})$. By comparison of dimensions, we deduce that $\Ad(x)(A_{M'})=A_{M'}$, so $x\in\Norm_{H(F)}(M')$. On the other hand, suppose that $x\in\Norm_{H(F)}(M')$. Since $x\in H(F)$, we have $\Ad(x)(\fm'\cap\fs)=(\fm'\cap\fs)$. But $\fc'\subseteq\fm'\cap\fs$, so we obtain $\Ad(x)(\fc')\subseteq\fm'\cap\fs$. In sum, we have proved our claim. 

From this claim, the number of $\fc'$ is
$$ |M'_H(F)\bs \Norm_{H(F)}(M')/\Norm_{H(F)}(\fc')|. $$
Since $M'_H(F)$ is a normal subgroup of $\Norm_{H(F)}(M')$, we know that the number of the double cosets is equal to
\[\begin{split}
 &|\Norm_{H(F)}(M')/M'_H(F)||\Norm_{H(F)}(\fc')/(\Norm_{H(F)}(\fc')\cap M'_H(F))|^{-1} \\
=&|\Norm_{H(F)}(M')/M'_H(F)||\Norm_{H(F)}(\fc')/\Norm_{M'_H(F)}(\fc')|^{-1}. 
\end{split}\]
For $x\in\Norm_{H(F)}(M')$, we have $\Ad(x)(A_{M'})=A_{M'}$. Because $M'\in\msl^{G,\omega}(M_0)$, there exists $w\in W_0^{H,\omega}$ such that $w^{-1}x\in\Cent_{H(F)}(A_{M'})=M'_H(F)$, where $\Cent_{H(F)}(A_{M'})$ denotes the centraliser of $A_{M'}$ in $H(F)$. Since $x\in\Norm_{H(F)}(M')$, we have $w\in\Norm_{W_0^{H,\omega}}(M')$. That is to say, 
$$ \Norm_{H(F)}(M')=\Norm_{W_0^{H,\omega}}(M') M'_H(F). $$
Therefore, 
\[\begin{split}
 |\Norm_{H(F)}(M')/M'_H(F)|=&|\Norm_{W_0^{H,\omega}}(M')||(\Norm_{W_0^{H,\omega}}(M')\cap M'_H(F))|^{-1} \\
=&|\Norm_{W_0^{GL_{n,D}}}(M'_{n,D})||(W_0^{H,\omega}\cap M'_H(F))|^{-1} \\ 
=&|\Norm_{W_0^{GL_{n,D}}}(M'_{n,D})||W_0^{M'_{n,D}}|^{-1}. 
\end{split}\]
Since $M'$ and $M$ are $W_0^{H,\omega}$-conjugate, $M'_{n,D}$ and $M_{n,D}$ are $W_0^{GL_{n,D}}$-conjugate. Hence, 
$$ |\Norm_{H(F)}(M')/M'_H(F)|=|\Norm_{W_0^{GL_{n,D}}}(M_{n,D})||W_0^{M_{n,D}}|^{-1}. $$
We also have
\[\begin{split}
 |\Norm_{H(F)}(\fc')/\Norm_{M'_H(F)}(\fc')|=&|\Norm_{H(F)}(\fc')/T_{\fc'}(F)||\Norm_{M'_H(F)}(\fc')/T_{\fc'}(F)|^{-1} \\
=&|W(H,\fc')||W(M'_H,\fc')|^{-1}. 
\end{split}\]
Since $(M',\fc')$ and $(M,\fc)$ are $W_0^{H,\omega}$-conjugate, we obtain
$$ |\Norm_{H(F)}(\fc')/\Norm_{M'_H(F)}(\fc')|=|W(H,\fc)||W(M_H,\fc)|^{-1}. $$
To sum up, the number of $\fc'$ is
$$ |\Norm_{W_0^{GL_{n,D}}}(M_{n,D})||W_0^{M_{n,D}}|^{-1}|W(M_H,\fc)||W(H,\fc)|^{-1}. $$

Finally, combining the numbers of $M'$ and $\fc'$, we obtain the number of pairs $(M',\fc')$ in 2).
\end{proof}

\begin{proof}[Proof of Proposition \ref{repr1}]
First of all, suppose that $X\in(\fm\cap\fs_\rs)(F)_\el$. Then $\fc:=\fs_X\subseteq\fm\cap\fs$ is an $M$-elliptic Cartan subspace and $X\in\fc_\reg(F)$. Fix an open compact subgroup $r\subseteq\fs(F)$ and set $r^\ast:=\{Y\in\fs(F):\forall Z\in r, \Psi(\langle Y,Z\rangle)=1\}$, which is also an open compact subgroup of $\fs(F)$. For all $L\in\msl^G(M)$, fix an open compact subgroup $r_L\subseteq(\fl\cap\fs)(F)$ such that if $Q\in\msp^G(L)$ and if $f\in\CC_c^\infty(\fs(F))$ satisfies $\Supp(f)\subseteq r$, then $\Supp(f_Q^\eta)\subseteq r_L$, where $f_Q^\eta$ is defined by (\ref{pardes1}). Define $r_L^\ast$ in the same way as $r^\ast$. 

There exists a neighbourhood $\sigma$ of $X$ in $\fc_\reg(F)$ such that for all $L\in\msl^G(M)$ and all $f\in\CC_c^\infty((\fl\cap\fs)(F)/r_L^\ast)$, the function $J_M^L(\eta,\cdot,f)$ is constant on $\sigma$. In fact, for $L$ and $f$ fixed, this results from Proposition \ref{propwoi1}.2) (actually its product form is needed). It suffices to apply Howe's finiteness (the product form of Corollary \ref{corhowe1}) to each symmetric pair $(L,L_H,\Ad(\epsilon))$ and an arbitrary compact neighbourhood of $X$ in $\fc_\reg(F)$, and then take the intersection of a finite number of neighbourhoods involved. 

We shall fix a $\sigma$ satisfying the above condition and such that if two elements of $\sigma$ are $H(F)$-conjugate (or equivalently $W(H,\fc)$-conjugate), then they are the same. The latter condition is achievable since the $W(H,\fc)$-conjugates of an element in $\fc_\reg(F)$ form a finite subset, which is discrete. Consider the local isomorphism $\beta: (T_\fc(F)\bs H(F))\times\fc_\reg(F)\ra\fs_\rs(F)$ of $F$-analytic manifolds induced by the adjoint action. Its restriction to $(T_\fc(F)\bs H(F))\times\sigma$ is injective. Choose a neighbourhood $\varepsilon$ of $1$ in $T_\fc(F)\bs H(F)$ such that $\eta(\Nrd(\varepsilon))=1$. The set $\beta(\varepsilon, \sigma)$ is a neighbourhood of $X$ in $\fs(F)$. Fix a function $f'\in\CC_c^\infty(\fs(F))$ such that $\Supp(f')\subseteq\beta(\varepsilon, \sigma), f'\geq0$ and $f'(X)\neq0$. 

Let $f\in\CC_c^\infty(\fs(F))$ with $\Supp(f)\subseteq r$. We shall calculate $J^G(\eta,\hat{f},f')$, which is defined by (\ref{geomnoninvltf1}). 

Consider $M'\in\msl^{G,\omega}(M_0)$ and $\fc'\in\mst_\el(\fm'\cap\fs)$. If $\fc'$ and $\fc$ are not $H(F)$-conjugate, by our choice of $f'$, the function $J_{M'}^G(\eta,\cdot,\hat{f},f')$ vanishes on $\fc'_\reg(F)$. Now suppose that $\fc'$ and $\fc$ are $H(F)$-conjugate. Let $x\in H(F)$ be such that $\Ad(x)(\fc)=\fc'$. By Lemma \ref{countlemforrepr1}.1), there exist elements $m'\in M'_H(F)$ and $w\in\Norm_{H(F)}(M_0)$ such that $x=m'w$. By Proposition \ref{propV.1.1}.4), for $X'\in\fc'_\reg(F)$, we have
\[\begin{split}
 J_{M'}^{G}(\eta, X', \hat{f}, f')&=J_{M'}^{G}(\eta, \Ad({m'}^{-1})(X'), \hat{f}, f') \\
&=J_{\Ad(w)(M)}^{G}(\eta, \Ad(wx^{-1})(X'), \hat{f}, f') \\
&=J_M^G(\eta, \Ad(x^{-1})(X'), \hat{f}, f'). 
\end{split}\]
From our choices of Haar measures, we obtain
$$ \int_{\fc'_\reg(F)} J_{M'}^{G}(\eta, X', \hat{f}, f')=\int_{\fc_\reg(F)} J_M^G(\eta, Y, \hat{f}, f') dY. $$
By Lemma \ref{countlemforrepr1}.2), the number of pairs $(M',\fc')$ with $\fc'$ being $H(F)$-conjugate to $\fc$ is
$$ |W_0^{GL_n}||W_0^{M_n}|^{-1}|W(M_H,\fc)||W(H,\fc)|^{-1}. $$ 
We deduce that
$$ J^G(\eta,\hat{f},f')=(-1)^{\dim(A_M/A_G)} |W(H,\fc)|^{-1} \int_{\fc_\reg(F)}  J_M^G(\eta, Y, \hat{f}, f') dY. $$
It follows from our choice of $f'$ that the support of the restriction to $\fc_\reg(F)$ of the function $J_M^G(\eta, \cdot, \hat{f}, f')$ is contained in $\coprod\limits_{w\in W(H,\fc)} \Ad(w)(\sigma)$. Then
$$ \int_{\fc_\reg(F)}  J_M^G(\eta, Y, \hat{f}, f') dY=\sum_{w\in W(H,\fc)} \int_{\fc_\reg(F)}  J_M^G(\eta, Y, \hat{f}, f') 1_{\Ad(w)(\sigma)}(Y) dY, $$
where $1_{\Ad(w)(\sigma)}$ denotes the characteristic function of $\Ad(w)(\sigma)$. By the change of variables $Y\mapsto\Ad(w)(Y)$, which does not modify the Haar measure, we have
$$ \int_{\fc_\reg(F)}  J_M^G(\eta, Y, \hat{f}, f') 1_{\Ad(w)(\sigma)}(Y) dY=\int_{\sigma} J_M^G(\eta, \Ad(w)(Y), \hat{f}, f') dY. $$
Since $w\in W(H,\fc)$, we have shown above that
$$ J_M^G(\eta, \Ad(w)(Y), \hat{f}, f')=J_M^G(\eta, Y, \hat{f}, f'), $$
which is independent of $w$. Therefore, 
$$ J^G(\eta,\hat{f},f')=(-1)^{\dim(A_M/A_G)} \int_{\sigma} J_M^G(\eta, Y, \hat{f}, f') dY. $$

Let $Y\in\sigma$. Applying the splitting formula for $J_M^G(\eta, Y, \hat{f}, f')$ (Proposition \ref{propV.1.1}.6)), we have 
$$ J_M^G(\eta, Y, \hat{f}, f')=\sum_{L_1,L_2\in\msl^G(M)} d_M^G(L_1,L_2) J_M^{L_1}(\eta,Y,\hat{f}_{\ov{Q_1}}^\eta) J_M^{L_2}(\eta,Y,{f'}_{Q_2}^\eta). $$
For all $Q\in\msf^{G}(M)$, since $\Supp(f_Q^\eta)\subseteq r_{M_Q}$, $\hat{f}_Q^\eta$ is invariant by translation of $r_{M_Q}^\ast$. In particular, $\hat{f}_{\ov{Q_1}}^\eta$ is invariant by $r_{L_1}^\ast$. Then by our assumption on $\sigma$, $J_M^{L_1}(\eta, \cdot, \hat{f}_{\ov{Q_1}}^\eta)$ is constant on $\sigma$ and thus equal to $J_M^{L_1}(\eta, X, \hat{f}_{\ov{Q_1}}^\eta)$. Therefore, 
$$ J^G(\eta,\hat{f},f')=\sum_{L_1,L_2\in\msl^G(M)} c(L_1,L_2) J_M^{L_1}(\eta, X, \hat{f}_{\ov{Q_1}}^\eta), $$
where
$$ c(L_1, L_2):=d_M^G(L_1,L_2) (-1)^{\dim(A_M/A_G)} \int_{\sigma} J_M^{L_2}(\eta,Y,{f'}_{Q_2}^\eta) dY. $$

We claim that $c(G, M)\neq0$. 
In fact, from (1) and (4) in Section \ref{gmmaps}, we have
$$ c(G,M)=(-1)^{\dim(A_M/A_G)} \int_{\sigma} J_M^{M}(\eta,Y,{f'}_{Q_2}^\eta) dY, $$
where $M_{Q_2}=M$. 
By Proposition \ref{propwoi1}.4), we have
$$ J_M^{M}(\eta,Y,{f'}_{Q_2}^\eta)=J_M^{Q_2}(\eta,Y,f'). $$
Since $v_M^{Q_2}=1$, we obtain
$$ J_M^{Q_2}(\eta,Y,f')=J_G^G(\eta,Y,f'). $$
Hence, 
$$ c(G,M)=(-1)^{\dim(A_M/A_G)} \int_{\sigma} |D^\fs(Y)|_F^{1/2} \int_{H_Y(F)\bs H(F)} f'(\Ad(x^{-1})(Y))\eta(\Nrd(x))dx dY. $$
If $\Ad(x^{-1})(Y)\in\Supp(f')\subseteq\beta(\varepsilon, \sigma)$, since the restriction of $\beta$ to $(H_Y(F)\bs H(F))\times\sigma$ is injective, we have $x\in\varepsilon$ and then $\eta(\Nrd(x))=1$. Since $f'\geq0$ and $f'(X)\neq0$, we deduce our claim. Now, because of (3) in Section \ref{gmmaps}, we have 
\begin{equation}\label{equV.10.(2)}
 \hat{J}_M^G(\eta, X, f)=J^G(\eta,\hat{f},f')-\sum_{L_1,L_2\in\msl^G(M), L_1\neq G} c(L_1,L_2) c(G,M)^{-1} \hat{J}_M^{L_1}(\eta, X, f_{\ov{Q_1}}^\eta) 
\end{equation}
for all $f\in\CC_c^\infty(\fs(F))$ with $\Supp(f)\subseteq r$. 

By induction on the dimension of $G$ and parabolic induction (Corollary \ref{corV.8.1}), one can suppose that for all $L_1\in\msl^G(M), L_1\neq G$, the distribution on $\fs(F)$ defined by $\forall f\in\CC_c^\infty(\fs(F)), f\mapsto\hat{J}_M^{L_1}(\eta, X, f_{\ov{Q_1}}^\eta)$ belongs to $\msd^\fs$. This is actually a product form of the induction hypothesis in lower dimensions. 

We claim that the distribution on $\fs(F)$ defined by $\forall f\in\CC_c^\infty(\fs(F)), f\mapsto J^G(\eta,\hat{f},f')$ belongs to $\msd^\fs$ (cf. \cite[(3) in \S V.10]{MR1344131}). In fact, thanks to the noninvariant trace formula (Theorem \ref{noninvltf1}), one can replace $J^G(\eta,\hat{f},f')$ with $J^G(\eta,f,\hat{f'})$. By its definition (\ref{geomnoninvltf1}), it suffices to fix $M'\in\msl^{G,\omega}(M_0), \fc'\in\mst_\el(\fm'\cap\fs)$ and prove that the distribution on $\fs(F)$ defined by
$$ \forall f\in\CC_c^\infty(\fs(F)), f\mapsto \int_{\fc'_\reg(F)} J_{M'}^G(\eta,Y,f,\hat{f'})dY $$
belongs to $\msd^\fs$. By the splitting formula (Proposition \ref{propV.1.1}.6)), it suffices to fix $L'_1,L'_2\in\msl^G(M')$ and prove that the distribution on $\fs(F)$ defined by
$$ \forall f\in\CC_c^\infty(\fs(F)), f\mapsto \int_{\fc'_\reg(F)} J_{M'}^{L'_1}(\eta,Y,f_{\ov{Q'_1}}^\eta) J_{M'}^{L'_2}(\eta,Y,\hat{f'}_{Q'_2}^\eta)dY, $$
where $(Q'_1,Q'_2):=s(L'_1,L'_2)$, belongs to $\msd^\fs$. By Proposition \ref{propwoi1}.4) and the definition (\ref{defwoi1}), we have
$$ J_{M'}^{L'_1}(\eta,Y,f_{\ov{Q'_1}}^\eta)=J_{M'}^{\ov{Q'_1}}(\eta,Y,f)=|D^\fs(Y)|_F^{1/2} \int_{H_Y(F)\bs H(F)} f(\Ad(x^{-1})(Y)) \eta(\Nrd(x)) v_{M'}^{\ov{Q'_1}}(x) dx. $$
Then
\[\begin{split}
 &\int_{\fc'_\reg(F)} J_{M'}^{L'_1}(\eta,Y,f_{\ov{Q'_1}}^\eta) J_{M'}^{L'_2}(\eta,Y,\hat{f'}_{Q'_2}^\eta)dY \\
=&\int_{\fc'_\reg(F)} |D^\fs(Y)|_F^{1/2} \int_{T_{\fc'}(F)\bs H(F)} f(\Ad(x^{-1})(Y)) \eta(\Nrd(x)) v_{M'}^{\ov{Q'_1}}(x) J_{M'}^{L'_2}(\eta,Y,\hat{f'}_{Q'_2}^\eta) dxdY. 
\end{split}\]
Define a function $e_3: (T_{\fc'}(F)\bs H(F))\times\fc'_\reg(F)\ra\BC$ by
$$ e_3(x,Y):=\eta(\Nrd(x)) v_{M'}^{\ov{Q'_1}}(x) J_{M'}^{L'_2}(\eta,Y,\hat{f'}_{Q'_2}^\eta). $$
It is locally constant by the product form of Proposition \ref{propwoi1}.2). Using Lemma \ref{lemIII.5.1} to dominate $ v_{M'}^{\ov{Q'_1}}(x)$ and Corollary \ref{corIII.6.1} to dominate $J_{M'}^{L'_2}(\eta,Y,\hat{f'}_{Q'_2}^\eta)$, we check that $e_3$ verifies the hypotheses of Lemma \ref{lemV.7.1}, which implies our claim. 

Now (\ref{equV.10.(2)}) shows that the distribution $\hat{J}_M^G(\eta, X, \cdot)$ conincides with some element in $\msd^\fs$ for all $f\in\CC_c^\infty(\fs(F))$ with $\Supp(f)\subseteq r$. By glueing (Remark \ref{gluermk}), the distribution $\hat{J}_M^G(\eta, X, \cdot)\in\msd^\fs$. 

Finally, consider a general $X\in(\fm\cap\fs_\rs)(F)$. There exists $x\in M_H(F)$ such that $Y:=\Ad(x^{-1})(X)\in(\fl\cap\fs)(F)_\el$ for some $L\in\msl^{G,\omega}(M_0), L\subseteq M$. Then $\hat{J}_M^G(\eta, X, \cdot)=\eta(\Nrd(x))\hat{J}_M^G(\eta, Y, \cdot)$. Applying the descent formula (Proposition \ref{propwoi1}.5)), the product form of the elliptic case that we have just proved (applied to $Y\in(\fl\cap\fs)(F)_\el$) and parabolic induction (Corollary \ref{corV.8.1}), we deduce that the distribution $\hat{J}_M^G(\eta, Y, \cdot)\in\msd^\fs$. Thus the distribution $\hat{J}_M^G(\eta, X, \cdot)\in\msd^\fs$. 
\end{proof}

\subsection{The case of $(G',H')$}

We define $\mse^{\fs'}$ in the same way as the previous case. For $e\in\mse^{\fs'}$, thanks to Corollary \ref{corII.2}, it defines a distribution on $\fs'(F)$:
\begin{equation}\label{defV.6.2}
 \forall f'\in\CC_c^\infty(\fs'(F)), f'\mapsto\int_{\fs'(F)} f'(Y)e(Y)|D^{\fs'}(Y)|_F^{-1/2}dY. 
\end{equation}
Denote by $\msd^{\fs'}$ the space of distributions obtained in this way. For $d\in\msd^{\fs'}$, we shall always denote by $e_d\in\mse^{\fs'}$ its preimage under the isomorphism $\mse^{\fs'}\ra\msd^{\fs'}$ defined above. One may extend these definitions to the symmetric pair $(\wt{M'}, M', \Ad(\alpha))$, where $M'\in\msl^{H'}(M'_0)$. If $d\in\msd^{\fs'}$ is invariant with respect to the adjoint action of $H'(F)$, then so is $e_d\in\mse^{\fs'}$ and by the Weyl integration formula (Proposition \ref{wif2.2}), we have the equality
\begin{equation}\label{equV.6inv2}
 d(f')=\sum_{M'\in\msl^{H'}(M'_0)} |W_0^{M'}||W_0^{H'}|^{-1} \sum_{\fc'\in\mst_\el(\wt{\fm'}\cap\fs')} |W(M', \fc')|^{-1} \int_{\fc'_\reg(F)} J_{H'}^{H'}(Y, f') e_d(Y) dY 
\end{equation}
for all $f'\in\CC_c^\infty(\fs'(F))$, where $J_{H'}^{H'}(Y, f')$ is defined by (\ref{defwoi2}). 

Let $M'\in\msl^{H'}(M'_0)$ and $Y\in(\wt{\fm'}\cap\fs'_\rs)(F)$. Denote by $\hat{J}_{M'}^{H'}(Y, \cdot)$ the distribution on $\fs'(F)$ defined by
$$ \hat{J}_{M'}^{H'}(Y, f'):=J_{M'}^{H'}(Y, \hat{f'}) $$
for all $f'\in\CC_c^\infty(\fs'(F))$, where the right hand side is defined by (\ref{defwoi2}). One also has a similar definition for the symmetric pair $(\wt{M'}, M', \Ad(\alpha))$, where $M'\in\msl^{H'}(M'_0)$. 

\begin{prop}[Representability]\label{repr2}
Let $M'\in\msl^{H'}(M'_0)$ and $Y\in(\wt{\fm'}\cap\fs'_\rs)(F)$. Then the distribution $\hat{J}_{M'}^{H'}(Y, \cdot)\in\msd^{\fs'}$. 
\end{prop}

The rest of this section is devoted to the proof of Proposition \ref{repr2}. Although it is similar to the proof of Proposition \ref{repr1}, we shall sketch some steps for later use. 

Let $\fc'$ be a Cartan subspace of $\fs'$. Recall that $T_{\fc'}$ denotes the centraliser of $\fc'$ in $H'$. Suppose that $e_0: (T_{\fc'}(F)\bs H'(F))\times \fc'_\reg(F)\ra\BC$ is a function such that

(1) $e_0$ is locally constant; 

(2) for any open compact subset $r'$ of $\fs'(F)$, there exist constants $c>0$ and $N\in\BN$ such that for all $x\in T_{\fc'}(F)\bs H'(F)$ and $Y\in\fc'_\reg(F)$ satisfying $\Ad(x^{-1})(Y)\in r'$, one has the inequality
$$ e_0(x, Y)\leq c\sup\{1, -\log |D^{\fs'}(Y)|_F\}^N. $$
For $f'\in\CC_c^\infty(\fs'(F))$, we define
\begin{equation}\label{defV.7.2}
 d_0(f'):=\int_{\fc'_\reg(F)} |D^{\fs'}(Y)|_F^{1/2} \int_{T_{\fc'}(F)\bs H'(F)} f'(\Ad(x^{-1})(Y)) e_0(x, Y) dx dY. 
\end{equation}

\begin{lem}[cf. Lemma \ref{lemV.7.1}]\label{lemV.7.2}
Let $\fc'$ be a Cartan subspace of $\fs'$. Suppose that $e_0$ satisfies  the above hypotheses. Then the integral (\ref{defV.7.2}) is absolutely convergent. Moreover, the distribution $d_0\in\msd^{\fs'}$. 
\end{lem}

\begin{coro}[Parabolic induction]\label{corV.8.2}
Let $M'\in\msl^{H'}(M'_0)$, $P'\in\msp^{H'}(M')$ and $d\in\msd^{\wt{\fm'}\cap\fs'}$. Then the distribution on $\fs'(F)$ defined by $\forall f'\in\CC_c^\infty(\fs'(F)), f'\mapsto d(f'_{P'})$ belongs to $\msd^{\fs'}$, where $f'_{P'}$ is defined by (\ref{pardes2}). 
\end{coro}

\begin{proof}
We may apply the argument of Corollary \ref{corV.8.1} with the aid of Lemma \ref{jac2}. 
\end{proof}

Let $M'\in\msl^{H'}(M'_0)$ and $d\in\msd^{\wt{\fm'}\cap\fs'}$. Suppose that $d$ is invariant with respect to the adjoint action of $M'(F)$. We define a distribution $\Ind_{M'}^{H', w}(d)$ on $\fs'(F)$ by
\begin{equation}\label{wtind2}
 \Ind_{M'}^{H', w}(d)(f'):=\sum_{\{L'\in\msl^{H'}(M'_0): L'\subseteq M'\}} |W_0^{L'}||W_0^{M'}|^{-1} \sum_{\fc'\in\mst_\el(\wt{\fl'}\cap\fs')} |W(L', \fc')|^{-1} \int_{\fc'_\reg(F)} J_{M'}^{H'}(Y, f') e_d(Y) dY
\end{equation}
for all $f'\in\CC_c^\infty(\fs'(F))$, where $J_{M'}^{H'}(Y, f')$ is defined by (\ref{defwoi2}). In particular, if $M'=H'$ and $d\in\msd^{\fs'}$ is invariant with respect to the adjoint action of $H'(F)$, we have $\Ind_{H'}^{H', w}(d)=d$ by (\ref{equV.6inv2}). 

\begin{coro}\label{corV.9.2}
Let $M'\in\msl^{H'}(M'_0)$ and $d\in\msd^{\wt{\fm'}\cap\fs'}$. Suppose that $d$ is invariant with respect to the adjoint action of $M'(F)$. Then the integral (\ref{wtind2}) is absolutely convergent. Moreover, the distribution $\Ind_{M'}^{H', w}(d)\in\msd^{\fs'}$. 
\end{coro}

\begin{remark}
This corollary is unnecessary for the proof of Proposition \ref{repr2} but useful in Section \ref{sectinvwoi2}. 
\end{remark}

\begin{proof}[Proof of Corollary \ref{corV.9.2}]
We may apply the argument of Corollary \ref{corV.9.1} thanks to Lemmas \ref{lemIII.5.2} and \ref{lemV.7.2}. 
\end{proof}

\begin{lem}[cf. Lemma \ref{countlemforrepr1}]\label{countlemforrepr2}
Let $M'\in\msl^{H'}(M'_0)$ and $\fc'\subseteq\wt{\fm'}\cap\fs'$ be an $M'$-elliptic Cartan subspace. 

1) Let $M\in\msl^{H'}(M'_0), \fc\in\mst_\el(\wt{\fm}\cap\fs')$ and $x\in H'(F)$ be such that $\Ad(x)(\fc')=\fc $. Then there exist elements $m\in M(F)$ and $w\in\Norm_{H'(F)}(M'_0)$ such that $x=mw$. 

2) The cardinality of
$$ \{(M,\fc): M\in\msl^{H'}(M'_0), \fc\in\mst_\el(\wt{\fm}\cap\fs'), \fc\text{ is }H'(F)\text{-conjugate to }\fc'\} $$
is
$$ |W_0^{H'}||W_0^{M'}|^{-1}|W(M',\fc')||W(H',\fc')|^{-1}. $$
\end{lem}

\begin{proof}[Proof of Proposition \ref{repr2}]
We may apply the argument of Proposition \ref{repr1} with obvious modifications. One needs almost all results that we have prepared in this and previous sections, notably Howe's finiteness (Corollary \ref{corhowe2}) and the noninvariant trace formula (Theorem \ref{noninvltf2}). 
\end{proof}


\section{\textbf{Invariant weighted orbital integrals}}\label{secinvwoi}

\subsection{The case of $(G,H)$}\label{sectinvwoi1}

Let $M\in\msl^{G,\omega}(M_0)$ and $X\in(\fm\cap\fs_\rs)(F)$. We shall define a distribution $\hat{I}_M^G(\eta, X, \cdot)\in\msd^\fs$ which is $\eta(\Nrd(\cdot))$-invariant with respect to the adjoint action of $H(F)$ by induction on $\dim(G)$. Suppose that for all $L\in\msl^G(M), L\neq G$, we have defined a distribution $\hat{I}_M^L(\eta, X, \cdot)$ by an analogue of \eqref{invftofwoi1} where $G$ is replaced by $L$ such that the following properties are verified: $\hat{I}_M^L(\eta, X, \cdot)\in\msd^{\fl\cap\fs}$ and is $\eta(\Nrd(\cdot))$-invariant with respect to the adjoint action of $L_H(F)$. This is actually a product form of the induction hypothesis in lower dimensions. Denote by $\hat{I}_M^{L,G,w}(\eta, X, \cdot)$ its image under $\Ind_L^{G, w}$ (see (\ref{wtind1})). As in \cite[(1) in \S VI.1]{MR1344131}, for $f\in\CC_c^\infty(\fs(F))$, we set
\begin{equation}\label{invftofwoi1}
 \hat{I}_M^G(\eta, X, f):=\hat{J}_M^G(\eta, X, f)-\sum_{L\in\msl^G(M), L\neq G} \hat{I}_M^{L,G,w}(\eta, X, f). 
\end{equation}

\begin{prop}\label{propVI.1.1}
The distribution $\hat{I}_M^G(\eta, X, \cdot)\in\msd^\fs$ and is $\eta(\Nrd(\cdot))$-invariant with respect to the adjoint action of $H(F)$. 
\end{prop}

\begin{proof}
The first statement results from the representability of $\hat{J}_M^G(\eta, X, \cdot)$ (Proposition \ref{repr1}), the induction hypothesis and Corollary \ref{corV.9.1}. Now let us consider the second one. 

Let $f\in\CC_c^\infty(\fs(F))$ and $y\in H(F)$. By the $H(F)$-invariance of $\langle\cdot, \cdot\rangle$, we see that $(\Ad(y^{-1})(f))^\string^=\Ad(y^{-1})(\hat{f})$. Applying Proposition \ref{propwoi1}.6), we have
$$ \hat{J}_M^G(\eta, X, \Ad(y^{-1})(f))=J_M^G(\eta, X, \Ad(y^{-1})(\hat{f}))=\eta(\Nrd(y)) \sum_{Q\in\msf^G(M)} J_M^{M_Q}(\eta, X, (\hat{f})_{Q,y}^\eta). $$
For all $Q\in\msf^G(M)$, we show that
$$ (\hat{f})_{Q,y}^\eta=(f_{Q,y}^\eta)^\string^ $$
by the same argument of an analogous property of (\ref{pardes1}). Then
\begin{equation}\label{equVI.1.(2)}
 \hat{J}_M^G(\eta, X, \Ad(y^{-1})(f))=\eta(\Nrd(y)) \sum_{Q\in\msf^G(M)} \hat{J}_M^{M_Q}(\eta, X, f_{Q,y}^\eta). 
\end{equation}

Let $L\in\msl^G(M), L\neq G$. Applying Proposition \ref{propwoi1}.6) again to $J_L^G(\eta, Y, \Ad(y^{-1})(f))$ in the integrand of the definition (\ref{wtind1}) of $\Ind_L^{G, w}$, we obtain
\begin{equation}\label{equVI.1.(4)}
 \hat{I}_M^{L,G,w}(\eta, X, \Ad(y^{-1})(f))=\eta(\Nrd(y)) \sum_{Q\in\msf^G(L)} \hat{I}_M^{L,M_Q,w}(\eta, X, f_{Q,y}^\eta), 
\end{equation}
where $\hat{I}_M^{L,M_Q,w}(\eta, X, \cdot)\in\msd^{\fm_Q\cap\fs}$ denotes the image of $\hat{I}_M^L(\eta, X, \cdot)$ under $\Ind_L^{M_Q, w}$, which is defined by a product form of (\ref{wtind1}). 

From (\ref{invftofwoi1}), (\ref{equVI.1.(2)}) and (\ref{equVI.1.(4)}), we deduce that
\[\begin{split}
 \hat{I}_M^G(\eta, X, \Ad(y^{-1})(f))=&\hat{J}_M^G(\eta, X, \Ad(y^{-1})(f))-\sum_{L\in\msl^G(M), L\neq G} \hat{I}_M^{L,G,w}(\eta, X, \Ad(y^{-1})(f)) \\
 =&\eta(\Nrd(y)) \sum_{Q\in\msf^G(M)} \left(\hat{J}_M^{M_Q}(\eta, X, f_{Q,y}^\eta)-\sum_{L\in\msl^{M_Q}(M), L\neq G}\hat{I}_M^{L,M_Q,w}(\eta, X, f_{Q,y}^\eta)\right). 
\end{split}\]
Consider $Q\neq G$ first. By the induction hypothesis, $\hat{I}_M^{M_Q}(\eta, X, \cdot)\in\msd^{\fm_Q\cap\fs}$ is $\eta(\Nrd(\cdot))$-invariant with respect to the adjoint action of $M_{Q_H}(F)$, so $\hat{I}_M^{M_Q,M_Q,w}(\eta, X, \cdot)=\hat{I}_M^{M_Q}(\eta, X, \cdot)$. By the definition of $\hat{I}_M^{M_Q}(\eta, X, \cdot)$ (a product form of (\ref{invftofwoi1})), the term in brackets is zero. Thus it remains the term for $Q=G$. Note that $f_{G,y}^\eta=f_G^\eta$. By Proposition \ref{propwoi1}.4) applied to $Q=G$, we see that the term in brackets is exactly $\hat{I}_M^G(\eta, X, f)$ defined by (\ref{invftofwoi1}). Therefore, we show that
$$ \hat{I}_M^G(\eta, X, \Ad(y^{-1})(f))=\eta(\Nrd(y))\hat{I}_M^G(\eta, X, f), $$
which is the second statement. 
\end{proof}

Let $M\in\msl^{G,\omega}(M_0)$ and $X\in(\fm\cap\fs_\rs)(F)$. Denote by $\hat{i}_M^G(\eta, X, \cdot)$ (resp. $\hat{j}_M^G(\eta, X, \cdot)$) the element of $\mse^\fs$ associated to $\hat{I}_M^G(\eta, X, \cdot)$ (resp. $\hat{J}_M^G(\eta, X, \cdot)$) $\in\msd^\fs$ by (\ref{defV.6.1}). That is to say, for all $f\in\CC_c^\infty(\fs(F))$, 
$$ \hat{I}_M^G(\eta, X, f)=\int_{\fs(F)} f(Y) \hat{i}_M^G(\eta, X, Y) |D^\fs(Y)|_F^{-1/2} dY $$
and
$$ \hat{J}_M^G(\eta, X, f)=\int_{\fs(F)} f(Y) \hat{j}_M^G(\eta, X, Y) |D^\fs(Y)|_F^{-1/2} dY. $$
We also have a similar definition for the symmetric pair $(M, M_H, \Ad(\epsilon))$, where $M\in\msl^{G,\omega}(M_0)$. 

\begin{lem}\label{lemVI.2.1}
Let $M\in\msl^{G,\omega}(M_0)$ and $X\in(\fm\cap\fs_\rs)(F)$. Let $L\in\msl^{G,\omega}(M_0)$ and $Y\in(\fl\cap\fs_\rs)(F)_\el$. Then $\hat{i}_M^G(\eta, X, Y)=\hat{j}_M^G(\eta, X, Y)$. 
\end{lem}

\begin{proof}
Let $L_2\in\msl^{G}(M), L_2\neq G$, $L_1\in\msl^{G, \omega}(M_0), L_1\subseteq L_2$ and $\fc\in\mst_\el(\fl_1\cap\fs)$. 
We define a distribution $d_{L_1,L_2,\fc}^{G, w}(\eta, X, \cdot)$ on $\fs(F)$ by
\begin{equation}\label{equlemVI.2}
 d_{L_1,L_2,\fc}^{G, w}(\eta, X, f):=\int_{\fc_\reg(F)} J_{L_2}^G(\eta, Z, f) \hat{i}_M^{L_2}(\eta, X, Z) dZ 
\end{equation}
for all $f\in\CC_c^\infty(\fs(F))$. By a product form of Proposition \ref{propVI.1.1}, the distribution $\hat{I}_M^{L_2}(\eta, X, \cdot)\in\msd^{\fl_2\cap\fs}$ and is $\eta(\Nrd(\cdot))$-invariant with respect to the adjoint action of $(L_2\cap H)(F)$. We see from the proof of Corollary \ref{corV.9.1} that the integral (\ref{equlemVI.2}) is absolutely convergent and the distribution $d_{L_1,L_2,\fc}^{G, w}(\eta, X, \cdot)\in\msd^\fs$. We denote by $e_{L_1,L_2,\fc}^{G, w}(\eta, X, \cdot)$ its associated element in $\mse^\fs$ by (\ref{defV.6.1}). From the definitions (\ref{invftofwoi1}) and (\ref{wtind1}), we have
\[\begin{split} 
\hat{i}_M^G(\eta, X, \cdot)=&\hat{j}_M^G(\eta, X, \cdot)-\sum_{L_2\in\msl^G(M), L_2\neq G} \sum_{\{L_1\in\msl^{G,\omega}(M_0): L_1\subseteq L_2\}} |W_0^{L_{1,n}}| |W_0^{L_{2,n}}|^{-1}\sum_{\fc\in\mst_\el(\fl_1\cap\fs)} \\ 
&|W(L_{1}\cap H,\fc)|^{-1} e_{L_1,L_2,\fc}^{G, w}(\eta, X, \cdot). 
\end{split}\]
To prove the lemma, it suffices to fix such a triple $(L_2, L_1, \fc)$ and prove that $e_{L_1,L_2,\fc}^{G, w}(\eta, X, Y)=0$. But (\ref{equpfcorV.9}) and (\ref{defe'}) in the proofs of Corollary \ref{corV.9.1} and Lemma \ref{lemV.7.1} respectively allow us to calculate $e_{L_1,L_2,\fc}^{G, w}(\eta, X, Y)$; explicitly, we have
\begin{equation}\label{equpfVI.2}
 e_{L_1,L_2,\fc}^{G, w}(\eta, X, Y)=\sum_{\{x\in T_\fc(F)\bs H(F):\Ad(x)(Y)\in \fc(F)\}} \eta(\Nrd(x)) v_{L_2}^G(x) \hat{i}_M^{L_2}(\eta, X, \Ad(x)(Y)). 
\end{equation}
Let $x\in T_\fc(F)\bs H(F)$ such that $\Ad(x)(Y)\in \fc(F)$. As $Y\in(\fl\cap\fs_\rs)(F)_\el$ and $\fc\in\mst_\el(\fl_1\cap\fs)$, from the proof of Lemma \ref{countlemforwif2.1}.1), there exist elements $l_1\in (L_1\cap H)(F)$ and $w\in\Norm_{H(F)}(M_0)$ such that $x=l_1 w$. Since any element in $W_0^H$ admits a representative in $K_H$, we can suppose that $w\in K_H$. Then $v_{L_2}^G(x)=v_{L_2}^G(1)$ since $L_1\subseteq L_2$. But $v_{L_2}^G(1)=0$ for $L_2\neq G$. Thus $e_{L_1,L_2,\fc}^{G, w}(\eta, X, Y)=0$ by (\ref{equpfVI.2}). 
\end{proof}

\begin{lem}\label{lemVI.3.1}
Let $M\in\msl^{G, \omega}(M_0)$. 

1) The function $(X, Y)\mapsto\hat{i}_M^G(\eta, X, Y)$ is locally constant on $(\fm\cap\fs_\rs)(F)\times\fs_\rs(F)$. 

2) If $w\in\Norm_{H(F)}(M_0)$, $x\in M_H(F)$ and $y\in H(F)$, we have the equality
$$ \hat{i}_{\Ad(w)(M)}^G(\eta, \Ad(wx)(X), \Ad(y)(Y))=\eta(\Nrd(wxy))\hat{i}_M^G(\eta, X, Y) $$
for all $(X, Y)\in(\fm\cap\fs_\rs)(F)\times\fs_\rs(F)$. 

3) If $\lambda\in F^\times$, we have the equality
$$ \hat{i}_M^G(\eta, \lambda X, Y)=\hat{i}_M^G(\eta, X, \lambda Y) $$
for all $(X, Y)\in(\fm\cap\fs_\rs)(F)\times\fs_\rs(F)$. 

4) Let $r_M\subseteq(\fm\cap\fs)(F)$ and $r\subseteq\fs(F)$ be two compact subsets. Then there exist constants $c>0$ and $N\in\BN$ such that
$$ |\hat{i}_M^G(\eta, X, Y)|\leq c\sup\{1, -\log|D^\fs(X)|_F\}^N \sup\{1, -\log|D^\fs(Y)|_F\}^N $$
for all $X\in r_M\cap\fs_\rs$ and $Y\in r\cap\fs_\rs$. 

5) Let $f\in\CC_c^\infty(\fs(F))$ and $r_M\subseteq(\fm\cap\fs)(F)$ be a compact subset. Then there exist constants $c>0$ and $N\in\BN$ such that
$$ |\hat{I}_M^G(\eta, X, f)|\leq c\sup\{1, -\log|D^\fs(X)|_F\}^N $$
for all $X\in r_M\cap\fs_\rs$. 
\end{lem}

\begin{proof}
Let $r_M\subseteq(\fm\cap\fs)(F)$ and $r\subseteq\fs(F)$ be two open compact subgroups. Set $r^\ast:=\{Y\in\fs(F): \forall Z\in r, \Psi(\langle Y,Z \rangle)=1\}$, which is an open compact subgroup of $\fs(F)$. Notice that if $f\in\CC_c^\infty(\fs(F))$ satisfies $\Supp(f)\subseteq r$, then $\hat{f}\in\CC_c^\infty(\fs(F)/r^\ast)$. Applying Howe's finiteness (Corollary \ref{corhowe1}) to $r^\ast$ and $r_M\cap\fs_\rs$, we know that there exists a finite subset $\{X_i: i\in I\}\subseteq r_M\cap\fs_\rs$ and a finite subset $\{f_i: i\in I\}\subseteq \CC_c^\infty(\fs(F)/r^\ast)$ such that for all $X\in r_M\cap\fs_\rs$ and all $f\in\CC_c^\infty(\fs(F))$ with $\Supp(f)\subseteq r$, we have
$$ J_M^G(\eta,X,\hat{f})=\sum_{i\in I} J_M^G(\eta,X_i,\hat{f}) J_M^G(\eta,X,f_i). $$
We deduce that
\begin{equation}\label{equpfVI.3}
 \hat{j}_M^G(\eta, X, Y)=\sum_{i\in I} \hat{j}_M^G(\eta, X_i, Y) J_M^G(\eta,X,f_i) 
\end{equation}
for all $X\in r_M\cap\fs_\rs$ and $Y\in r\cap\fs_\rs$. 

1) The local constancy of $(X,Y)\mapsto\hat{j}_M^G(\eta,X,Y)$ on $(\fm\cap\fs_\rs)(F)\times\fs_\rs(F)$ results from (\ref{equpfVI.3}), Proposition \ref{propwoi1}.2) and $\hat{j}_M^G(\eta, X_i, \cdot)\in\mse^\fs$ for $i\in I$. For $L\in\msl^{G,\omega}(M_0)$, we deduce from Lemma \ref{lemVI.2.1} the local constancy of $(X, Y)\mapsto\hat{i}_M^G(\eta, X, Y)$ on $(\fm\cap\fs_\rs)(F)\times(\fl\cap\fs_\rs)(F)_\el$. Let $(X,Y)\in(\fm\cap\fs_\rs)(F)\times\fs_\rs(F)$. Choose $L\in\msl^{G,\omega}(M_0)$ and $Y'\in(\fl\cap\fs_\rs)(F)_\el$ such that $Y'$ is $H(F)$-conjugate to $Y$. Fix a neighbourhood $V_1\times V_2$ of $(X,Y')$ in $(\fm\cap\fs_\rs)(F)\times(\fl\cap\fs_\rs)(F)_\el$ such that $(X, Y)\mapsto\kappa(Y)\hat{i}_M^G(\eta, X, Y)$ is constant on $V_1\times V_2$. Thanks to the $\eta(\Nrd(\cdot))$-invariance of $\hat{i}_M^G(\eta, X, \cdot)$ with respect to the adjoint action of $H(F)$ (Proposition \ref{propVI.1.1}), we know that $(X, Y)\mapsto\kappa(Y)\hat{i}_M^G(\eta, X, Y)$ is constant on $V_1\times \Ad(H(F))(V_2)$ which is a neighbourhood of $(X,Y)$ in $(\fm\cap\fs_\rs)(F)\times\fs_\rs(F)$. Since $\kappa(\cdot)$ is locally constant on $\fs_\rs(F)$, we show that $(X, Y)\mapsto\hat{i}_M^G(\eta, X, Y)$ is constant on a neighbourhood of $(X,Y)$ in $(\fm\cap\fs_\rs)(F)\times\fs_\rs(F)$. 

2) The effect of $\Ad(y)$ comes from Propostion \ref{propVI.1.1}. Then when considering the effects of $\Ad(w)$ and $\Ad(x)$, up to $H(F)$-conjugation, we may and shall suppose that $Y\in(\fl\cap\fs_\rs)(F)_\el$ for some $L\in\msl^{G,\omega}(M_0)$. That is to say, it suffices to prove the equality
$$ \hat{i}_{\Ad(w)(M)}^G(\eta, \Ad(wx)(X), Y)=\eta(\Nrd(wx))\hat{i}_M^G(\eta, X, Y) $$
for all $(X, Y)\in(\fm\cap\fs_\rs)(F)\times(\fl\cap\fs_\rs)(F)_\el$. By Lemma \ref{lemVI.2.1}, we may replace $\hat{i}_{\Ad(w)(M)}^G$ and $\hat{i}_M^G$ by $\hat{j}_{\Ad(w)(M)}^G$ and $\hat{j}_M^G$ respectively in the equality to be proved. Now the equality results from Proposition \ref{propwoi1}.3). 

3) Let $\lambda\in F^\times, X\in(\fm\cap\fs_\reg)(F)$ and $f\in\CC_c^\infty(\fs(F))$. From (\ref{weyldisc}), we have
$$ |D^\fs(\lambda X)|_F=|\lambda|_F^{(\dim(\fg)-\rank(\fg))/2}|D^\fs(X)|_F, $$
where $\dim(\fg)$ and $\rank(\fg)$ denote the dimension and rank (over an algebraic closure of $F$) of $\fg$ respectively. Then we have
\[\begin{split}
&\int_{\fs(F)} f(Y) \hat{j}_M^G(\eta, \lambda X, Y) |D^\fs(Y)|_F^{-1/2} dY=\hat{J}_M^G(\eta, \lambda X, f)=J_M^G(\eta, \lambda X, \hat{f}) \\ 
=&|D^\fs(\lambda X)|_F^{1/2} \int_{H_{\lambda X}(F)\bs H(F)} \hat{f}(\Ad(x^{-1})(\lambda X)) \eta(\Nrd(x)) v_M^Q(x) dx \\
=&|\lambda|_F^{(\dim(\fg)-\rank(\fg))/4} J_M^G(\eta, X, \hat{f}(\lambda\cdot)). 
\end{split}\]
But
\[\begin{split}
 \hat{f}(\lambda\cdot)=&c_\Psi(\fs(F)) \int_{\fs(F)} f(Z) \Psi(\langle \lambda\cdot,Z \rangle) dZ=c_\Psi(\fs(F)) \int_{\fs(F)} f(Z) \Psi(\langle \cdot, \lambda Z \rangle) dZ \\
=&|\lambda|_F^{-\dim(\fs)}c_\Psi(\fs(F)) \int_{\fs(F)} f(\lambda^{-1}Z) \Psi(\langle \cdot, Z \rangle) dZ \\
=&|\lambda|_F^{-\dim(\fg)/2}(f(\lambda^{-1}\cdot))^\string^. 
\end{split}\]
Thus we have
\[\begin{split}
J_M^G(\eta, X, \hat{f}(\lambda\cdot))=&|\lambda|_F^{-\dim(\fg)/2}\hat{J}_M^G(\eta, X, f(\lambda^{-1}\cdot))=|\lambda|_F^{-\dim(\fg)/2}\int_{\fs(F)} f(\lambda^{-1}Y) \hat{j}_M^G(\eta, X, Y) |D^\fs(Y)|_F^{-1/2} dY \\
=&\int_{\fs(F)} f(Y) \hat{j}_M^G(\eta, X, \lambda Y) |D^\fs(\lambda Y)|_F^{-1/2} dY \\
=&|\lambda|_F^{(\rank(\fg)-\dim(\fg))/4} \int_{\fs(F)} f(Y) \hat{j}_M^G(\eta, X, \lambda Y) |D^\fs(Y)|_F^{-1/2} dY. 
\end{split}\]
Therefore, we deduce the equality
$$ \hat{j}_M^G(\eta, \lambda X, Y)=\hat{j}_M^G(\eta, X, \lambda Y) $$
for all $\lambda\in F^\times$ and all $(X, Y)\in(\fm\cap\fs_\rs)(F)\times\fs_\rs(F)$. We obtain a similar equality for $\hat{i}_M^G$ thanks to Lemma \ref{lemVI.2.1} and the $\eta(\Nrd(\cdot))$-invariance of $\hat{i}_M^G(\eta, X, \cdot)$ with respect to the adjoint action of $H(F)$ (Proposition \ref{propVI.1.1}). 

4) A similar bound for $(X,Y)\mapsto\hat{j}_M^G(\eta,X,Y)$ on $(r_M\cap\fs_\rs)\times(r\cap\fs_\rs)$ results from (\ref{equpfVI.3}), Corollary \ref{corIII.6.1} (together with Proposition \ref{propwoi1}.4)) and $\hat{j}_M^G(\eta, X_i, \cdot)\in\mse^\fs$ for $i\in I$. For $L\in\msl^{G,\omega}(M_0)$, we deduce from Lemma \ref{lemVI.2.1} a similar bound of $(X, Y)\mapsto\hat{i}_M^G(\eta, X, Y)$ on $(r_M\cap\fs_\rs)\times(r\cap(\fl\cap\fs_\rs)(F)_\el)$. Let $(X,Y)\in(r_M\cap\fs_\rs)\times(r\cap\fs_\rs)$. Thanks to the $\eta(\Nrd(\cdot))$-invariance of $\hat{i}_M^G(\eta, X, \cdot)$ with respect to the adjoint action of $H(F)$ (Proposition \ref{propVI.1.1}), if we replace $Y$ by $\Ad(y)(Y)$, where $y\in H(F)$, the two sides in the inequality to be proved remain unchanged. Since any Cartan subspace in $\fs$ is $H(F)$-conjugate to an element in $\mst_\el(\fl\cap\fs)$ for some $L\in\msl^{G,\omega}(M_0)$, with the help of Lemma \ref{lem28}, it suffices to fix $L\in\msl^{G,\omega}(M_0), \fc\in\mst_\el(\fl\cap\fs)$ and $r_\fc\subseteq\fc(F)$ a compact subset, prove a similar bound for $(X,Y)\in(r_M\cap\fs_\rs)\times(r_\fc\cap\fc_\reg)$, and then obtain a uniform bound for $(X,Y)\in(r_M\cap\fs_\rs)\times(r\cap\fs_\rs)$ by the finiteness of $\mst_\el(\fl\cap\fs)$. But this is what we have established. 

5) It is a consequence of 4) applied to $r:=\Supp(f)$ and Corollary \ref{corII.1}. 
\end{proof}

For $M\in\msl^{G, \omega}(M_0)$ and $X\in(\fm\cap\fs_\rs)(F)$, we define a distribution $I_M^G(\eta, X, \cdot)$ on $\fs(F)$ by
\begin{equation}\label{invwoi1}
 I_M^G(\eta, X, \hat{f}):=\hat{I}_M^G(\eta, X, f) 
\end{equation}
for all $f\in\CC_c^\infty(\fs(F))$. 

\begin{remark}
For $M=G$, it is evident that $I_G^G(\eta, X, f)=J_G^G(\eta, X, f)$ for all $X\in\fs_\rs(F)$ and $f\in\CC_c^\infty(\fs(F))$. 
\end{remark}

One may easily extend the definitions (\ref{invftofwoi1}) and (\ref{invwoi1}) to the symmetric pair $(L, L_H, \Ad(\epsilon))$, where $L\in\msl^{G,\omega}(M_0)$, since it appears as the product of some copies of the form $(G, H, \Ad(\epsilon))$ in lower dimensions. 

\begin{lem}\label{lemVI.5.1}
Let $M\in\msl^{G, \omega}(M_0)$ and $X\in(\fm\cap\fs_\rs)(F)$. The distribution $I_M^G(\eta, X, \cdot)$ on $\fs(F)$ is independent of the choice of the $H(F)$-invariant non-degenerate symmetric bilinear form $\langle\cdot, \cdot\rangle$ on $\fs(F)$ or the continuous nontrivial unitary character $\Psi$ of $F$. 
\end{lem}

\begin{proof}
Suppose that $\langle\cdot,\cdot\rangle'$ is another bilinear form and that $\Psi'$ is another character. Denote by $f\mapsto\tilde{f}$ the associated Fourier transform and by $\tilde{I}_M^G(\eta, X, \cdot)$ (resp. $\tilde{J}_M^G(\eta, X, \cdot)$) the associated analogue of $\hat{I}_M^G(\eta, X, \cdot)$ (resp. $\hat{J}_M^G(\eta, X, \cdot)$). For $f\in\CC_c^\infty(\fs(F))$, define $f_{-}\in\CC_c^\infty(\fs(F))$ by $f_{-}(Y):=f(-Y)$ for all $Y\in\fs(F)$. Since $\tilde{\tilde{f}}=\hat{\hat{f}}=f_{-}$ for all $f\in\CC_c^\infty(\fs(F))$, it suffices to prove the equality
$$ \tilde{I}_M^G(\eta, X, \tilde{f})=\hat{I}_M^G(\eta, X, \hat{f}) $$
for all $f\in\CC_c^\infty(\fs(F))$. 

Let $\tau'$ be the linear automorphism of $\fs(F)$ such that
$$ \forall Y,Z\in\fs(F), \langle Y,Z\rangle'=\langle \tau'(Y),Z\rangle. $$
Let $a\in F^\times$ such that $\Psi'(\cdot)=\Psi(a\cdot)$. Set $\tau:=a\tau'$. Then
$$ \tilde{f}(\cdot)=\frac{c_{\Psi'}(\fs(F))}{c_{\Psi}(\fs(F))}\hat{f}(\tau(\cdot)) $$
for all $f\in\CC_c^\infty(\fs(F))$. One may check that $\tau$ is an $H(F)$-equivariant linear automorphism of $\fs(F)$ thanks to $H(F)$-invariance of two bilinear forms. One also deduces that
$$ \forall Y,Z\in\fs(F), \langle\tau(Y), Z\rangle=\langle Y,\tau(Z)\rangle $$
from the symmetry of two bilinear forms. Now for all $f\in\CC_c^\infty(\fs(F))$ and all $Y\in\fs(F)$, we have
\[\begin{split}
f(-Y)&=\tilde{\tilde{f}}(Y)=c_{\Psi'}(\fs(F)) \int_{\fs(F)} \tilde{f}(Z) \Psi'(\langle Y,Z \rangle') dZ =c_{\Psi'}(\fs(F)) \int_{\fs(F)} \frac{c_{\Psi'}(\fs(F))}{c_{\Psi}(\fs(F))}\hat{f}(\tau(Z)) \Psi(\langle \tau(Y),Z \rangle) dZ \\
&=\frac{c_{\Psi'}(\fs(F))^2}{c_{\Psi}(\fs(F))} \int_{\fs(F)} \hat{f}(\tau(Z)) \Psi(\langle Y,\tau(Z) \rangle) dZ=\frac{c_{\Psi'}(\fs(F))^2}{c_{\Psi}(\fs(F))|\det_{\fs(F)}(\tau)|_F} \int_{\fs(F)} \hat{f}(Z') \Psi(\langle Y,Z' \rangle) dZ' \\
&=\frac{c_{\Psi'}(\fs(F))^2}{c_{\Psi}(\fs(F))^2|\det_{\fs(F)}(\tau)|_F} \hat{\hat{f}}(Y)=\frac{c_{\Psi'}(\fs(F))^2}{c_{\Psi}(\fs(F))^2|\det_{\fs(F)}(\tau)|_F} f(-Y).
\end{split}\]
Therefore, we obtain $\frac{c_{\Psi'}(\fs(F))}{c_{\Psi}(\fs(F))}=|\det_{\fs(F)}(\tau)|_F^{1/2}$. Then for all $f\in\CC_c^\infty(\fs(F))$, we have
$$ \tilde{f}(\cdot)=|\det\nolimits_{\fs(F)}(\tau)|_F^{1/2}\hat{f}(\tau(\cdot)). $$

Denote by $\tilde{i}_M^G(\eta, X, \cdot)$ (resp. $\tilde{j}_M^G(\eta, X, \cdot)$) the element of $\mse^\fs$ associated to $\tilde{I}_M^G(\eta, X, \cdot)$ (resp. $\tilde{J}_M^G(\eta, X, \cdot)$) $\in\msd^\fs$ by (\ref{defV.6.1}). For $f\in\CC_c^\infty(\fs(F))$, we have
\[\begin{split}
 \tilde{I}_M^G(\eta, X, \tilde{f})&=\int_{\fs(F)} \tilde{f}(Y) \tilde{i}_M^G(\eta, X, Y) |D^\fs(Y)|_F^{-1/2} dY \\
&=|\det\nolimits_{\fs(F)}(\tau)|_F^{1/2}\int_{\fs(F)} \hat{f}(\tau(Y)) \tilde{i}_M^G(\eta, X, Y) |D^\fs(Y)|_F^{-1/2} dY \\
&=|\det\nolimits_{\fs(F)}(\tau)|_F^{-1/2}\int_{\fs(F)} \hat{f}(Y) \tilde{i}_M^G(\eta, X, \tau^{-1}(Y)) |D^\fs(\tau^{-1}(Y))|_F^{-1/2} dY. 
\end{split}\]
We reduce ourselves to proving the equality
\begin{equation}\label{equVI.5(1)}
|\det\nolimits_{\fs(F)}(\tau)|_F^{-1/2} \tilde{i}_M^G(\eta, X, \tau^{-1}(Y)) |D^\fs(\tau^{-1}(Y))|_F^{-1/2}=\hat{i}_M^G(\eta, X, Y) |D^\fs(Y)|_F^{-1/2}
\end{equation}
for all $Y\in\fs_\rs(F)$. But we have the equality
$$ \tilde{J}_M^G(\eta, X, \tilde{f})=\hat{J}_M^G(\eta, X, \hat{f}) $$
since both sides equal $J_M^G(\eta, X, f_{-}(\cdot))$, which is defined by (\ref{defwoi1}). The same computation as above shows that the equality (\ref{equVI.5(1)}) is true when one replaces $\tilde{i}_M^G$ and $\hat{i}_M^G$ with $\tilde{j}_M^G$ and $\hat{j}_M^G$ respectively. Recall that $\tau$ is $H(F)$-equivariant, so $H_Y=H_{\tau^{-1}(Y)}$ for $Y\in\fs_\rs(F)$. As a consequence, for $L\in\msl^{G,\omega}(M_0)$, $Y\in(\fl\cap\fs_\rs)(F)_\el$ if and only if $\tau^{-1}(Y)\in(\fl\cap\fs_\rs)(F)_\el$. One may conclude by Lemma \ref{lemVI.2.1} together with the $\eta(\Nrd(\cdot))$-invariance of $\hat{i}_M^G(\eta, X, \cdot)$ with respect to the adjoint action of $H(F)$ (Proposition \ref{propVI.1.1}). 
\end{proof}

\subsection{The case of $(G',H')$}\label{sectinvwoi2}

Let $M'\in\msl^{H'}(M'_0)$ and $Y\in(\wt{\fm'}\cap\fs'_\rs)(F)$. We shall define a distribution $\hat{I}_{M'}^{H'}(Y, \cdot)\in\msd^{\fs'}$ which is invariant with respect to the adjoint action of $H'(F)$ by induction on $\dim(H')$. Suppose that we have defined a distribution $\hat{I}_{M'}^{L'}(Y, \cdot)\in\msd^{\wt{\fl'}\cap\fs'}$ which is invariant with respect to the adjoint action of $L'(F)$ for all $L'\in\msl^{H'}(M'), L'\neq H'$. This is actually a product form of the induction hypothesis in lower dimensions. Denote by $\hat{I}_{M'}^{L',H',w}(Y, \cdot)$ its image under $\Ind_{L'}^{H', w}$ (see (\ref{wtind2})). For $f'\in\CC_c^\infty(\fs'(F))$, we set
\begin{equation}\label{invftofwoi2}
 \hat{I}_{M'}^{H'}(Y, f'):=\hat{J}_{M'}^{H'}(Y, f')-\sum_{L'\in\msl^{H'}(M'), L'\neq H'} \hat{I}_{M'}^{L',H',w}(Y, f'). 
\end{equation}

\begin{prop}\label{propVI.1.2}
The distribution $\hat{I}_{M'}^{H'}(Y, \cdot)\in\msd^{\fs'}$ and is invariant with respect to the adjoint action of $H'(F)$. 
\end{prop}

\begin{proof}
We may apply the argument of Proposition \ref{propVI.1.1} thanks to the representability of $\hat{J}_{M'}^{H'}(Y, \cdot)$ (Proposition \ref{repr2}), Corollary \ref{corV.9.2} and Proposition \ref{propwoi2}.6). 
\end{proof}

Let $M'\in\msl^{H'}(M'_0)$ and $Y\in(\wt{\fm'}\cap\fs'_\rs)(F)$. Denote by $\hat{i}_{M'}^{H'}(Y, \cdot)$ (resp. $\hat{j}_{M'}^{H'}(Y, \cdot)$) the element of $\mse^{\fs'}$ associated to $\hat{I}_{M'}^{H'}(Y, \cdot)$ (resp. $\hat{J}_{M'}^{H'}(Y, \cdot)$) $\in\msd^{\fs'}$ by (\ref{defV.6.2}). That is to say, for all $f'\in\CC_c^\infty(\fs'(F))$, 
$$ \hat{I}_{M'}^{H'}(Y, f')=\int_{\fs'(F)} f'(X) \hat{i}_{M'}^{H'}(Y, X) |D^{\fs'}(X)|_F^{-1/2} dX $$
and
$$ \hat{J}_{M'}^{H'}(Y, f')=\int_{\fs'(F)} f'(X) \hat{j}_{M'}^{H'}(Y, X) |D^{\fs'}(X)|_F^{-1/2} dX. $$
One has a similar definition for the symmetric pair $(\wt{M'}, M', \Ad(\alpha))$, where $M'\in\msl^{H'}(M'_0)$. 

\begin{lem}\label{lemVI.2.2}
Let $M'\in\msl^{H'}(M'_0)$ and $Y\in(\wt{\fm'}\cap\fs'_\rs)(F)$. Let $L'\in\msl^{H'}(M'_0)$ and $X\in(\wt{\fl'}\cap\fs'_\rs)(F)_\el$. Then $\hat{i}_{M'}^{H'}(Y, X)=\hat{j}_{M'}^{H'}(Y, X)$. 
\end{lem}

\begin{proof}
We may apply the argument of Lemma \ref{lemVI.2.1} by using Proposition \ref{propVI.1.2} and consulting the proofs of Corollary \ref{corV.9.2}, Lemmas \ref{lemV.7.2} and \ref{countlemforwif2.2}.1). 
\end{proof}

\begin{lem}\label{lemVI.3.2}
Let $M'\in\msl^{H'}(M'_0)$. 

1) The function $(Y, X)\mapsto\hat{i}_{M'}^{H'}(Y, X)$ is locally constant on $(\wt{\fm'}\cap\fs'_\rs)(F)\times\fs'_\rs(F)$. 

2) If $w\in\Norm_{H'(F)}(M'_0)$, $x\in M'(F)$ and $y\in H'(F)$, we have the equality
$$ \hat{i}_{\Ad(w)(M')}^{H'}(\Ad(wx)(Y), \Ad(y)(X))=\hat{i}_{M'}^{H'}(Y, X) $$
for all $(Y, X)\in(\wt{\fm'}\cap\fs'_\rs)(F)\times\fs'_\rs(F)$. 

3) If $\lambda\in F^\times$, we have the equality
$$ \hat{i}_{M'}^{H'}(\lambda Y, X)=\hat{i}_{M'}^{H'}(Y, \lambda X) $$
for all $(Y, X)\in(\wt{\fm'}\cap\fs'_\rs)(F)\times\fs'_\rs(F)$. 

4) Let $r'_{M'}\subseteq(\wt{\fm'}\cap\fs')(F)$ and $r'\subseteq\fs'(F)$ be two compact subsets. Then there exist constants $c>0$ and $N\in\BN$ such that
$$ |\hat{i}_{M'}^{H'}(Y, X)|\leq c\sup\{1, -\log|D^{\fs'}(Y)|_F\}^N \sup\{1, -\log|D^{\fs'}(X)|_F\}^N $$
for all $Y\in r'_{M'}\cap\fs'_\rs$ and $X\in r'\cap\fs'_\rs$. 

5) Let $f'\in\CC_c^\infty(\fs'(F))$ and $r'_{M'}\subseteq(\wt{\fm'}\cap\fs')(F)$ be a compact subset. Then there exist constants $c>0$ and $N\in\BN$ such that
$$ |\hat{I}_{M'}^{H'}(Y, f')|\leq c\sup\{1, -\log|D^{\fs'}(Y)|_F\}^N $$
for all $Y\in r'_{M'}\cap\fs'_\rs$. 
\end{lem}

\begin{proof}
It is almost the same as the proof of Lemma \ref{lemVI.3.1}, except that one needs to use Howe's finiteness (Corollary \ref{corhowe2}), Proposition \ref{propVI.1.2} and Lemma \ref{lemVI.2.2}. We also need Proposition \ref{propwoi2}.2) for 1), Proposition \ref{propwoi2}.3) for 2), Corollary \ref{corIII.6.2} for 4) and Corollary \ref{corII.2} for 5). 
\end{proof}

For $M'\in\msl^{H'}(M'_0)$ and $Y\in(\wt{\fm'}\cap\fs'_\rs)(F)$, we define a distribution $I_{M'}^{H'}(Y, \cdot)$ on $\fs'(F)$ by
\begin{equation}\label{invwoi2}
 I_{M'}^{H'}(Y, \hat{f'}):=\hat{I}_{M'}^{H'}(Y, f') 
\end{equation}
for all $f'\in\CC_c^\infty(\fs'(F))$. One may easily extend the definitions (\ref{invftofwoi2}) and (\ref{invwoi2}) to the symmetric pair $(\wt{L'}, L', \Ad(\alpha))$, where $L'\in\msl^{H'}(M'_0)$. 

\begin{remark}
For $M'=H'$, it is evident that $I_{H'}^{H'}(Y, f')=J_{H'}^{H'}(Y, f')$ for all $Y\in\fs'_\rs(F)$ and $f'\in\CC_c^\infty(\fs'(F))$. 
\end{remark}

\begin{lem}\label{lemVI.5.2}
Let $M'\in\msl^{H'}(M'_0)$ and $Y\in(\wt{\fm'}\cap\fs'_\rs)(F)$. The distribution $I_{M'}^{H'}(Y, \cdot)$ on $\fs'(F)$ is independent of the choice of the $H'(F)$-invariant non-degenerate symmetric bilinear form $\langle\cdot, \cdot\rangle$ on $\fs'(F)$ or the continuous nontrivial unitary character $\Psi$ of $F$. 
\end{lem}

\begin{proof}
We may apply the argument of Lemma \ref{lemVI.5.1} thanks to Proposition \ref{propVI.1.2} and Lemma \ref{lemVI.2.2}. 
\end{proof}


\section{\textbf{The invariant trace formula}}\label{secinvtf}

\subsection{The case of $(G,H)$}

For $f, f'\in\CC_c^\infty(\fs(F))$, we define
\begin{equation}\label{geominvltf1}
\begin{split}
 I^G(\eta, f, f'):=&\sum_{M\in\msl^{G, \omega}(M_0)} |W_0^{M_n}| |W_0^{GL_n}|^{-1} (-1)^{\dim(A_M/A_G)} \sum_{\fc\in\mst_\el(\fm\cap\fs)} |W(M_H, \fc)|^{-1} \int_{\fc_\reg(F)} \\
 &\hat{I}_M^G(\eta, X, f) I_G^G(\eta, X, f') dX. 
\end{split}
\end{equation}
From Proposition \ref{propwoi1}.2), for any $\fc\in\mst_\el(\fm\cap\fs)$, $I_G^G(\eta, \cdot, f')$ vanishes outside a compact subset of $\fc(F)$, so one may apply Lemma \ref{lemVI.3.1}.5) to show that this expression is absolutely convergent with the help of Proposition \ref{bdoi1} and Corollary \ref{cor20.2}. 

\begin{thm}[Invariant trace formula]\label{invltf1}
For all $f, f'\in\CC_c^\infty(\fs(F))$, we have the equality
$$ I^G(\eta, f, f')=I^G(\eta, f', f). $$
\end{thm}

The rest of this section is devoted to the proof of Theorem \ref{invltf1}. We shall follow the main steps in \cite[\S VII.2-3]{MR1344131}. The theorem will be proved by induction on the dimension of $G$. 

Let $f\in\CC_c^\infty(\fs_\rs(F))$ and $M\in\msl^{G,\omega}(M_0)$. By Proposition \ref{propwoi1}.2) and 3), the function $\kappa(\cdot)J_M^G(\eta, \cdot, f): (\fm\cap\fs_\rs)(F)\ra\BC$ is locally constant and invariant by the adjoint action of $M_H(F)$, where $\kappa$ is defined by \eqref{eqdeftf}. Moreover, the support of its restriction to $\fc(F)$ for any $\fc\in\mst(\fm\cap\fs)$ is included in the compact subset $\fc_\reg(F)\cap\Ad(H(F))(\Supp(f))$. Then $\kappa(\cdot)J_M^G(\eta, \cdot, f)$ defines a locally constant function on $(\fm\cap\fs)_\rs(F)$ (see Section \ref{sectsympar1} for the notation) via the extension by zero on the complement. From Harish-Chandra's submersion principle (Lemma \ref{thm11}), there exists $f'\in\CC_c^\infty((\fm\cap\fs)_\rs(F))$ such that
$$ \kappa(X)J_M^G(\eta, X, f)=|D^{\fm\cap\fs}(X)|_F^{1/2}\int_{M_{H,X}(F)\bs M_H(F)}f'(\Ad(x^{-1})(X))dx $$
for all $X\in(\fm\cap\fs_\rs)(F)$. Let $f'':=\kappa f'\in\CC_c^\infty((\fm\cap\fs)_\rs(F))$, where we extend the definition of $\kappa$ to the product form. Then we have
\begin{equation}\label{defphiMGf}
 J_M^G(\eta, X, f)=J_M^M(\eta, X, f''), \forall X\in(\fm\cap\fs_\rs)(F). 
\end{equation}
We have shown that for $f\in\CC_c^\infty(\fs_\rs(F))$ and $M\in\msl^{G,\omega}(M_0)$, there exists a function $f''\in\CC_c^\infty((\fm\cap\fs)_\rs(F))$ such that \eqref{defphiMGf} holds. We shall fix such an $f''$ and denote it by $\phi_M^G(f)$. 

As before, one may extend in the obvious way the definition (\ref{geominvltf1}) and the notation $\phi_M^G(f)$ to the symmetric pair $(M, M_H, \Ad(\epsilon))$, where $M\in\msl^{G,\omega}(M_0)$, since it appears as the product of some copies of the form $(G, H, \Ad(\epsilon))$ in lower dimensions. 

\begin{lem}\label{lemVII.2}
Let $M\in\msl^{G,\omega}(M_0), X\in(\fm\cap\fs_\rs)(F)_\el$ and $f,f'\in\CC_c^\infty(\fs_\rs(F))$. Then we have the equality
\begin{equation}\label{equlemVII.2}
 J_M^G(\eta, X, \hat{f}, f')=\sum_{L\in\msl^G(M)} \sum_{L_1,L_2\in\msl^G(L)} d_L^G(L_1,L_2) \hat{I}_M^L(\eta, X, \phi_L^{L_1}(f_{\ov{Q_1}}^\eta)) I_L^L(\eta, X,\phi_L^{L_2}({f'}_{Q_2}^\eta)), 
\end{equation}
where $J_M^G(\eta, X, \hat{f}, f')$ is defined by (\ref{defV.1.1}), and $(Q_1,Q_2):=s(L_1,L_2)$ (see Section \ref{gmmaps}). 
\end{lem}

\begin{proof}
By definition, 
$$ J_M^G(\eta, X, \hat{f}, f')=|D^\fs(X)|_F^{1/2} \int_{A_M(F)\bs H(F)} f'(\Ad(y^{-1})(X)) \eta(\Nrd(y)) \varphi_1(y) dy, $$
where
$$ \varphi_1(y):=|D^\fs(X)|_F^{1/2} \int_{A_M(F)\bs H(F)} \hat{f}(\Ad(x^{-1})(X)) \eta(\Nrd(x^{-1})) v_M(x, y) dx. $$
For $L\in\msl^G(M)$ and $L_2\in\msl^G(L)$, since $Q_2\in\msp^G(L_2)$, by Proposition \ref{propwoi1}.4), we have 
\begin{equation}\label{equ0pfVII.2}
 I_L^L(\eta, X,\phi_L^{L_2}({f'}_{Q_2}^\eta))=J_L^{L_2}(\eta, X, {f'}_{Q_2}^\eta)=J_L^{Q_2}(\eta, X, f'). 
\end{equation}
Since $X\in(\fm\cap\fs_\rs)(F)_\el$, the right hand side of (\ref{equlemVII.2}) is
$$ |D^\fs(X)|_F^{1/2} \int_{A_M(F)\bs H(F)} f'(\Ad(y^{-1})(X)) \eta(\Nrd(y)) \varphi_2(y) dy, $$
where
$$ \varphi_2(y):=\sum_{L\in\msl^G(M)} \sum_{L_1,L_2\in\msl^G(L)} d_L^G(L_1,L_2) \hat{I}_M^L(\eta, X, \phi_L^{L_1}(f_{\ov{Q_1}}^\eta)) v_L^{Q_2}(y). $$
It suffices to fix $y\in H(F)$ and prove that $\varphi_1(y)=\varphi_2(y)$. 

Let $L\in\msl^G(M)$ and
$$ h_L:=\sum_{L_1,L_2\in\msl^G(L)} d_L^G(L_1,L_2) \phi_L^{L_1}(f_{\ov{Q_1}}^\eta) v_L^{Q_2}(y). $$
Then
$$ \varphi_2(y)=\sum_{L\in\msl^G(M)} \hat{I}_M^L(\eta, X, h_L). $$
For $Y\in(\fl\cap\fs_\rs)(F)$, we have
$$ J_L^L(\eta, Y, h_L)=\sum_{L_1,L_2\in\msl^G(L)} d_L^G(L_1,L_2) J_L^L(\eta, Y, \phi_L^{L_1}(f_{\ov{Q_1}}^\eta)) v_L^{Q_2}(y). $$
For $L_1\in\msl^G(L)$, as in (\ref{equ0pfVII.2}), we have
$$ J_L^L(\eta, Y, \phi_L^{L_1}(f_{\ov{Q_1}}^\eta))=J_L^{L_1}(\eta, Y, f_{\ov{Q_1}}^\eta)=J_L^{\ov{Q_1}}(\eta, Y, f). $$
Then
$$ J_L^L(\eta, Y, h_L)=|D^\fs(Y)|_F^{1/2} \int_{H_Y(F)\bs H(F)} f(\Ad(x^{-1})(Y)) \eta(\Nrd(x)) h(x, y) dx, $$
where
$$ h(x, y):=\sum_{L_1,L_2\in\msl^G(L)} d_L^G(L_1,L_2) v_L^{\ov{Q_1}}(x) v_L^{Q_2}(y). $$
It is shown in the proof of \cite[Lemme VII.2]{MR1344131} that
$$ h(x, y)=v_L(x,y)=\sum_{Q\in\msf^G(L)} v'_Q(y) v_L^{\ov{Q}}(x). $$
Thus 
$$ J_L^L(\eta, Y, h_L)=\sum_{Q\in\msf^G(L)} v'_Q(y) J_L^{\ov{Q}}(\eta, Y, f). $$
As in (\ref{equ0pfVII.2}), we have
$$ J_L^{\ov{Q}}(\eta, Y, f)=J_L^{M_Q}(\eta, Y, f_{\ov{Q}}^\eta)=J_L^L(\eta, Y, \phi_L^{M_Q}(f_{\ov{Q}}^\eta)). $$
Let
$$ h'_L:=\sum_{Q\in\msf^G(L)} v'_Q(y) \phi_L^{M_Q}(f_{\ov{Q}}^\eta). $$
Then we obtain
\begin{equation}\label{equh_LpfVII.2}
 J_L^L(\eta, Y, h_L)=J_L^L(\eta, Y, h'_L) 
\end{equation}
for all $Y\in(\fl\cap\fs_\rs)(F)$. 

By a product form of Proposition \ref{propVI.1.1}, the distribution $d:=\hat{I}_M^L(\eta, X, \cdot)\in\msd^{\fl\cap\fs}$ and is $\eta(\Nrd(\cdot))$-invariant with respect to the adjoint action of $L_H(F)$. By a product form of (\ref{equV.6inv1}), we deduce from (\ref{equh_LpfVII.2}) that
$$ \hat{I}_M^L(\eta, X, h_L)=\hat{I}_M^L(\eta, X, h'_L). $$
Therefore, 
\[\begin{split}
 \varphi_2(y)&=\sum_{L\in\msl^G(M)} \hat{I}_M^L(\eta, X, h_L)=\sum_{L\in\msl^G(M)} \hat{I}_M^L(\eta, X, h'_L) \\
&=\sum_{Q\in\msf^G(M)} v'_Q(y) \sum_{L\in\msl^{M_Q}(M)} \hat{I}_M^L(\eta, X, \phi_L^{M_Q}(f_{\ov{Q}}^\eta)). 
\end{split}\]
By (\ref{wtind1}) and (\ref{equV.6inv1}) (actually their product forms are needed), we have
\[\begin{split}
&\hat{I}_M^{L, M_Q, w}(\eta, X, f_{\ov{Q}}^\eta)=\Ind_L^{M_Q,w}(d)(f_{\ov{Q}}^\eta) \\
=&\sum_{\{L'\in\msl^{G,\omega}(M_0): L'\subseteq L\}} |W_0^{L'_n}| |W_0^{L_n}|^{-1} \sum_{\fc\in\mst_\el(\fl'\cap\fs)} |W(L'_H, \fc)|^{-1} \int_{\fc_\reg(F)} J_L^{M_Q}(\eta, Z, f_{\ov{Q}}^\eta) e_d(Z) dZ \\
=&\sum_{\{L'\in\msl^{G,\omega}(M_0): L'\subseteq L\}} |W_0^{L'_n}| |W_0^{L_n}|^{-1} \sum_{\fc\in\mst_\el(\fl'\cap\fs)} |W(L'_H, \fc)|^{-1} \int_{\fc_\reg(F)} J_L^L(\eta, Z, \phi_L^{M_Q}(f_{\ov{Q}}^\eta)) e_d(Z) dZ \\
=&d(\phi_L^{M_Q}(f_{\ov{Q}}^\eta))=\hat{I}_M^L(\eta, X, \phi_L^{M_Q}(f_{\ov{Q}}^\eta)). 
\end{split}\]
Then by (\ref{invftofwoi1}), we get
\[\begin{split}
 &\sum_{L\in\msl^{M_Q}(M)} \hat{I}_M^L(\eta, X, \phi_L^{M_Q}(f_{\ov{Q}}^\eta)) =\sum_{L\in\msl^{M_Q}(M)} \hat{I}_M^{L, M_Q, w}(\eta, X, f_{\ov{Q}}^\eta) \\ 
=&\hat{J}_M^{M_Q}(\eta, X, f_{\ov{Q}}^\eta)=J_M^{M_Q}(\eta, X, \hat{f}_{\ov{Q}}^\eta)=J_M^{\ov{Q}}(\eta, X, \hat{f}). 
\end{split}\]
Hence, 
\[\begin{split}
 \varphi_2(y)=&\sum_{Q\in\msf^G(M)} v'_Q(y) J_M^{\ov{Q}}(\eta, X, \hat{f}) \\
=&|D^\fs(X)|_F^{1/2} \int_{A_M(F)\bs H(F)} \hat{f}(\Ad(x^{-1})(X)) \eta(\Nrd(x)) \sum_{Q\in\msf^G(M)} v'_Q(y) v_M^{\ov{Q}}(x) dx. 
\end{split}\]
But
$$ \sum_{Q\in\msf^G(M)} v'_Q(y) v_M^{\ov{Q}}(x)=v_M(x,y), $$
which implies that $\varphi_1(y)=\varphi_2(y)$. 
\end{proof}

\begin{proof}[Proof of Theorem \ref{invltf1}]
We use induction on the dimension of $G$. Suppose that the equality is true for $L\in\msl^{G,\omega}(M_0), L\neq G$, which is actually a product form in lower dimensions. Now we would like to prove the equality for $G$. The argument below is also valid for the case $\msl^{G,\omega}(M_0)=\{G\}$, i.e., $n=1$. 

First of all, suppose that $f, f'\in\CC_c^\infty(\fs_\rs(F))$. Applying Lemma \ref{lemVII.2} to the definition (\ref{geomnoninvltf1}) of $J^G(\eta, \hat{f}, f')$, we obtain
\[\begin{split}
 J^G(\eta, \hat{f}, f')=&\sum_{M\in\msl^{G,\omega}(M_0)} |W_0^{M_n}| |W_0^{GL_n}|^{-1} (-1)^{\dim(A_M/A_G)} \sum_{\fc\in\mst_\el(\fm\cap\fs)} |W(M_H, \fc)|^{-1} \int_{\fc_\reg(F)} \\
&J_M^G(\eta, X, \hat{f}, f') dX \\
=&\sum_{M\in\msl^{G,\omega}(M_0)} |W_0^{M_n}| |W_0^{GL_n}|^{-1} (-1)^{\dim(A_M/A_G)} \sum_{\fc\in\mst_\el(\fm\cap\fs)} |W(M_H, \fc)|^{-1} \int_{\fc_\reg(F)} \\
&\sum_{L\in\msl^G(M)} \sum_{L_1,L_2\in\msl^G(L)} d_L^G(L_1,L_2) \hat{I}_M^L(\eta, X, \phi_L^{L_1}(f_{\ov{Q_1}}^\eta)) I_L^L(\eta, X,\phi_L^{L_2}({f'}_{Q_2}^\eta)) dX \\
=&\sum_{L\in\msl^{G,\omega}(M_0)} |W_0^{L_n}| |W_0^{GL_n}|^{-1} (-1)^{\dim(A_L/A_G)} B_L(\eta, f, f'), 
\end{split}\]
where
\[\begin{split}
B_L(\eta, f, f'):=&\sum_{L_1,L_2\in\msl^G(L)} d_L^G(L_1,L_2) \sum_{\{M\in\msl^{G,\omega}(M_0): M\subseteq L\}} |W_0^{M_n}| |W_0^{L_n}|^{-1} (-1)^{\dim(A_M/A_L)} \sum_{\fc\in\mst_\el(\fm\cap\fs)} \\
&|W(M_H, \fc)|^{-1} \int_{\fc_\reg(F)} \hat{I}_M^L(\eta, X, \phi_L^{L_1}(f_{\ov{Q_1}}^\eta)) I_L^L(\eta, X,\phi_L^{L_2}({f'}_{Q_2}^\eta)) dX \\
=&\sum_{L_1,L_2\in\msl^G(L)} d_L^G(L_1,L_2) I^L(\eta, \phi_L^{L_1}(f_{\ov{Q_1}}^\eta), \phi_L^{L_2}({f'}_{Q_2}^\eta)). 
\end{split}\]
Here we have used the absolute convergence of the expressions above to exchange the order of sums, and $I^L(\eta, \phi_L^{L_1}(f_{\ov{Q_1}}^\eta), \phi_L^{L_2}({f'}_{Q_2}^\eta))$ is defined by a product form of (\ref{geominvltf1}). 

By the noninvariant trace formula (Theorem \ref{noninvltf1}) and Remark \ref{rmkV.1.1}, we have the equality $J^G(\eta, \hat{f}, f')=J^G(\eta, \hat{f'}, f)$. Therefore, 
\begin{equation}\label{equVII.3(1)}
 \sum_{L\in\msl^{G,\omega}(M_0)} |W_0^{L_n}| |W_0^{GL_n}|^{-1} (-1)^{\dim(A_L/A_G)} (B_L(\eta, f, f')-B_L(\eta, f', f))=0. 
\end{equation}
Let $L\in\msl^{G,\omega}(M_0), L\neq G$. Applying the induction hypothesis, we have
$$ B_L(\eta, f, f')=\sum_{L_1,L_2\in\msl^G(L)} d_L^G(L_1,L_2) I^L(\eta, \phi_L^{L_2}({f'}_{Q_2}^\eta), \phi_L^{L_1}(f_{\ov{Q_1}}^\eta)). $$
By exchanging $L_1$ and $L_2$ and by using (2) and (5) in Section \ref{gmmaps}, we obtain $B_L(\eta, f, f')=B_L(\eta, f', f)$. We deduce from (\ref{equVII.3(1)}) that $B_G(\eta, f, f')=B_G(\eta, f', f)$. But
$$ B_G(\eta, f, f')=I^G(\eta, f, f'), $$
which implies $I^G(\eta, f, f')=I^G(\eta, f', f)$. 

Now consider $f, f'\in\CC_c^\infty(\fs(F))$ in general. Let $\{\Omega_i\}_{i\geq1}$ be a sequence of increasing open compact subsets of $\fs_\rs(F)$ such that $\bigcup\limits_{i\geq1}^\infty \Omega_i=\fs_\rs(F)$. Such a consequence exists. For example, one may take $\Omega_i:=\{X\in\fs_\rs(F): \|X\|\leq i\}$ for all $i\geq 1$, where $\|\cdot\|$ denotes the abstract norm on $\fs_\rs(F)$ defined by \cite[(18.2.1) in \S18.2]{MR2192014}. From \cite[Proposition 18.1.(3)]{MR2192014}, since $\|\cdot\|$ is continuous, we deduce that $\Omega_i$ is compact for all $i\geq1$. It is obvious that $\Omega_i$ is open for all $i\geq1$ and that $\bigcup\limits_{i\geq1}^\infty \Omega_i=\fs_\rs(F)$. For all $i\geq1$, denote by $1_{\Omega_i}$ the characteristic function of $\Omega_i$. Let $f_i:=f 1_{\Omega_i}$ and $f'_i:=f' 1_{\Omega_i}$. 

Let $M\in\msl^{G,\omega}(M_0)$ and $\fc\in\msl_\el(\fm\cap\fs)$. For all $X\in\fc_\reg(F)$, by Lebesgue's dorminated convergence theorem, we have $\lim\limits_{i\ra\infty} I_G^G(\eta, X, f'_i)=I_G^G(\eta, X, f')$. For $X\in(\fm\cap\fs_\rs)(F)$, because $\hat{I}_M^G(\eta, X, \cdot)\in\msd^\fs$ (see Proposition \ref{propVI.1.1}), again by Lebesgue's dorminated convergence theorem, we have $\lim\limits_{i\ra\infty} \hat{I}_M^G(\eta, X, f_i)=\hat{I}_M^G(\eta, X, f)$. Because of Lemma \ref{lem28} applied to $\Supp(f')$, there exists a compact subset $r\subseteq\fc(F)$ such that for all $X\in\fc_\reg(F)-r, I_G^G(\eta, X, f'_i)=0$ for all $i\geq1$. By Lemma \ref{lemVI.3.1}.4) applied to $r$ and $\Supp(f')$, there exist constants $c>0$ and $N\in\BN$ such that
\[\begin{split}
 |\hat{I}_M^G(\eta, X, f_i)|&= \int_{\fs(F)} |f_i(Y) \hat{i}_M^G(\eta, X, Y)| |D^\fs(Y)|_F^{-1/2} dY \\
&\leq c\sup\{1, -\log|D^\fs(X)|_F\}^N \int_{\fs(F)} |f(Y)| \sup\{1, -\log|D^\fs(Y)|_F\}^N |D^\fs(Y)|_F^{-1/2} dY 
\end{split}\]
for all $i\geq1$ and $X\in r\cap\fc_\reg$. For all $X\in\fs_\rs(F)$, we also have $|I_G^G(\eta, X, f'_i)|\leq I_X(|f'|)$, where $I_X$ is defined by (\ref{oiwochar}). Combining Corollary \ref{corII.1}, Corollary \ref{cor20.2} and Proposition \ref{bdoi1}, we deduce that $\{\hat{I}_M^G(\eta, X, f_i) I_G^G(\eta, X, f'_i)\}_{i\geq1}$ is bounded by an integrable function on $\fc_\reg(F)$. Using Lebesgue's dorminated convergence theorem once again, we obtain
$$ \int_{\fc_\reg(F)} \hat{I}_M^G(\eta, X, f) I_G^G(\eta, X, f') dX=\lim_{i\ra\infty} \int_{\fc_\reg(F)} \hat{I}_M^G(\eta, X, f_i) I_G^G(\eta, X, f'_i) dX. $$
Therefore, 
$$ I^G(\eta, f, f')=\lim_{i\ra\infty} I^G(\eta, f_i, f'_i). $$
By exchanging $f$ and $f'$ and using the regular semi-simple support case that we have proved, we draw the conclusion. 
\end{proof}

\begin{coro}\label{corinvltf1}
Let $M, L\in\msl^{G, \omega}(M_0)$, $X\in(\fm\cap\fs_\rs)(F)_\el$ and $Y\in(\fl\cap\fs_\rs)(F)_\el$. Then we have the equality
$$ (-1)^{\dim(A_M/A_G)} \hat{i}_M^G(\eta, X, Y)=(-1)^{\dim(A_L/A_G)} \hat{i}_L^G(\eta, Y, X). $$
\end{coro}

\begin{proof}
By Lemma \ref{lemVI.3.1}.2), up to $M_H(F)$-conjugation on $X$ and $L_H(F)$-conjugation on $Y$, we may and shall suppose that there exist Cartan subspaces $\fc_1\in\mst_\el(\fm\cap\fs)$ and $\fc_2\in\mst_\el(\fl\cap\fs)$ such that $X\in\fc_{1, \reg}(F)$ and $Y\in\fc_{2, \reg}(F)$. As in the proof of Proposition \ref{repr1}, we can choose an open compact neighbourhood $V_1$ of $X$ in $\fc_{1, \reg}$ (resp. $V_2$ of $Y$ in $\fc_{2, \reg}$) such that if two elements in $V_1$ (resp. $V_2$) are $H(F)$-conjugate, then they are the same. Let $f, f'\in\CC_c^\infty(\fs(F))$ with $\Supp(f)\subseteq \Ad(H(F))(V_2)$ and $\Supp(f')\subseteq \Ad(H(F))(V_1)$. By an analogous calculation to that of $J^G(\eta, \hat{f}, f')$ in the proof of Proposition \ref{repr1}, with the help of Lemma \ref{lemVI.3.1}.2) and Proposition \ref{propwoi1}.3), we show the equalities
$$ I^G(\eta, f, f')=(-1)^{\dim(A_M/A_G)}\int_{V_1} \hat{I}_M^G(\eta, X_1, f) I_G^G(\eta, X_1, f') dX_1 $$
and
$$ \hat{I}_M^G(\eta, X_1, f)=\int_{V_2} \hat{i}_M^G(\eta, X_1, Y_2) I_G^G(\eta, Y_2, f) dY_2 $$
for all $X_1\in V_1$ by (\ref{equV.6inv1}). Then
$$ I^G(\eta, f, f')=(-1)^{\dim(A_M/A_G)}\int_{V_1\times V_2} \hat{i}_M^G(\eta, X_1, Y_2) I_G^G(\eta, Y_2, f) I_G^G(\eta, X_1, f') dY_2 dX_1. $$
Similarly, we have
$$ I^G(\eta, f', f)=(-1)^{\dim(A_L/A_G)}\int_{V_2\times V_1} \hat{i}_L^G(\eta, Y_2, X_1) I_G^G(\eta, X_1, f') I_G^G(\eta, Y_2, f) dX_1 dY_2. $$
By Harish-Chandra's submersion principle (Lemma \ref{thm11}), when $f'$ varies, the function $X_1\mapsto I_{X_1}(f')=\kappa(X_1) I_G^G(\eta, X_1, \kappa f')$ on $V_1$ runs over all $\CC_c^\infty(V_1)$, so the function $I_G^G(\eta, \cdot, f')$ on $V_1$ also runs over all $\CC_c^\infty(V_1)$. Similarly, when $f$ varies, the function $I_G^G(\eta, \cdot, f)$ on $V_2$ runs over all $\CC_c^\infty(V_2)$. Then from the invariant trace formula (Theorem \ref{invltf1}), we deduce that
$$ (-1)^{\dim(A_M/A_G)} \hat{i}_M^G(\eta, X_1, Y_2)=(-1)^{\dim(A_L/A_G)} \hat{i}_L^G(\eta, Y_2, X_1) $$
for all $(X_1, Y_2)\in V_1\times V_2$. We conclude by $(X, Y)\in V_1\times V_2$. 
\end{proof}

\subsection{The case of $(G',H')$}

For $f, f'\in\CC_c^\infty(\fs'(F))$, we define
\begin{equation}\label{geominvltf2}
\begin{split}
 I^{H'}(f, f'):=&\sum_{M'\in\msl^{H'}(M'_0)} |W_0^{H'}||W_0^{M'}|^{-1} (-1)^{\dim(A_{M'}/A_{H'})} \sum_{\fc'\in\mst_\el(\wt{\fm'}\cap\fs')} |W(M', \fc')|^{-1} \int_{\fc'_\reg(F)} \\
 &\hat{I}_{M'}^{H'}(Y, f) I_{H'}^{H'}(Y, f') dY. 
\end{split}
\end{equation}
From Proposition \ref{propwoi2}.2), for any $\fc'\in\mst_\el(\wt{\fm'}\cap\fs')$, $I_{H'}^{H'}(\cdot, f')$ vanishes outside a compact subset of $\fc'(F)$, so one may apply Lemma \ref{lemVI.3.2}.5) to show that this expression is absolutely convergent with the help of Proposition \ref{bdoi2} and Corollary \ref{cor20.2}. One may extend in the obious way the definition \eqref{geominvltf2} to the symmetric pair $(\wt{M'}, M', \Ad(\alpha))$, where $M'\in\msl^{H'}(M'_0)$. 

\begin{thm}[Invariant trace formula]\label{invltf2}
For all $f, f'\in\CC_c^\infty(\fs'(F))$, we have the equality
$$ I^{H'}(f, f')=I^{H'}(f', f). $$
\end{thm}

\begin{proof}
We may apply the argument of Theorem \ref{invltf1} with obvious modifications. It is deduced from the noninvariant trace formula (Theorem \ref{noninvltf2}) and other results that we have prepared in previous sections. 
\end{proof}

\begin{coro}\label{corinvltf2}
Let $M', L'\in\msl^{H'}(M'_0)$, $Y\in(\wt{\fm'}\cap\fs'_\rs)(F)_\el$ and $X\in(\wt{\fl'}\cap\fs'_\rs)(F)_\el$. Then we have the equality
$$ (-1)^{\dim(A_{M'}/A_{H'})} \hat{i}_{M'}^{H'}(Y, X)=(-1)^{\dim(A_{L'}/A_{H'})} \hat{i}_{L'}^{H'}(Y, X). $$
\end{coro}

\begin{proof}
We may apply the argument of Corollary \ref{corinvltf1} by using the invariant trace formula (Theorem \ref{invltf2}) and consulting the proof of Proposition \ref{repr2}. 
\end{proof}


\section{\textbf{A vanishing property at infinity}}\label{secvanpro}

\subsection{The case of $(G,H)$}

The following proposition is an analogue of \cite[Proposition 2.2]{MR2164623}. 

\begin{prop}\label{limitformula1}
Let $M\in\msl^{G, \omega}(M_0), M\neq G$. Let $X\in(\fm\cap\fs_\rs)(F)$ and $Y\in\fs_\rs(F)$. Then there exists $N\in\BN$ such that if $\lambda\in F^\times$ satisfies $v_F(\lambda)<-N$, we have
$$ \hat{i}_M^G(\eta, \lambda X, Y)=0. $$
\end{prop}

\begin{remark}
A limit formula at infinity for $\hat{i}_G^G(\eta, \lambda X, Y)$ in the spirit of Laplace transform is given in \cite[Proposition 7.1]{MR3414387} (see also \cite[Proposition 6.4]{MR3299843}), which is an analogue of \cite[Proposition VIII.1]{MR1344131}. 
\end{remark}

\begin{proof}[Proof of Proposition \ref{limitformula1}]
We shall imitate the proof of \cite[Proposition 2.2]{MR2164623}. 

By Lemma \ref{lemVI.3.1}.2), up to $H(F)$-conjugation on $Y$, we may and shall suppose that there exists $L\in\msl^{G, \omega}(M_0)$ and an $L$-elliptic Cartan subspace $\fc\subseteq\fl\cap\fs$ such that $Y\in\fc_\reg(F)$. By Lemma \ref{lemVI.2.1}, we have the equality
$$ \hat{i}_M^G(\eta, \lambda X, Y)=\hat{j}_M^G(\eta, \lambda X, Y). $$
Thus it suffices to prove that there exists $N\in\BN$ such that if $\lambda\in F^\times$ satisfies $v_F(\lambda)<-N$, we have
$$ \hat{j}_M^G(\eta, \lambda X, Y)=0. $$

Fix an $\CO_F$-lattice $k_\fh$ (resp. $k_\fs$) of $\fh(F)$ (resp. $\fs(F)$). Denote by $\wt{k_\fs}$ the dual $\CO_F$-lattice of $k_\fs$ in $\fs(F)$, i.e., 
$$ \wt{k_\fs}:=\{Z\in\fs(F): \forall Z'\in k_\fs, \Psi(\langle Z,Z'\rangle)=1\}. $$
Set
$$ \fc(X):=\{X'\in\fc(F): \exists x\in H(F), \Ad(x)(X')=X\}, $$
which is a finite (perhaps empty) set. For $\lambda\in F^\times$, choose $h_\lambda\in\BN$ such that both of the functions $\hat{j}_M^G(\eta, \lambda X, \cdot)$ and $|D^\fs(\cdot)|_F$ are constant on $Y+\varpi^{h_\lambda} k_\fs$. 

Let $f$ (resp. $f'$) $\in\CC_c^\infty(\fs(F))$ be the characteristic function of $Y+\varpi^{h_\lambda}k_\fs$ (resp. $\varpi^{-h_\lambda}\wt{k_\fs}$). Then for $Z\in\fs(F)$, we see that
\[\begin{split}
\hat{f}(Z)=&c_\Psi(\fs(F))\int_{Y+\varpi^{h_\lambda}k_\fs} \Psi(\langle Z,Z' \rangle)dZ'=c_\Psi(\fs(F))\Psi(\langle Z,Y \rangle)\int_{\varpi^{h_\lambda}k_\fs} \Psi(\langle Z,Z' \rangle)dZ' \\
=&c_\Psi(\fs(F)) \vol(\varpi^{h_\lambda}k_\fs) \Psi(\langle Z,Y \rangle) f'(Z). 
\end{split}\]
Now there are two expressions for $J_M^G(\eta, \lambda X, \hat{f})$. One the one hand, 
\begin{equation}\label{prop2.2(2.5)exp1}
 J_M^G(\eta, \lambda X, \hat{f})=\int_{\fs(F)} f(Z) \hat{j}_M^G(\eta, \lambda X, Z) |D^\fs(Z)|_F^{-1/2} dZ=\vol(\varpi^{h_\lambda}k_\fs) \hat{j}_M^G(\eta, \lambda X, Y) |D^\fs(Y)|_F^{-1/2}. 
\end{equation}
On the other hand, 
\begin{equation}\label{prop2.2(2.6)exp2}
\begin{split}
 &J_M^G(\eta, \lambda X, \hat{f})=|D^\fs(\lambda X)|_F^{1/2} \int_{H_{\lambda X}(F)\bs H(F)} \hat{f}(\Ad(x^{-1})(\lambda X)) \eta(\Nrd(x)) v_M^G(x) dx \\
=&c_\Psi(\fs(F)) \vol(\varpi^{h_\lambda}k_\fs) |D^\fs(\lambda X)|_F^{1/2} \int_{H_X(F)\bs H(F)} \Psi(\langle \Ad(x^{-1})(\lambda X),Y \rangle) f'(\Ad(x^{-1})(\lambda X)) \eta(\Nrd(x)) v_M^G(x) dx. 
\end{split}
\end{equation}

Fix an open neighbourhood $\msv_\fh$ of $0$ in $\fh(F)$ which is invariant by the adjoint action of $H(F)$ such that a homeomorphic exponential map is defined on $\msv_\fh$. Choose $a\in\BN$ verifying the following conditions: 

(1) $\varpi^a k_\fh\subseteq\msv_\fh$; 

(2) $K_a:=\exp(\varpi^a k_\fh)$ is a subgroup of $K_H$; 

(3) $\eta(\Nrd(K_a))=1$; 

(4) the adjoint action of $K_a$ stabilises $k_\fs$ (and thus $\wt{k_\fs}$). 

Fix a set $\Gamma$ of representatives in $H(F)$ of double cosets $H_X(F)\bs H(F)/K_a$. 
We may and shall suppose that if $x\in\Gamma$ and $y\in H_X(F)xK_a$ verify $\Ad(y^{-1})(X)\in\fc(F)$, then $\Ad(x^{-1})(X)\in\fc(F)$. 

The integral in (\ref{prop2.2(2.6)exp2}) can be decomposed as
$$ \sum_{x\in\Gamma} \int_{H_X(F)\bs H_X(F)xK_a} \Psi(\langle \Ad(y^{-1})(\lambda X),Y \rangle) f'(\Ad(y^{-1})(\lambda X)) \eta(\Nrd(y)) v_M^G(y) dy. $$
By the conditions (2), (3) and (4) on $a$ respectively, the factors $v_M^G$, $\eta$ and $f'$ can be extracted from the integral. By comparing (\ref{prop2.2(2.5)exp1}) and (\ref{prop2.2(2.6)exp2}), since $\langle\cdot,\cdot\rangle$ is invariant by the adjoint action of $K_a$, we obtain
\begin{equation}\label{prop2.2(2.7)}
\begin{split}
\hat{j}_M^G(\eta, \lambda X, Y)=&c_\Psi(\fs(F)) |D^\fs(\lambda X)D^\fs(Y)|_F^{1/2} \sum_{x\in\Gamma} f'(\Ad(x^{-1})(\lambda X)) \eta(\Nrd(x)) v_M^G(x) \int_{H_X(F)\bs H_X(F)xK_a} \\
&\Psi(\langle \Ad(y^{-1})(\lambda X),Y \rangle) dy \\
=&c_\Psi(\fs(F)) |D^\fs(\lambda X)D^\fs(Y)|_F^{1/2} \vol(K_a)^{-1} \sum_{x\in\Gamma} \vol(H_X(F)\bs H_X(F)xK_a) f'(\Ad(x^{-1})(\lambda X)) \\ 
&\eta(\Nrd(x)) v_M^G(x) i(x), 
\end{split}
\end{equation}
where
$$ i(x):=\int_{K_a} \Psi(\langle \Ad(x^{-1})(\lambda X),\Ad(y)(Y) \rangle) dy. $$

For $x\in\Gamma$, consider the map $K_a\ra F$ defined by
\begin{equation}\label{equprop2.2(2.8)}
 \forall y\in K_a, y\mapsto \langle \Ad(x^{-1})(X),\Ad(y)(Y) \rangle. 
\end{equation}
Its differential at the point $y_0\in K_a$ is the map $\fh(F)\ra F$ defined by
\begin{equation}\label{equprop2.2(2.9)}
 \forall Z\in \fh(F), Z\mapsto \langle \Ad(x^{-1})(X),\Ad(y_0)([Z,Y]) \rangle. 
\end{equation}
Since $\langle\cdot,\cdot\rangle$ is invariant by the adjoint action of $G(F)$, we see that
$$ \langle \Ad(x^{-1})(X),\Ad(y_0)([Z,Y]) \rangle=\langle [Y, \Ad(xy_0)^{-1}(X)] , Z \rangle. $$ 
Because the restriction of $\langle\cdot,\cdot\rangle$ to $\fh(F)$ is non-degenerate, the map (\ref{equprop2.2(2.9)}) is not surjective if and only if
$$ [Y, \Ad(xy_0)^{-1}(X)]=0. $$
Since $Y\in\fc_\reg(F)$, this condition is equivalent to
$$ \Ad(xy_0)^{-1}(X)\in\fc(F). $$
From our choice of $\Gamma$, as $y_0\in K_a$, it implies that
$$ \Ad(x^{-1})(X)\in\fc(F). $$

Let
$$ \Gamma':=\{x\in\Gamma: \Ad(x^{-1})(X)\in\fc(F)\}, $$
which is a finite (perhaps empty) subset of $\Gamma$. Then for $x\in\Gamma-\Gamma'$, the map (\ref{equprop2.2(2.8)}) is a submersion. Define
$$ \Omega:=\bigcup_{x\in\Gamma-\Gamma'} H_X(F)xK_a, $$
which is an open and closed subset of $H(F)$. 
Fix a basis of the $F$-linear space $\fs(F)$. For $Z\in\fs(F)$, define its norm $\|Z\|\in\BR_{\geq0}$ as the maximum of normalised absolute values of coefficients of $Z$ with respect to the fixed basis. For $Z\in\fs(F)-\{0\}$, define $\nu(Z)\in\BZ$ by $\|Z\|=|\varpi^{\nu(Z)}|_F$. 
Let $S_X$ be the closure of
$$ S_X^0:=\{\varpi^{-\nu(\Ad(y^{-1})(X))}\Ad(y^{-1})(X): y\in\Omega\} $$
in the unit sphere $S_\fs:=\{Z\in\fs(F): \|Z\|=1\}$. 
Then $S_X$ is compact. Recall that we denote by $\CN^\fs$ the set of nilpotent elements in $\fs(F)$. 

\begin{lem}\label{HClem7.4}
We have
$$ S_X-S_X^0\subseteq\CN^\fs-\{0\}. $$
\end{lem}

\begin{proof}[Proof of Lemma \ref{HClem7.4}]
Since $S_X\subseteq S_\fs$, it is obvious that $\{0\}\notin S_X$. Let $Z\in S_X$. There exists a sequence $\{y_i\}$ in $\Omega$ such that when $i\ra\infty$, 
$$ \varpi^{-\nu(\Ad(y_i^{-1})(X))}\Ad(y_i^{-1})(X)\ra Z. $$
We distinguish two cases. 

i) Suppose that the sequence $\{\|\Ad(y_i^{-1})(X)\|\}$ remains bounded. By Harish-Chandra's compactness lemma for symmetric spaces (Lemma \ref{lem25}), the projection of the sequence $\{y_i\}$ to $H_X(F)\bs H(F)$ is contained in a compact subset. By taking a subsequence, since the projection of $\Omega$ to $H_X(F)\bs H(F)$ is closed, we may assume that when $i\ra\infty$, $\Ad(y_i^{-1})(X)\ra \Ad(y^{-1})(X)$ with $y\in\Omega$. Thus $Z\in S_X^0$ in this case. 

ii) Suppose that the sequence $\{\|\Ad(y_i^{-1})(X)\|\}$ is unbounded. By taking a subsequence, we may assume that when $i\ra\infty$, $\|\Ad(y_i^{-1})(X)\|\ra+\infty$. The eigenvalues of $\ad(\Ad(y_i^{-1})(X))$ are the same as those of $\ad(X)$; here $\ad(\Ad(y_i^{-1})(X))$ and $\ad(X)$ are viewed as linear endomorphisms of $\fg$. Thus the eigenvalues of $\ad(\varpi^{-\nu(\Ad(y_i^{-1})(X))}\Ad(y_i^{-1})(X))$ tend to zero when $i\ra\infty$. Hence $\ad(Z)$ is nilpotent. We shall prove that $Z\in\CN^\fs$ in this case. 

Since $\fg$ is reductive, one has $\fg=\fz\oplus\fg_\der$, where $\fz$ denotes the centre of $\fg$ and $\fg_\der$ denotes the derived algebra of $\fg$, and $\fg_\der$ is semi-simple. Let $Z=Z_1+Z_2$ with $Z_1\in\fz(F)$ and $Z_2\in\fg_\der(F)$. Since $\ad(Z)$ is nilpotent as a linear endomorphism of $\fg$, we deduce that $\ad(Z_2)$ is nilpotent as a linear endomorphism of $\fg_\der$. As $\fg_\der$ is semi-simple, we obtain that $Z_2$ is a nilpotent element in $\fg$. Let $X=X_1+X_2$ with $X_1\in\fz(F)$ and $X_2\in\fg_\der(F)$. The projection of $\varpi^{-\nu(\Ad(y_i^{-1})(X))}\Ad(y_i^{-1})(X)$ to $\fz(F)$ is equal to $\varpi^{-\nu(\Ad(y_i^{-1})(X))}X_1$, which tends to zero when $i\ra\infty$. Thus $Z_1=0$, and $Z=Z_2$ is a nilpotent element in $\fg$. Hence $Z\in\CN^\fs$. 
\end{proof}

For $U\in\CN^\fs-\{0\}$, consider the map $K_a\ra F$ defined by
\begin{equation}\label{equprop2.2nilp}
 \forall y\in K_a, y\mapsto \langle U,\Ad(y)(Y) \rangle. 
\end{equation}
Its differential at the point $y_0\in K_a$ is the map $\fh(F)\ra F$ defined by
$$ \forall Z\in \fh(F), Z\mapsto \langle U,\Ad(y_0)([Z,Y]) \rangle=\langle [Y, \Ad(y_0^{-1})(U)] , Z \rangle $$
by the $G(F)$-invariance of $\langle \cdot,\cdot \rangle$. Since $Y\in\fc_\reg(F)$ and $\Ad(y_0^{-1})(U)\in\CN^\fs-\{0\}$, we have
$$ [Y, \Ad(y_0^{-1})(U)]\neq 0. $$
Then the map (\ref{equprop2.2nilp}) is a submersion by the non-degeneration of $\langle \cdot,\cdot \rangle$ on $\fh(F)$. 

Using Lemma \ref{HClem7.4} and combining our discussion on the maps (\ref{equprop2.2(2.8)}) and (\ref{equprop2.2nilp}), we deduce that there exists an open compact neighbourhood $\wt{S_X}$ of $S_X$ in $S_\fs$ such that the map $\varphi: K_a\times \wt{S_X}\ra F\times \wt{S_X}$ defined by
$$ \forall (y,Z)\in K_a\times \wt{S_X}, (y,Z)\mapsto (\langle Z,\Ad(y)(Y) \rangle, Z) $$
is a submersion. Since any submersion is open, the image of $\varphi$ (denoted by $\Im(\varphi)$) is an open compact subset of $F\times \wt{S_X}$. 
Then the map $\varphi$ induces a surjective submersion $\varphi': K_a\times \wt{S_X}\ra \Im(\varphi)$. 
Applying Harish-Chandra's submersion principle \cite[Theorem 11]{MR0414797} to $\varphi'$, there exists a function $\phi\in\CC_c^\infty(\Im(\varphi))$ such that for all $\Phi'\in\CC_c^\infty(\Im(\varphi))$, 
$$ \int_{K_a\times\wt{S_X}} \Phi'(\langle Z,\Ad(y)(Y) \rangle, Z) dZdy=\int_{\Im(\varphi)} \phi(t,Z)\Phi'(t, Z) dZdt. $$
Fix such a $\phi$. Denote by $\CC^\infty(F\times \wt{S_X})$ the space of locally constant, complexed-valued functions on $F\times \wt{S_X}$. For $\Phi\in\CC^\infty(F\times \wt{S_X})$, the restriction of $\Phi$ to $\Im(\varphi)$ belongs to $\CC_c^\infty(\Im(\varphi))$, so we obtain
$$ \int_{K_a\times\wt{S_X}} \Phi(\langle Z,\Ad(y)(Y) \rangle, Z) dZdy=\int_{F\times \wt{S_X}} \phi(t,Z)\Phi(t, Z) dZdt. $$
By taking $\Phi(t,Z):=\Psi(\mu t)\beta(Z)$ with $\mu\in F$ and $\beta\in\CC_c^\infty(\wt{S_X})$, we deduce that for all $Z\in\wt{S_X}$, 
$$ \int_{K_a} \Psi(\mu \langle Z,\Ad(y)(Y) \rangle) dy=\int_F \phi(t,Z)\Psi(\mu t) dt. $$
Since $\Im(\phi)$ is an open compact subset of $F\times \wt{S_X}$, we see that $\phi\in\CC_c^\infty(\Im(\varphi))\subseteq\CC_c^\infty(F\times \wt{S_X})=\CC_c^\infty(F)\otimes\CC_c^\infty(\wt{S_X})$. Suppose that $\phi=\sum\limits_{1\leq j\leq m} c_j\cdot\xi_j\otimes\chi_j$ with $c_j\in \BC, \xi_j\in\CC_c^\infty(F)$ and $\chi_j\in\CC_c^\infty(\wt{S_X})$. Then
$$ \int_F \phi(t,Z)\Psi(\mu t) dt=\sum_{1\leq j\leq m}c_j \hat{\xi_j}(\mu)\chi_j(Z), $$
where $\hat{\xi_j}\in\CC_c^\infty(F)$ is the Fourier transform of $\xi_j$. We see that there exists $N_0\in\BN$ such that for all $\mu\in F^\times$ satisfying $v_F(\mu)<-N_0$ and all $Z\in S_X$, we have
$$ \int_{K_a} \Psi(\mu \langle Z,\Ad(y)(Y) \rangle) dy=0. $$
Fix such an $N_0$. 

Recall that for $x\in\Gamma$, 
$$ i(x)=\int_{K_a} \Psi(\mu \langle Z,\Ad(y)(Y) \rangle) dy, $$
where $\mu:=\lambda\varpi^{\nu(\Ad(x^{-1})(X))}$ and $Z:=\varpi^{-\nu(\Ad(x^{-1})(X))}\Ad(x^{-1})(X)$. For $x\in\Gamma-\Gamma'$, we have $Z\in S_X$, so $i(x)=0$ if
$$ v_F(\lambda)+\nu(\Ad(x^{-1})(X))<-N_0. $$
Set
$$ \nu_0:=\sup_{x\in\Gamma} \nu(\Ad(x^{-1})(X)), $$
which is finite thanks to Harish-Chandra's compactness lemma for symmetric spaces (Lemma \ref{lem25}). Now let
$$ N:=N_0+\nu_0. $$
Suppose that $v_F(\lambda)<-N$. From (\ref{prop2.2(2.7)}), to show $\hat{j}_M^G(\eta, \lambda X, Y)=0$, it suffices to prove $v_M^G(x)=0$ for all $x\in\Gamma'$. 

For $x\in\Gamma'$, we have $\Ad(x^{-1})(X)\in\fc_\reg(F)$. Then $\Ad(x^{-1})(H_X)=T_\fc$. Since $X\in(\fm\cap\fs_\rs)(F)$, we see that $\Ad(x^{-1})(A_M)$ is an $F$-split torus in $T_\fc$. As $\fc\subseteq\fl\cap\fs$ is $L$-elliptic, $A_L$ is the maximal $F$-split torus in $T_\fc$. Thus $\Ad(x^{-1})(A_M)\subseteq A_L$. Then $A_M\subseteq \Ad(x)(A_L)\subseteq \Ad(x)(A_{M_0})$. We deduce that $\Ad(x)(A_{M_0})$ is a maximal $F$-split torus in $M_H$, so it is $M_H(F)$-conjugate to $A_{M_0}$. Therefore, $x\in M_H(F)\Norm_{H(F)}(M_0)\subseteq M_H(F)K_H$. Consequently, we have $v_M^G(x)=0$ and conclude. 
\end{proof}

\subsection{The case of $(G',H')$}

\begin{prop}\label{limitformula2}
Let $M'\in\msl^{H'}(M'_0), M'\neq H'$. Let $Y\in(\wt{\fm'}\cap\fs'_\rs)(F)$ and $X\in\fs'_\rs(F)$. Then there exists $N\in\BN$ such that if $\lambda\in F^\times$ satisfies $v_F(\lambda)<-N$, we have
$$ \hat{i}_{M'}^{H'}(\lambda Y, X)=0. $$
\end{prop}

\begin{proof}
It is almost the same as the proof of Proposition \ref{limitformula1}, except that one needs to use Lemma \ref{lemVI.3.2}.2) and Lemma \ref{lemVI.2.2}. 
\end{proof}


\bibliography{References}
\bibliographystyle{plain}

\medskip

\begin{flushleft}
Max-Planck-Institut für Mathematik, Vivatsgasse 7, 53111 Bonn, Germany \\
\medskip
E-mail: huajie.li@mpim-bonn.mpg.de \\
\end{flushleft}

\end{document}